\documentclass{amsart}
\usepackage{amscd,amsmath,amssymb,mathtools,combelow}

\makeatletter
\@addtoreset{equation}{section}
\makeatother

\newtheorem{theorem}[subsection]{Theorem}

\newtheorem{proposition}[subsection]{Proposition}
\newtheorem{lemma}[subsection]{Lemma}
\newtheorem{corollary}[subsection]{Corollary}
\newtheorem{conjecture}[subsection]{Conjecture}

\theoremstyle{definition}

\newtheorem{definition}[subsection]{Definition}

\theoremstyle{remark}

\newtheorem{claim}[subsection]{Claim}

\newtheorem{remark}[subsection]{Remark}

\def\loccitt{\emph{loc. cit.}}
\def\loccit{\emph{loc. cit. }}

\def\fsl{{\mathfrak{sl}}}
\def\fgl{{\mathfrak{gl}}}
\def\hsl{{\widehat{\fsl}}}
\def\hgl{{\widehat{\fgl}}}

\def\BA{{\mathbb{A}}}
\def\BC{{\mathbb{C}}}

\def\BN{{\mathbb{N}}}
\def\BF{{\mathbb{F}}}

\def\BR{{\mathbb{R}}}
\def\BQ{{\mathbb{Q}}}
\def\BZ{{\mathbb{Z}}}

\def\woo{\widehat{\otimes}}

\def\CA{{\mathcal{A}}}
\def\CB{{\mathcal{B}}}

\def\CN{{\mathcal{N}}}
\def\CO{{\mathcal{O}}}

\def\CR{{\mathcal{R}}}

\def\CT{{\mathcal{T}}}

\def\CV{{\mathcal{V}}}

\def\ph{\varphi}
\def\coop{{\textrm{coop}}}

\def\e{\varepsilon}

\def\col{\textrm{col }}

\def\vs{\varsigma}

\def\prim{\textrm{prim}}

\def\and{\textrm{ }\&\textrm{ }}

\def\sym{\textrm{Sym}}

\def\CC{{{\mathcal{C}}}}

\def\sym{\textrm{Sym}}

\def\nn{{{\BN}}^n}
\def\zz{{{{\mathbb{Z}}}^n}}

\def\su{{U_q(\dot{\fsl}_n)}}
\def\sug{{U_q^\geq(\dot{\fsl}_n)}}
\def\sul{{U_q^\leq(\dot{\fsl}_n)}}
\def\sup{{U_q^+(\dot{\fsl}_n)}}

\def\uui{{U_q(\dot{\fgl}_1)}}
\def\uuig{{U_q^\geq(\dot{\fgl}_1)}}
\def\uuil{{U_q^\leq(\dot{\fgl}_1)}}

\def\uu{{U_q(\dot{\fgl}_n)}}

\def\uux{{U_q(\dot{\fgl}_{\frac ng})}}

\def\uup{{U_q^+(\dot{\fgl}_n)}}
\def\uum{{U_q^-(\dot{\fgl}_n)}}
\def\uupm{{U_q^\pm(\dot{\fgl}_n)}}
\def\uul{{U_q^\leq(\dot{\fgl}_n)}}

\def\uug{{U_q^\geq(\dot{\fgl}_n)}}

\def\UU{{U_{q,\oq}(\ddot{\fgl}_n)}}
\def\UUo{{U^0_{q,\oq}(\ddot{\fgl}_n)}}
\def\UUm{{U^-_{q,\oq}(\ddot{\fgl}_n)}}
\def\UUp{{U^+_{q, \oq}(\ddot{\fgl}_n)}}
\def\UUl{{U^\leq_{q, \oq}(\ddot{\fgl}_n)}}
\def\UUg{{U^\geq_{q, \oq}(\ddot{\fgl}_n)}}

\def\A{{\CA}}
\def\Ap{{\CA^+}}
\def\Am{{\CA^-}}

\def\Ag{{\CA^\geq}}
\def\Ao{{\CA^0}}
\def\Al{{\CA^\leq}}

\def\Bp{\CB^+}

\def\Ups{{\Upsilon}}

\def\ba{{\mathbf{a}}}

\def\bk{{\mathbf{k}}}
\def\bl{{\mathbf{l}}}
\def\bs{{\boldsymbol{\vs}}}

\def\bde{{\boldsymbol{\delta}}}

\def\bideg{{\text{bideg }}}
\def\homdeg{{\text{hom deg }}}

\def\oq{{\overline{q}}}

\def\bari{\bar{i}}
\def\barj{\bar{j}}
\def\bark{\bar{k}}
\def\barii{\bar{i}'}
\def\barjj{\bar{j}'}

\begin{document}

\title[Quantum toroidal and shuffle algebras]{\Large{\textbf{Quantum toroidal and shuffle algebras}}}

\author[Andrei Negu\cb t]{Andrei Negu\cb t}
\address{MIT, Department of Mathematics, Cambridge, MA, USA}
\address{Simion Stoilow Institute of Mathematics, Bucharest, Romania}
\email{andrei.negut@@gmail.com}

\maketitle

\begin{abstract}

In this paper, we prove that the quantum toroidal algebra $\UU$ is isomorphic to the double shuffle algebra of Feigin and Odesskii for the cyclic quiver. The shuffle algebra viewpoint will allow us to prove a factorization formula for the universal $R-$matrix of the quantum toroidal algebra. \\

\end{abstract}

\section{Introduction}

\noindent The quantum toroidal algebra $\UU$ is defined, in the Drinfeld presentation\footnote{To avoid double hats above our symbols, we will denote hats by points in this paper}, by certain generators and relations  (see \cite{FJMM})\footnote{We set $q^c=1$ in the notation of \emph{loc. cit.}, which makes their algebra into a central extension of ours. The central element in question acts trivially in all representations coming from geometry}. It admits a triangular decomposition:
$$
\UU \ = \ \UUm \otimes U_{q,\oq}^0(\ddot{\fgl}_n) \otimes \UUp
$$
The shuffle algebra was defined by Feigin and Odesskii (\cite{FO}) as the space of certain symmetric elliptic functions. The case we consider in the present paper is a particular trigonometric degeneration $\CA^+$ of their algebra. Elements of $\CA^+$ are symmetric rational functions with prescribed poles, which furthermore satisfy the vanishing properties \eqref{eqn:wheel}, and with multiplication given by \eqref{eqn:mult}. There is a natural homomorphism between these two algebras, defined akin to the construction of Enriquez in the case of affine quantum groups (\cite{E}):
$$
\Upsilon^+ \ : \ \UUp \longrightarrow \Ap
$$
Feigin conjectured that the above map is an isomorphism, and one of the main goals of the present paper is to prove this. The difficult part is surjectivity, and proving this boils down to the fact that the shuffle algebra is generated by degree 1 elements, which we establish in Proposition \ref{prop:surj}. We study the Drinfeld double $\A$ of the shuffle algebra $\Ap$ (with respect to a coproduct that will be introduced in Proposition \ref{prop:coproduct} and the bialgebra pairing defined in Proposition \ref{prop:pairshuf}) and show that the required isomorphism extends to the Drinfeld doubles. \\

\begin{theorem}
\label{thm:iso}

There exists a bialgebra isomorphism $\Upsilon:\UU \stackrel{\sim}\longrightarrow \CA$. 

\end{theorem}

\text{} \\
We prove the above theorem by constructing a certain \underline{slope filtration} of $\CA$ (generalizing that of \cite{F, Shuf}) which is not directly visible in the quantum toroidal picture. More concretely, we construct a factorization:
\begin{equation}
\label{eqn:product}
\Ap \ = \ \prod^{\rightarrow}_{\mu \in \BQ} \Bp_{\mu}
\end{equation}
by which we understand that elements of $\Ap$ can be written uniquely as finite sums of products of elements of the subalgebras $\Bp_{\mu} \subset \Ap$, in increasing order of $\mu\in \BQ$. This factorization respects the bialgebra pairings on both sides, and this implies that the universal$^*$ $R-$matrix \footnote{The object we call the universal$^*$ $R$--matrix differs from the actual universal $R$--matrix by $q$ raised to certain Cartan elements, see Subsection \ref{sub:R} for more details. In the classical theory of quantum groups, the analogous notion is referred to as ``quasi $R$--matrix" (\cite{L})} of $\CA$ is the product of the universal$^*$ $R-$matrices for the Drinfeld doubles $\CB_\mu$ of the subalgebras in \eqref{eqn:product}:
\begin{equation}
\label{eqn:R0}
R_{\CA} \ = \ \prod^{\rightarrow}_{\mu \in \BQ \cup \{\infty\}} R_{\CB_{\mu}} \in \CA \widehat{\otimes} \CA
\end{equation}
The subalgebra $\CB_\infty$ consists of Cartan elements that must be added to the shuffle algebra in order to make it a well-defined bialgebra. For finite $\mu$, the subalgebras $\Bp_{\mu}$ are defined by the limit properties \eqref{eqn:limit1}, and we will show in Lemma \ref{lem:sub} that:
\begin{equation}
\label{eqn:smalliso0}
\Xi \ : \ \uux^{\otimes g} \stackrel{\sim}\longrightarrow \CB_{\frac ba}
\end{equation}
where $\gcd(a,b)=1$ and $g = \gcd(n,a)$. Here and throughout the paper, $\uu$ denotes the quantum affine algebra. Taking into account Theorem \ref{thm:iso}, formula \eqref{eqn:smalliso0} constructs embeddings of quantum affine algebras into the quantum toroidal algebra, associated to any rational slope $\mu = b/a$. Moreover, the isomorphisms $\Upsilon$ and $\Xi$ respect the bialgebra structures, and thus \eqref{eqn:R0} and \eqref{eqn:smalliso0} imply the following formula for the universal$^*$ $R-$matrix of the quantum toroidal algebra: \\

\begin{corollary}

The universal$^*$ $R$--matrix of the quantum toroidal algebra factors:
\begin{equation}
\label{eqn:R}
R_{\UU} = \prod_{\frac ba \in \BQ \cup \{\infty\}}^{\rightarrow} R_{\uux^{\otimes g}} \ \in \ \UU \widehat{\otimes} \UU
\end{equation}
as a product of universal$^*$ $R$--matrices for quantum affine algebras. \\

\end{corollary}

\noindent Formula \eqref{eqn:R} can be thought of as a toroidal analogue of the factorization of quantum (affine) algebra $R-$matrices from \cite{KT, LSS}. The shuffle algebra viewpoint also has the advantage that one can write down elements explicitly as symmetric rational functions, and in particular we will show that the isomorphism \eqref{eqn:smalliso0} sends the root generators of the various $U_q(\dot{\fgl}_{\frac ng})$ to particular shuffle elements of the form:
\begin{align}
&S_m^\pm = \sym \left[ \frac {m(z_i,...,z_{j-1})}{\left(1 - \frac {z_{i+1}}{z_{i}q^2}  \right) ... \left(1 - \frac {z_{j-1}}{z_{j-2}q^2} \right)} \prod_{i\leq a < b < j} \zeta \left( \frac {z_b}{z_a} \right)  \right] \in \CA^\pm \label{eqn:defs} \\
&T_m^\pm \ = \ \sym \left[ \frac {m(z_i,...,z_{j-1})}{\left(1 - \frac {z_{i}}{z_{i+1}} \right) ... \left(1 - \frac {z_{j-2}}{z_{j-1}} \right)} \prod_{i\leq a < b < j} \zeta \left( \frac {z_a}{z_b} \right)  \right] \in \CA^\pm \label{eqn:deft}
\end{align}
where $m$ is a Laurent polynomial, $\zeta$ is given by \eqref{eqn:bigzeta} and the notation for rational functions is that of Subsection \ref{sub:min}. We will show in \cite{Aff} that the shuffle elements \eqref{eqn:defs} and \eqref{eqn:deft} play a very important role in the geometric representation theory of affine Laumon spaces. We use this fact in \cite{CM} to compute the Nekrasov partition function of ${{\mathcal{N}}}=2$ supersymmetric $U(n)$ gauge theory with bifundamental matter in the presence of a complete surface operator. \\ 

\noindent I would like to thank Alexander Braverman, Boris Feigin, Michael Finkelberg, Sachin Gautam, Andrei Okounkov, Valerio Toledano Laredo and Alexander Tsymbaliuk for many wonderful discussions that taught me a lot of mathematics. I gratefully acknowledge the support of NSF grant DMS--1600375. \\

\section{Quantum Algebras}
\label{sec:tor}

\subsection{}\label{sub:reptheory}

All algebraic structures $A$ that will be studied in this paper are \underline{bialgebras} over a field $\BF$, meaning that they are endowed with a product and a coproduct:
$$
A \otimes A \stackrel{*}\longrightarrow A \qquad \qquad A \stackrel{\Delta}\longrightarrow A \otimes A
$$
which are associative and coassociative, respectively. All our bialgebras will be endowed with a unit $1 : \BF \rightarrow A$ and a counit $\e : A \rightarrow \BF$, although these will usually be obvious from the situation and we will not bother writing them down explicitly. There are a host of properties one expects from the above data, but the most important one is the compatibility between product and coproduct:
\begin{equation}
\label{eqn:compatibility}
\Delta(a * a') = \Delta(a) * \Delta(a')
\end{equation}
$\forall a,a' \in A$. We will often use Sweedler notation for the coproduct, namely:
\begin{equation}
\label{eqn:sweedler}
\Delta(a) = a_{1} \otimes a_{2}
\end{equation}
$\forall a \in A$, which implies the existence of a hidden summation sign in front of the tensor in the right-hand side (so the full notation would be $\Delta(a) = \sum_i a_{1,i} \otimes a_{2,i}$). Then \eqref{eqn:compatibility} can be written as:
$$
(aa')_1 \otimes (aa')_2 = a_1a_1' \otimes a_2 a_2'
$$
A \underline{Hopf algebra} is a bialgebra that is endowed with an antipode map:
$$
S : A \rightarrow A
$$
which satisfies certain compatibility properties with the product, coproduct, unit and counit. All the bialgebras in this paper are Hopf algebras, and in fact the antipode map will be determined from the compatibility properties. However, we will not explicitly write down the antipode, since it will be of no use to us. \\

\subsection{}\label{sub:drinfeldouble}

Given bialgebras $A^-$ and $A^+$, a \underline{bialgebra pairing} between them:
\begin{equation}
\label{eqn:bialgpair}
\langle \cdot, \cdot \rangle : A^- \otimes A^+ \rightarrow \BF
\end{equation}
is an $\BF$--linear pairing which satisfies the properties:
\begin{align}
&\left \langle a * a', b \right \rangle = \left \langle a\otimes a', \Delta(b) \right \rangle \label{eqn:bialg 1} \\
&\left \langle a,b * b' \right \rangle  = \left \langle \Delta^{\text{op}}(a), b \otimes b' \right \rangle \label{eqn:bialg 2}
\end{align}
for all $a,a' \in A^-$ and $b,b' \in A^+$ (where $\Delta^{\text{op}}$ denotes the opposite coproduct) and:
\begin{align*}
&\langle a,1 \rangle = \e(a) \\
&\langle 1,b \rangle = \e(b)
\end{align*}
If there is no danger of confusion, we will take the liberty to write $\langle b,a \rangle := \langle a,b\rangle$ for any $a \in A^-$ and $b \in A^+$. If $A^-$ and $A^+$ are also Hopf algebras, then the pairing is called Hopf if it satisfies the additional property:
\begin{equation}
\label{eqn:hopfy}
\left \langle S(a),b \right \rangle = \left \langle a,S^{-1}(b) \right \rangle
\end{equation}
$\forall a \in A^-, b\in A^+$. All the bialgebra pairings that will feature in this paper will be Hopf pairings, although we will not be concerned with this extra structure. \\

\begin{definition}
\label{def:drinfeld}
	
To any two bialgebras $A^-$ and $A^+$ with a bialgebra pairing \eqref{eqn:bialgpair} between them, one associates their \underline{Drinfeld double} (\cite{Drin}):
$$
A = A^- \otimes A^+
$$ 
It has the property that $A^- \otimes 1$ and $1 \otimes A^+$ are both sub-bialgebras of $A$, and they generate $A$ subject to the relations:
\begin{equation}
\label{eqn:drinfeld}
\langle a_1,b_1 \rangle a_2 * b_2 = b_1 * a_1 \langle a_2,b_2 \rangle
\end{equation}
$\forall a \in A^-, b\in A^+$, where in \eqref{eqn:drinfeld} we use Sweedler notation \eqref{eqn:sweedler} for $\Delta(a)$ and $\Delta(b)$. \\	
\end{definition}

\begin{remark}
\label{rem:double}

If the above does not look like the usual definition of the Drinfeld double (as appears, for example, in \cite{K}), it is because of convenience: we have sought a formulation which does not involve writing down the antipode map. Indeed, if we multiply \eqref{eqn:drinfeld} on the right with $\langle a_3, S^{-1}(b_3) \rangle$, we obtain:
$$
\langle a_1,b_1 \rangle a_2 * b_2 \langle a_3, S^{-1}(b_3) \rangle = b_1 * a_1 \langle a_2,b_2 \rangle \langle a_3, S^{-1}(b_3) \rangle = 
$$
$$
= b_1 * a_1 \langle \Delta(a_2), b_2 \otimes S^{-1}(b_3) \rangle = b_1 * a_1 \langle a_2 , S^{-1}(b_3)* b_2 \rangle =  
$$
$$
= b_1 * a_1 \langle a_2 , \e(b_2) \cdot 1 \rangle = b_1\e(b_2) * a_1 \e(a_2) = b * a
$$
where in the middle equality on the second row we have used \eqref{eqn:bialg 2}. The equation above is the more standard definition of the Drinfeld double, see \loccitt \\

\end{remark}

\subsection{}\label{sub:rmatrix}

A universal $R$--matrix of a bialgebra $A$ is an element $R \in A \otimes A$ such that:
\begin{equation}
\label{eqn:prop1}
R * \Delta(a) = \Delta^{\text{op}}(a) * R
\end{equation}
for all $a\in A$, and:
\begin{align}
&(\Delta \otimes 1)R = R_{13} * R_{23} \label{eqn:prop2} \\
&(1 \otimes \Delta)R = R_{13} * R_{12} \label{eqn:prop3}
\end{align}
where $R_{12} = R \otimes 1$, $R_{23} = 1 \otimes R$, and $R_{13}$ is defined analogously. Property \eqref{eqn:prop1} implies that for any representations $V,W\in \text{Rep}(A)$, the operator $R_{VW}$ given by: 
$$
A \otimes A \longrightarrow \text{End}(V \otimes W), \qquad R \leadsto R_{VW}
$$ 
intertwines the $A$--modules $V \otimes W$ and $W \otimes V$, which explains the terminology ``universal" and ``matrix". When $A$ is presented as a Drinfeld double $A^- \otimes A^+$ as in Definition \ref{def:drinfeld} with a non-degenerate pairing, a universal $R-$matrix always exists: \\

\begin{proposition} 
\label{prop:univR}
	
Let $\{a_i\}$ and $\{b_i\}$ be dual bases of $A^-$ and $A^+$ with respect to the bialgebra pairing \eqref{eqn:bialgpair}, which we assume to be non-degenerate. Then the tensor:
\begin{equation}
\label{eqn:rmatrix}
R = \sum_i b_i \otimes a_i \ \in \ A \otimes A
\end{equation}
is a universal $R-$matrix. This definition does not depend on the choice of dual bases, since \eqref{eqn:rmatrix} is nothing but the canonical tensor of the bialgebra pairing. \\
	
\end{proposition}

\begin{remark} 
\label{rem:completion}

Proposition \ref{prop:univR} tacitly assumed that $A^\pm$ are finite-dimensional over $\BF$, which will not be the case in the present paper. Instead, the algebras $A^\pm$ studied herein are graded by a certain cone in $\zz$, and the graded pieces are endowed with a filtration by finite-dimensional vector spaces. Therefore, it makes sense to consider the completion $A^- \woo A^+$, and to define the right-hand side \eqref{eqn:rmatrix} in this completion. \\

\end{remark}

\subsection{}\label{sub:intro}

Let $n>1$. Consider the semigroup $\nn$ of $n-$tuples of non-negative integers, whose elements will be denoted by $\bk = (k_1,...,k_n)$. Consider the partial ordering:
\begin{equation}
\label{eqn:partialorder}
\bl \leq \bk \qquad \Leftrightarrow \qquad l_i \leq k_i \quad \forall i\in \{1,...,n\}
\end{equation}
Elements $\bk \in \nn$ will sometimes be called \underline{degree vectors}. A particularly important example of degree vector is $\bk=[i;j)$ for integers $i < j$, defined by:
\begin{equation}
\label{eqn:part}
k_a = \# \Big\{ \textrm{integers} \equiv a \textrm{ mod }n \textrm{ in } \{i,...,j-1\} \Big\}
\end{equation}
Let us introduce the bilinear forms: 
\begin{align}
&\langle \cdot, \cdot \rangle : \nn \otimes \nn \longrightarrow \BZ, \qquad \langle \bk, \bl \rangle = \sum_{i=1}^n k_il_i - k_{i-1}l_{i} \label{eqn:bin} \\
&( \cdot, \cdot ) : \nn \otimes \nn \longrightarrow \BZ, \qquad (\bk, \bl) = \sum_{i=1}^n 2k_il_i - k_{i-1} l_i - k_i l_{i-1} \label{eqn:symbin}
\end{align}
where we identify $k_0 = k_n$ and $l_0 = l_n$. Also define the degree vectors:
\begin{equation}
\label{eqn:real}
\bs^i \ = \ \underbrace{(0,...,0,1,0,...,0)}_{1\text{ on the }i-\text{th position}} \ \in \nn
\end{equation}
and:
\begin{equation}
\label{eqn:imaginary}
\bde = (1,...,1) = \bs^1 + ... + \bs^n
\end{equation}
Note that $\bde$ spans the kernels of the forms \eqref{eqn:bin} and \eqref{eqn:symbin}. \\

\subsection{}\label{sub:quantumgroup}

In this Section, we will introduce the algebras $\su$, $\uu$ and $\UU$ for $n\geq 2$. The base field will be implicit in the choice of the subscripts: we choose $\BC(q)$ in the first two cases and $\BC(q,\oq)$ in the last case, where $q, \oq$ are formal symbols. The algebra $\sup$ is generated by symbols $x_1,...,x_n$ under the relations:
\begin{equation}
\label{eqn:sug1}
[x_i,x_j] = 0 \qquad \qquad \qquad \qquad \qquad \forall i - j \not \equiv \{-1,1\} \text{ mod } n
\end{equation}
\begin{equation}
\label{eqn:sug2}
x_i^2x_{i\pm 1} - (q+q^{-1})x_ix_{i\pm 1}x_i + x_{i\pm 1} x_i^2 = 0 \ \quad \forall i \in \{1,...,n\}
\end{equation}
(above, we identify $x_0 = x_n$ and $x_1 = x_{n+1}$). We define the extended algebra: 
$$
\sug = \Big \langle \sup, \ph_1,...,\ph_n \Big \rangle \Big/_{[\ph_i,\ph_j] = 0 \text{ and relation \eqref{eqn:sug3}}}
$$
where for all $i,j \in \{1,...,n\}$ we set:
\begin{equation}
\label{eqn:sug3}
\ph_j x_i = q^{(\bs^i,\bs^j)} x_i \ph_j
\end{equation}
In classical notation, $\ph_i = K_i$. As a consequence of \eqref{eqn:sug3}, note that the element:
\begin{equation}
\label{eqn:central}
c = \ph_1...\ph_n
\end{equation}
is central. There is a bialgebra structure on $\sug$ with coproduct given by:
\begin{align}
&\Delta(x_i) = \ph_i \otimes x_i + x_i \otimes 1 \label{eqn:cop 1} \\
&\Delta(\ph_i) = \ph_i \otimes \ph_i \label{eqn:cop 2}
\end{align}
Note the bialgebra automorphism $\sug \rightarrow \sug$ that sends $x_{i} \mapsto x_{i+1}$ and $\ph_i \mapsto \ph_{i+1}$. Because $x_{i+n} = x_i$ and $\ph_{i+n} = \ph_i$, this automorphism has order $n$. \\

\subsection{}\label{sub:quantumdouble} 

We set $\sul = \sug^\coop$, namely the same algebra with the opposite coproduct, and use the notation $x^-$ (resp. $\ph^-$) and $x^+$ (resp. $\ph^+$) to differentiate between elements of $\sul$ and $\sug$. There exists a bialgebra pairing:
$$
\langle \cdot, \cdot \rangle \ : \ \sul \otimes \sug \ \longrightarrow \ \BQ(q)
$$
completely determined by the assignments:
\begin{align}
&\left \langle x_i^-, x_j^+ \right \rangle = \frac {\delta_j^i}{q^{-1} - q} \label{eqn:pair 1} \\ 
&\left \langle \ph_i^-,\ph_j^+ \right \rangle = q^{(\bs^i,\bs^j)} \label{eqn:pair 2}
\end{align}
and relations \eqref{eqn:bialg 1}, \eqref{eqn:bialg 2}. Then we may define the Drinfeld double associated to this data, according to Definition \ref{def:drinfeld}. We set:
$$
\su := \sul \otimes \sug \Big |_{\ph^+_i \ph^-_i =1}
$$
We will often write $\ph_i = \ph_i^+$. Note that \eqref{eqn:drinfeld}, \eqref{eqn:cop 1}, \eqref{eqn:cop 2}, \eqref{eqn:pair 1}, \eqref{eqn:pair 2} imply the well-known commutation relation between the positive and negative generators: 
\begin{equation}
\label{eqn:basiccomm}
\left [x_i^+,x_j^- \right ] = \delta_j^i \cdot \frac {\ph_i - \ph_i^{-1}}{q-q^{-1}} 
\end{equation}
The algebra $\su$ defined above is nothing but the \underline{quantum group} associated to the Cartan matrix of the cyclic quiver with $n$ vertices. The bialgebra automorphism:
\begin{align*}
&\su \longrightarrow \su \\
&x_i^\pm \mapsto x_{i+1}^\pm, \ph_i \mapsto \ph_{i+1}
\end{align*}
comes from the order $n$ rotation of the cyclic quiver. We will write $\BZ/n\BZ \curvearrowright \su$. \\

\subsection{}\label{sub:alg2}

The larger quantum affine algebra $\uu$ admits a description as above, but we will prefer to start with its RTT presentation (see \cite{FRT}, \cite{RS}, \cite{FD}, \cite{GM}). In order to present its generators and relations, we appeal to an equivalent reformulation of \loccit Consider the following tensor product of $n \times n$ matrices valued in $\BQ(q,z)$:
\begin{equation}
\label{eqn:r}
R \left( z \right) = \sum_{1 \leq i, j \leq n} E_{ii} \boxtimes E_{jj} \left( \frac {q - zq^{-1}}{1 - z} \right)^{\delta_j^i} + (q - q^{-1}) \sum_{1\leq i \neq j \leq n } E_{ij} \boxtimes E_{ji} \frac {z^{\delta_{i>j}}}{1 - z} \qquad \quad
\end{equation}
where $E_{ij}$ denotes the elementary matrix with a single entry 1 at the intersection of row $i$ and column $j$. We denote the tensor product of matrices by $\boxtimes$ in order to distinguish it from the tensor product $\otimes$ in the definition of the coproduct $\Delta$. We will call \eqref{eqn:r} an $R-$matrix, because it satisfies the Yang-Baxter equation:
\begin{equation}
\label{eqn:yb}
R_{12} \left(\frac xy \right) R_{13} \left(\frac xz \right) R_{23} \left(\frac yz \right) = R_{23} \left(\frac yz \right) R_{13} \left(\frac xz \right) R_{12} \left(\frac xy \right)
\end{equation}
among the triple tensors $R_{12} = R \boxtimes \textrm{Id}$, $R_{23} = \text{Id} \boxtimes R$ etc. We leave \eqref{eqn:yb} as an easy exercise for the interested reader, and note that it reduces to the analogous computation carried out in \cite{GM}. Indeed, \eqref{eqn:r} equals the $R-$matrix given in (2.41) of \loccitt, upon multiplication by $(u-v)^{-1}$ and the substitution $u/v = z$. \\

\subsection{}\label{sub:rtt}
 
Let us define the following algebras: 
\begin{align}
&\uug = \Big \langle e_{[i;j)}, \psi_1,...,\psi_n, c \Big \rangle^{1\leq i \leq n}_{i<j} \Big/_{\text{relation \eqref{eqn:rtt} for } \pm = +} \label{eqn:gen1} \\
&\uul = \Big \langle f_{[i;j)}, \psi^{-1}_1,...,\psi^{-1}_n, c^{-1} \Big \rangle^{1\leq i \leq n}_{i<j} \Big/_{\text{relation \eqref{eqn:rtt} for } \pm = -} \label{eqn:gen2}
\end{align}
where $c$ is central and the $\psi_i$ all commute between themselves. We set: 
$$
e_{[i;j)}=e_{[i-n;j-n)} \qquad \quad f_{[i;j)}=f_{[i-n;j-n)} \qquad \quad \psi_{i} = c \psi_{i-n}
$$
for all $i < j \in \BZ$, so our indices may be arbitrary integers. We will call $e_{[i;j)}$ and $f_{[i;j)}$ the \underline{root generators}, and in order to present the relations between them, let us introduce the power-series valued matrices:
\begin{align*}
&T^+(z) = \sum_{1\leq i \leq n}^{i\leq j} e_{[i;j)} \psi_i \cdot E_{\bar{j} i} \ z^{-\left \lfloor \frac {j-1}n \right \rfloor} \quad \in \uug \otimes \text{Mat}_{n \times n} [[z^{-1}]] \\
&T^-(z) = \sum_{1 \leq i \leq n}^{i\leq j} f_{[i;j)} \psi_i^{-1} \cdot E_{i \bar{j}} \ z^{\left \lfloor \frac {j-1}n \right \rfloor} \quad \in \uul \otimes \text{Mat}_{n \times n} [[z]] 
\end{align*}
where:
\begin{equation}
\label{eqn:residue mod n}
\bar{i} = \left( i \text{ mod }n \right) \in \{1,...,n\}
\end{equation}
for all $i \in \BZ$. For convenience, we will write: 
\begin{align*}
&e_{[i;i)} = f_{[i;i)} = 1 \\
&e_{[i;j)} = f_{[i;j)} = 0
\end{align*}
for all integers $i > j$. \\

\begin{definition}
\label{def:rtt}

The algebras $\uug$ and $\uul$ are given by generators as in \eqref{eqn:gen1} and \eqref{eqn:gen2}, modulo the so-called RTT relations:
\begin{equation}
\label{eqn:rtt}
R\left(\frac xy \right) T^\pm_1(x) T^\pm_2(y) = T^\pm_2(y) T^\pm_1(x) R \left(\frac xy \right)
\end{equation}
where $T^\pm_1(x) = T^\pm(x) \boxtimes \text{Id}$ and $T^\pm_2(y) = \text{Id} \boxtimes T^\pm(y)$. \\

\end{definition}

\noindent As shown in \cite{thesis}, \eqref{eqn:rtt} implies the following explicit relations:
\begin{align}
&\psi_k e_{[i;j)} = q^{\delta_{\barj}^{\bark} - \delta_{\bari}^{\bark}} e_{[i;j)} \psi_k \label{eqn:expre 1} \\
&\psi_k f_{[i;j)} = q^{\delta_{\bari}^{\bark} - \delta_{\barj}^{\bark}} f_{[i;j)} \psi_k \label{eqn:expre 2} 
\end{align}
\begin{multline}
\label{eqn:expre 3} \frac {e_{[i;j)} e_{[i';j')}}{q^{\delta_{\bari}^{\barj} - \delta_{\bari}^{\barjj} - \delta_{\barj}^{\barjj}}} - \frac {e_{[i';j')} e_{[i;j)}}{q^{\delta_{\bari}^{\barj} - \delta_{\barj}^{\barii} - \delta_{\bari}^{\barii}}} = \\ = (q-q^{-1}) \left(\sum_{i' < a \leq j'}^{a \equiv i} e_{[a,j')} e_{[i+i'-a;j)} - \sum_{i' \leq a < j'}^{a \equiv j} e_{[i;j+j'-a)} e_{[i',a)}\right)
\end{multline}
\begin{multline}
\label{eqn:expre 4} \frac {f_{[i;j)} f_{[i';j')}}{q^{\delta_{\barjj}^{\barii} - \delta_{\barii}^{\bari} - \delta_{\barjj}^{\bari}}} - \frac {f_{[i';j')} f_{[i;j)}}{q^{\delta_{\barjj}^{\barii} - \delta_{\barii}^{\barj} - \delta_{\barjj}^{\barj}}} = \\ = (q-q^{-1}) \left(\sum_{i \leq a < j}^{a \equiv j'} f_{[i';j+j'-a)} f_{[i;a)} - \sum_{i<a\leq j}^{a \equiv i'} f_{[a;j)} f_{[i+i'-a;j')}\right) 
\end{multline}
for all integer indices $k$, $i < j$, $i' < j'$. Recall that $\bar{i} = i \text{ mod }n$ as in \eqref{eqn:residue mod n}. \\

\subsection{}\label{sub:drinfeld}

The algebras $\uug$ and $\uul$ are bialgebras, with coproduct given by:
\begin{equation}
\label{eqn:copy1}
\Delta(T^+(z)) = T^+(z) \otimes T^+(z/c_1) 
\end{equation}
\begin{equation}
\label{eqn:copy2}
\Delta(T^-(z)) = T^-(z/c_2) \otimes T^-(z) 
\end{equation}
and counit $\e(T^\pm(z)) = \text{Id}$, where the central elements are $c_1 = c\otimes 1$ and $c_2 = 1\otimes c$. We define a pairing by $\langle c,a \rangle = \langle a, c \rangle = \e(a)$ for all $a \in \uug, \uul$, and:
\begin{equation}
\label{eqn:bialg2}
\Big \langle T_1^-(x), T_2^+(y) \Big \rangle = R\left( \frac xy \right)
\end{equation}
where we expand in non-negative powers of $x/y$. As we will show in Proposition \ref{prop:rtt}, these choices extend to a well-defined bialgebra pairing between $\uug$ and $\uul$. Explicitly in terms of the root generators, relations \eqref{eqn:copy1}--\eqref{eqn:copy2} entail:
\begin{align}
&\Delta(\psi_i) = \psi_i \otimes \psi_i \label{eqn:quant1} \\
&\Delta \left( e_{[i;j)} \right) =  \sum_{a=i}^j e_{[a;j)}\frac {\psi_a}{\psi_i} \otimes e_{[i;a)} \label{eqn:quant2}  \\
&\Delta \left( f_{[i;j)} \right) = \sum_{a=i}^j f_{[i;a)} \otimes f_{[a;j)}\frac {\psi_i}{\psi_a} \label{eqn:quant3} 
\end{align}
for all $i<j$, while \eqref{eqn:bialg2} implies:
\begin{equation}
\label{eqn:quant4} 
\langle \psi_i^{-1}, \psi_j \rangle = q^{\delta_j^i} 
\end{equation} 
\begin{equation} 
\label{eqn:quant5}
\Big \langle  e_{[i;j)}, f_{[i';j')} \Big \rangle = \delta_{[i';j')}^{[i;j)} \left(1 - q^{-2} \right) 
\end{equation}
for all $i < j$, $i' < j'$. The Kronecker delta in \eqref{eqn:quant5} is 1 iff $(i,j)-(i',j') \in \BZ(n,n)$. The following result will be proved in the Appendix. \\

\begin{proposition}
\label{prop:rtt}

The pairing \eqref{eqn:bialg2} generates a well-defined bialgebra pairing. This allows us to define the Drinfeld double:
$$
\uu := \uul \otimes \uug \Big|_{(\psi_i \otimes 1)(1 \otimes \psi_i) = (c \otimes 1)(1\otimes c) = 1}
$$
We will write $\psi_i = \psi_i \otimes 1$ and $c = c \otimes 1$. The defining relations among the two tensor factors of the Drinfeld double are:
\begin{equation}
\label{eqn:unwind}
R \left( \frac x{yc} \right) T_1^-(x) T_2^+(y) = T_2^+(y) T_1^-(x) R\left(\frac {xc}y\right)
\end{equation}
\text{}
 
\end{proposition}

\noindent Explicitly, \eqref{eqn:unwind} implies the following commutation relations, for all $i<j, i'<j'$:
\begin{multline}
\label{eqn:expre 5}
\Big[ e_{[i;j)}, f_{[i';j')} \Big] = (q-q^{-1}) \cdot \\ \sum_{k=1}^{\min(j-i,j'-i')} \left( \delta_{\bark+\bari'}^{\barj} \frac {f_{[i'+k,j')} e_{[i;j-k)}}{q^{\delta_{\barii}^{\bari} - \delta_{\barii}^{\barj} + \delta_{\barj}^{\bari}}} \frac {\psi_{i'}}{\psi_{i'+k}} - \delta_{\bark+\bari}^{\barjj} \frac {e_{[i+k,j)} f_{[i';j'-k)}}{q^{\delta_{\bari}^{\barii} + \delta_{\bari}^{\barjj} - \delta_{\barii}^{\barjj}}} \frac {\psi_{i+k}}{\psi_{i}} \right)
\end{multline}
Formula \eqref{eqn:expre 5} was proved in \cite{thesis}, and will not be used in the present paper. \\

\subsection{}\label{sub:pbw} 

Note the action $\BZ/n\BZ \curvearrowright \uupm$ generated by the bialgebra automorphism:
$$
e_{[i;j)} \mapsto e_{[i+1;j+1)} \qquad \qquad \qquad f_{[i;j)} \mapsto f_{[i+1;j+1)}
$$
which has order $n$. Moreover, the algebra $\uu$ is $\zz$--graded, with:
$$
\deg e_{[i;j)} = [i;j) \in \nn \quad \ \qquad \deg f_{[i;j)} = -[i;j) \in -\nn \quad \ \qquad \deg \psi_i = 0
$$
A PBW basis of $\uu$ was constructed in \cite{GM}, and we will give a related construction in the Appendix. Specifically, one unwinds relations \eqref{eqn:expre 3}--\eqref{eqn:expre 4} to show that any product of root generators $e_{[i;j)}$ (respectively $f_{[i;j)}$) is equal to a linear combination of products where the root generators are placed in a certain order. This implies that the ordered products of root generators form a linear generating set of the following subalgebras:
\begin{align*}
&\uup = \Big \langle e_{[i;j)} \Big \rangle_{1\leq i \leq n}^{i<j} \subset \uug \\
&\uum = \Big \langle f_{[i;j)} \Big \rangle_{1\leq i \leq n}^{i<j} \subset \uul
\end{align*}
This allows us to estimate the dimensions of the graded pieces (see the Appendix): \\

\begin{proposition} 
\label{prop:dim}

For all degree vectors $\bk\in \nn$ we have:
\begin{equation}
\label{eqn:dimboundaff}
\dim \ U_q^\pm(\dot{\fgl}_n)_{\pm \bk} \ \leq \ \# \text{ partitions of } \bk
\end{equation}
into degree vectors $[i;j)$ for various $1\leq i \leq n$, $i<j$. Alternatively, we may consider all integers $i < j$, but identify the degree vectors $[i;j)$ and $[i-n;j-n)$. \\

\end{proposition}

\noindent In fact, it will follow from Proposition \ref{prop:iso} and Remark \ref{rem:iso} that the inequality \eqref{eqn:dimboundaff} is actually an equality. This will be proved in Section \ref{sec:shuf} rather indirectly, by appealing to the shuffle algebra incarnation of quantum affine algebras. \\

\subsection{}\label{sub:heisenberg}

\noindent When $n=1$, $\uui$ is called the \underline{quantum Heisenberg algebra}. Let us write $g_k = e_{[1;k+1)}$ and $g_{-k} = f_{[1;k+1)}$ for all $k \geq 0$, and observe that relation \eqref{eqn:rtt} allows one to show that the $g_{k}$'s and the $g_{-k}$'s all commute. We conclude that:
\begin{align}
&\uuig = \BQ(q) \left [ g_k, c \right ]_{k\in \BN} \label{eqn:heis1} \\
&\uuil = \BQ(q) \left [ g_{-k}, c^{-1} \right ]_{k\in \BN} \label{eqn:heis2}
\end{align}
\footnote{Strictly speaking, the algebras $\uuig$ and $\uuil$ should also contain the Cartan element $\psi_1$. However, since this element is central and plays no role in the coproduct, we slightly abuse notation by removing it from the definition of the aforementioned algebras}
are both commutative algebras. The coproduct relations \eqref{eqn:quant2} and \eqref{eqn:quant3} imply: 
\begin{align}
\Delta(g_k) &= \sum^{a,b\geq 0}_{a+b = k} g_a c^b \otimes g_b \label{eqn:group1} \\
\Delta(g_{-k}) &= \sum^{a,b\geq 0}_{a+b = k} g_{-a} \otimes g_{-b} c^{-a} \label{eqn:group2}
\end{align}
for any $k\in \BN$. We call the elements $g_{\pm k}$ \underline{group-like}\footnote{Usually, group-like and primitive elements are ones which satisfy $\Delta(p) = p \otimes 1 + 1 \otimes p$ and $\Delta(g) = g \otimes g$, respectively. In the case of quantum algebras, we must relax this notion to allow various powers of Cartan elements, such as $c$ in the two factors of the coproduct}, although we remark that it would be more appropriate to assign this term to the formal sum $1+g_{\pm 1}+g_{\pm 2}+...$. We may replace the generators $g_{\pm k}$ by the generators $p_{\pm k}$ defined by:
\begin{align}
\sum_{k = 0}^\infty g_{k} z^k &= \exp \left[\sum_{k=1}^\infty \frac {p_k z^k}k \left(  q^{-k} - q^{k} \right) \right] \label{eqn:generating 1} \\
\sum_{k = 0}^\infty g_{- k} z^k &= \exp \left[\sum_{k=1}^\infty \frac {p_{- k} z^k}k \left( 1-q^{-2k} \right) \right] \label{eqn:generating 2}
\end{align}
Then \eqref{eqn:group1} and \eqref{eqn:group2} imply that:	
\begin{align}
&\Delta(p_k)  = c^k \otimes p_k + p_k \otimes 1 \label{eqn:prim1} \\
&\Delta(p_{-k}) = 1 \otimes p_{-k} + p_{-k} \otimes c^{-k} \label{eqn:prim2}
\end{align}
for all $k\in \BN$. We call the elements $p_{\pm k}$ \underline{primitive} (although this term is imprecise, see the footnote about group-like elements). The pairing \eqref{eqn:quant5} takes the form:
$$
\langle g_{-k}, g_{l} \rangle = \delta_k^l ( 1-q^{-2} )
$$
for all $(k,l) \neq (0,0)$, which implies (by a straightforward computation that we leave to the interested reader):
\begin{equation}
\label{eqn:pie}
\langle p_{-k}, p_{l} \rangle = \frac {\delta_k^l k}{q^{-k} - q^{k}}
\end{equation}
for all $k > 0$. Therefore, relation \eqref{eqn:drinfeld} together with \eqref{eqn:prim1}--\eqref{eqn:pie} imply the following well-known commutation relation for the $q$--Heisenberg algebra:
\begin{equation}
\label{eqn:comm heis}
[p_k, p_l ] = \delta_{k+l}^0 k \cdot \frac {c^k - c^{-k}}{q^k - q^{-k}}
\end{equation}
The property of an element being primitive is preserved under rescaling. In other words, if we perform the substitution:
\begin{equation}
\label{eqn:plethysm}
p_k \leadsto p_k a_k \quad \text{and} \quad p_{-k} \leadsto p_{-k} b_k
\end{equation}
for arbitrary scalars $a_k,b_k \in \BF \backslash 0$, then the only thing that changes in the present Subsection is that the right-hand sides of \eqref{eqn:pie} and \eqref{eqn:comm heis} need to be multiplied by the constant $a_kb_k$. Properties \eqref{eqn:prim1} and \eqref{eqn:prim2} would still hold, as would \eqref{eqn:group1} and \eqref{eqn:group2} if $g_{\pm k}$ are defined by applying the substitution \eqref{eqn:plethysm} to \eqref{eqn:generating 1}--\eqref{eqn:generating 2}. \\

\subsection{}\label{sub:nondeg} 

In the present Subsection, we will show that there exists an embedding:
\begin{equation}
\label{eqn:iota}
\su \ \hookrightarrow \ \uu
\end{equation}
obtained by sending:
\begin{equation}
\label{eqn:iota2}
x_i^+ \mapsto \frac {e_{[i;i+1)}}{q^{-1} - q} \qquad \qquad x_i^- \mapsto \frac {f_{[i;i+1)}}{1 - q^{-2}} \qquad \qquad \ph_i \mapsto \frac {\psi_{i+1}}{\psi_i} 
\end{equation}
The fact that \eqref{eqn:iota2} defines an algebra homomorphism \eqref{eqn:iota} is easy to prove, simply by comparing relations \eqref{eqn:sug1}--\eqref{eqn:sug2} with \eqref{eqn:expre 3}--\eqref{eqn:expre 4}, and relation \eqref{eqn:basiccomm} with \eqref{eqn:expre 5}. It is also easy to show that \eqref{eqn:iota} preserves the coproduct and pairing. The injectivity of \eqref{eqn:iota} follows from a well-known and easy to prove exercise: \\

\begin{lemma}
\label{lem:inj}

Any linear map $i: A \rightarrow B$ which respects pairings on $A$ and $B$:
$$
\langle x, y \rangle_A = \langle i(x), i(y) \rangle_B
$$
$\forall x,y\in A$, is injective if the pairing $\langle \cdot, \cdot \rangle_A$ is non-degenerate. \\

\end{lemma}

\noindent The non-degeneracy of the bialgebra pairing of $\su$ is a well-known result, see for example \cite{Jan}. Then Lemma \ref{lem:inj} implies that the restrictions of \eqref{eqn:iota} to the halves $\sug$ and $\sul$ of $\su$ are injective. Due to the triangular decompositions of $\su$ and $\uu$, this implies that \eqref{eqn:iota} is injective. \\

\noindent In Subsection \ref{sub:iso}, we will prove the following result: \\

\begin{proposition}
\label{prop:iso}

For each $k\in \BN$, there is a primitive, $\BZ/n\BZ-$invariant element:
$$
p_{\pm k} \in U_q(\dot{\fgl}_n)_{\pm k\bde}
$$
which is uniquely defined up to a constant multiple. We have a subalgebra: 
$$
\uui \cong \Big \langle p_{\pm k}, c^{\pm 1} \Big \rangle_{k\in \BN} \subset \uu
$$
which commutes with $\su \subset \uu$ of \eqref{eqn:iota}, yielding a bialgebra isomorphism:
\begin{equation}
\label{eqn:decomposition}
\su \otimes \uui \stackrel{\sim}\longrightarrow \uu
\end{equation}

\end{proposition}

\begin{remark}

Although the isomorphism \eqref{eqn:decomposition} was known for a long time as the deformation of the usual decomposition:
$$
\hsl_n \oplus \hgl_1 \cong \hgl_n 
$$
we believe that the construction and proof provided in the present paper are new. \\

\end{remark} 

\begin{remark}
\label{rem:iso}

Once one shows that a bialgebra homomorphism \eqref{eqn:decomposition} exists, the fact that it is bijective can be proved as follows. Because: 
$$
\left \langle U_q^\pm(\dot{\fsl}_n), U_q^\mp(\dot{\fgl}_1) \right \rangle = 0
$$
with respect to the bialgebra pairing on $U_q^\pm(\dot{\fgl}_n)$ (because their generators live in different degrees, and therefore have pairing 0 for trivial reasons), then the injectivity of \eqref{eqn:decomposition} is a consequence of Lemma \ref{lem:inj}. Meanwhile, surjectivity follows from Proposition \ref{prop:dim} and the fact that:
\begin{equation}
\label{eqn:mugabe}
\dim \ \Big[ U_q^\pm(\dot{\fsl}_n) \otimes U_q^\pm(\dot{\fgl}_1) \Big]_{\pm \bk} = \# \text{ partitions of } \bk
\end{equation}
into $[i;j)$ for various $i < j$. Indeed, the dimensions of the graded pieces of $U_q^\pm(\dot{\fsl}_n)$ are given by the number of partitions of $\bk$ into positive roots, counted with multiplicities. For the root system of affine type $A$, the multiplicity is 1 for roots of the form $[i;j)$ if $j \not \equiv i$ modulo $n$, and $n-1$ for roots of the form $k\bde$. Since $U_q^\pm(\dot{\fgl}_1)$ contributes an additional root in degree $k\bde$, this produces exactly the count \eqref{eqn:mugabe}. \\ 

\end{remark}

\subsection{}\label{sub:colorvar}

In the remainder of this Section, we will define the algebra $\UU$. In order to do so, we consider variables $z$ with certain ``colors":
$$
\col z  \in \{1,...,n\}
$$
In our formulas, we might often see expressions in terms of variables of any integer color. Whenever this happens, we make the convention that:
\begin{equation}
\label{eqn:quasi}
f(...,z \text{ of color } k,...) := f \left(...,  z\oq^{-2\left \lfloor \frac {k-1}n \right \rfloor} \text{ of color } \bark,... \right)
\end{equation}
\footnote{The parameter $\oq$ will be identified with one of the equivariant parameters of affine Laumon spaces in \cite{Aff}, the other equivariant parameter being the quantum parameter $q$} The following color-dependent rational function will be very important for us:
\begin{equation}
\label{eqn:bigzeta}
\zeta \left( \frac zw \right) = \left( \frac {zq \oq^{2 \left \lceil \frac {i-j}n \right \rceil} - w q^{-1}}{z \oq^{2 \left \lceil \frac {i-j}n \right \rceil} - w} \right)^{\delta_{\bari}^{\barj} - \delta_{\overline{i+1}}^{\barj}}
\end{equation}
for variables $z$, $w$ of any colors $i,j \in \BZ$. Note the following simple identity:
\begin{equation}
\label{eqn:identity}
\zeta \left(\frac zw \right) = \zeta \left(\frac w{z q^2} \right)
\end{equation}
if $z$ and $w$ have colors $i$ and $j$ in the LHS, but colors $i+1$ and $j$ in the RHS. \\

\subsection{}\label{sub:tor}

Consider the algebras:
$$
U^\pm_{q,\oq}(\ddot{\fgl}_n) = \BC(q,\oq) \Big \langle x^\pm_{i,d} \Big \rangle^{1 \leq i \leq n}_{d\in \BZ} \Big/_{\text{relations \eqref{eqn:reltor2} and \eqref{eqn:serre}}}
$$
as well as the extended algebras:
\begin{align*}
&\UUl = \left \langle \UUm, \psi_{i,d}^-, c^- \right \rangle^{1 \leq i \leq n}_{d \in \BN \sqcup 0} \Big /^{\ [\psi_{i,d}^-, \psi_{i',d'}^-] = 0, \ \psi_{i,0}^- \text{ invertible}, \ c^- \text{ central}}_{\text{and relation \eqref{eqn:reltor1} for }\pm = -} \\
&\UUg = \left \langle \UUp, \psi_{i,d}^+, c^+ \right \rangle^{1 \leq i \leq n}_{d \in \BN \sqcup 0}  \Big /^{\ [\psi_{i,d}^+, \psi_{i',d'}^+] = 0, \ \psi_{i,0}^+ \text{ invertible}, \ c^+ \text{ central}}_{\text{and relation \eqref{eqn:reltor1} for }\pm = +}
\end{align*}
To write down the relations between the generators, it makes sense to collect them into \underline{currents}, by which we mean bi-infinite power series:
\begin{align}
&x_i^\pm(z) = \sum_{d\in \BZ} \frac {x^\pm_{i,d}}{z^d} \label{eqn:currents1} \\
&\psi^\pm_i(z) = \sum_{d = 0}^\infty \frac {\psi^\pm_{i,d}}{z^{\pm d}} \label{eqn:currents2}
\end{align}
$\forall i \in \{1,...,n\}$, such that the leading terms $\psi^\pm_i = \psi^\pm_{i,0}$ are invertible. We will write:
\begin{align}
&x^\pm_{i}(z) = x^\pm_{i-n}(z \oq^2) \label{eqn:period 1} \\ 
&\psi^\pm_{i}(z) = \psi^\pm_{i-n}(z \oq^2) c^{\pm 1} \label{eqn:period 2} 
\end{align}
for all $i$, and hence we may think of the indices as being arbitrary integers. Set:
\begin{align}
\psi^\pm_j(w) x^\pm_i(z) &= x^\pm_i(z) \psi^\pm_j(w) \cdot \zeta\left( \frac zw \right)^{\pm 1} \label{eqn:reltor1} \\ 
x^\pm_i(z) x^\pm_{j}(w) \cdot \zeta\left( \frac wz \right)^{\pm 1} &= x^\pm_{j}(w) x^\pm_i(z) \cdot \zeta\left( \frac zw \right)^{\pm 1} \label{eqn:reltor2}
\end{align}
and the Serre relation:
$$
x^\pm_{i \pm' 1}(w) x^\pm_i(z) x^\pm_i(z') - (q+q^{- 1}) x^\pm_i(z) x^\pm_{i\pm' 1}(w) x^\pm_i(z') + x^\pm_i(z) x^\pm_i(z') x^\pm_{i\pm' 1}(w) +
$$
\begin{equation}
\label{eqn:serre}
+ \{ \text{same expression with } z \text{ and } z' \text{ switched} \} = 0
\end{equation}
for all $i,j$ and all signs $\pm, \pm'$. The function $\zeta$ that appears in formulas \eqref{eqn:reltor1}--\eqref{eqn:reltor2} is defined by giving the variables $z$ and $w$ colors $i$ and $j$, respectively, and we note that formulas \eqref{eqn:reltor1}, \eqref{eqn:reltor2}, \eqref{eqn:serre} are compatible with \eqref{eqn:period 1}--\eqref{eqn:period 2}. \\

\begin{remark}
\label{rem:note}

To make sense of \eqref{eqn:reltor1}, \eqref{eqn:reltor2} and \eqref{eqn:serre}, one must equate the coefficients of all $z^aw^b$ in the left and right hand sides. To make sense of this, in \eqref{eqn:reltor1}, one must first expand in non-positive powers of $w^{\pm 1}$, while in \eqref{eqn:reltor2}, one multiplies the equation by all denominators of the functions $\zeta^{\pm 1}$ and then equates coefficients. \\

\end{remark}

\noindent The algebras $\UUl$ and $\UUg$ are bigraded by $(\pm \nn) \times \BZ$:
\begin{align*}
&\bideg x^\pm_{i,d} = (\pm \bs^i , d) \\
&\bideg \psi^\pm_{i,d} = (0,\pm d)
\end{align*} 
We will write $\UUl_{\bk,d}$ and $\UUg_{\bk,d}$ for the bigraded components. \\

\subsection{}\label{sub:bitor}

$\UUl$ and $\UUg$ are bialgebras with coproducts given by: 
\begin{align}
&\Delta(\psi^\pm_i(z)) = \psi^\pm_i(z) \otimes \psi^\pm_i(z) \label{eqn:deltator1} \\
&\Delta(x^+_i(z)) = \ph_i^+(z) \otimes x^+_i(z) + x^+_i(z) \otimes 1 \label{eqn:deltator2} \\
&\Delta(x_i^-(z)) = 1 \otimes x^-_i(z) + x^-_i(z) \otimes \ph_i^-(z) \label{eqn:deltator3}
\end{align}
where we write: 
\begin{equation}
\label{eqn:ph}
\ph_i^\pm(z) := \frac {\psi_{i+1}^\pm(zq^2)}{\psi_i^\pm(z)}
\end{equation}
for all integers $i$. Recall the \underline{delta-function} defined as $\delta(y) = \sum_{d\in \BZ} y^d$. The following Proposition will be proved in the Appendix. \\

\begin{proposition}
\label{prop:pairtor}

There exists a bialgebra pairing:
\begin{equation}
\label{eqn:pair}
\langle \cdot, \cdot \rangle : \UUl \otimes \UUg \longrightarrow \BC(q,\oq) 
\end{equation}
generated by the formulas:
\begin{equation}
\label{eqn:pairtor1}
\langle x_j^-(w), x_i^+(z) \rangle = \frac {\delta_j^i}{q^{-1} - q} \cdot \delta\left(\frac zw \right) 
\end{equation}
\begin{equation}
\label{eqn:pairtor2}
\langle \psi^-_j(w),\psi^+_i(z) \rangle = \zeta \left( \frac {z_i}w \right) \zeta \left( \frac {z_{i+1}q^2}w \right) \zeta \left( \frac {z_{i+2}q^4}w \right)... \ \Big|^{z_i, z_{i+1}, z_{i+2}, ... \mapsto z} 
\end{equation}
for variables $z$ and $w$ of any colors $i$ and $j$. In the right-hand side of \eqref{eqn:pairtor2}, the variables $z_i$, $z_{i+1}$,... are assumed to have colors $i$, $i+1$,... The infinite product is a ratio of infinite $q-$Pochhammer symbols, which converge for $|q| \gg 1$. \\

\end{proposition}

\begin{definition}
\label{def:tor} The \underline{quantum toroidal algebra} is the Drinfeld double:
\begin{equation}
\label{eqn:triangular u}
\UU \ := \ \UUl \otimes \UUg \ \Big |_{\psi_i^+ \psi_i^- = c^+ c^- = 1} 
\end{equation}
with respect to the pairing of Proposition \ref{prop:pairtor}. Recall from \eqref{eqn:currents2} that $\psi^\pm_i = \psi^\pm_{i,0}$ denote the leading terms of the series $\psi^\pm_i(z)$. We set $\psi_i = \psi_i^+$ and $c = c^+$. \\

\end{definition}

\noindent Unwinding relation \eqref{eqn:drinfeld} in the case at hand, we obtain the following relations in the quantum toroidal algebra, on top of \eqref{eqn:reltor1}, \eqref{eqn:reltor2}, \eqref{eqn:serre}:
\begin{align}
\psi^\mp_j(w)\psi^\pm_i(z) &= \psi^\pm_i(z)\psi^\mp_j(w) \label{eqn:reltor2.5} \\
\psi^\mp_j(w) x^\pm_i(z) &= x^\pm_i(z) \psi^\mp_j(w) \cdot \zeta\left( \frac zw \right)^{\pm 1} \label{eqn:reltor3} \\
\Big[x_i^+(z), x_j^-(w)\Big] &= \delta_j^i \delta\left(\frac zw \right) \cdot \frac {\ph_i^+(z) - \ph_i^-(w)}{q-q^{-1}} \label{eqn:reltor4} 
\end{align}
We write $\UUo \subset \UU$ for the subalgebra generated by the elements $\psi_{i,d}^\pm$. \\

\subsection{}\label{sub:simple}

By analogy with Section 4.17 of \cite{T}, we have an embedding of bialgebras:
$$
\su \hookrightarrow \UU, \qquad \quad x^\pm_i \mapsto x^\pm_{i,0}, \qquad \ph_i \mapsto \frac {\psi_{i+1}}{\psi_i}
$$
where we write $\psi_i = \psi_i^+ = (\psi_i^-)^{-1}$. The above claim also follows from the more general statement that there exists an embedding of bialgebras:
\begin{equation}
\label{eqn:hor}
\uu \hookrightarrow \UU
\end{equation}
It's not immediately obvious how to present \eqref{eqn:hor} in terms of $e_{[i;j)}, f_{[i;j)} \in \uu$. However, the $(a,b)=(1,0)$ case of Lemma \ref{lem:sub} shows how $\uu$ embeds into the double shuffle algebra $\CA$, which by Theorem \ref{thm:iso} is isomorphic to $\UU$. \\ 

\section{The Shuffle Algebra}
\label{sec:shuf}

\subsection{}\label{sub:shufprod}

Inspired by \cite{FO}, we will now present the shuffle algebra. Specifically, we define a realization of the trigonometric version of the shuffle algebra from \loccit Recall the discussion of colored variables from Subsection \ref{sub:colorvar}, and consider an infinite family of variables $z_{i1},z_{i2},...$ of color $i$, for all $1 \leq i \leq n$. We call a rational function:
\begin{equation}
\label{eqn:rat func}
R(...,z_{ia},...)^{1 \leq i \leq n}_{1\leq a \leq k_i} 
\end{equation}
\underline{color-symmetric} if it is symmetric in the variables $z_{i1},...,z_{ik_i}$ for each $i$ separately. The vector $\bk = (k_1,...,k_n) \in \nn$ keeps track of the number of variables of $R$, and will be called the degree of $R$. Often, we will write explicit formulas for rational functions \eqref{eqn:rat func} that include $z_{ia}$ for all $i \in \BZ$, with the convention that:
\begin{equation}
\label{eqn:identify}
z_{ia} \quad \text{ should be replaced with } \quad z_{\bar{i}a} \oq^{-2\left \lfloor \frac {i-1}n \right \rfloor}
\end{equation}
where $\bar{i}$ is the residue of $i$ in the set $\{1,...,n\}$ (see also \eqref{eqn:quasi}). Let $\BF=\BQ(q,\oq)$ and consider the vector space of color-symmetric rational functions:
\begin{equation}
\label{eqn:big}
\CV = \bigoplus_{\bk \in \nn} \BF(...,z_{i1},...,z_{ik_i},...)^{\sym}
_{1 \leq i \leq n} 
\end{equation}
We make the above vector space into a $\BF-$algebra via the \underline{shuffle product}:
\begin{equation}
\label{eqn:mult}
R(...,z_{i1},...,z_{ik_i},...) * R'(...,z_{i1},...,z_{ik'_i},...) = \frac 1{\bk! \cdot \bk'!} \cdot
\end{equation}
$$
\textrm{Sym} \left[ R(...,z_{i1},...,z_{ik_i},...) R'(...,z_{i,k_i+1},...,z_{i,k_i+k'_i},...) \prod_{i,i'=1}^{n} \prod^{j \leq k_i}_{j' > k'_{i'}} \zeta \left( \frac {z_{ij}}{z_{i'j'}} \right) \right] 
$$
for all rational functions $R$ and $R'$ in $\bk$ and $\bk'$ variables, respectively. In \eqref{eqn:mult}, $\sym$ denotes symmetrization with respect to the: 
\begin{equation}
\label{eqn:deffactorial}
(\bk+\bk')! := \prod_{i=1}^{n} (k_i+k_i')!
\end{equation}
permutations that preserve the color of the variables modulo $n$. The algebra $\CV$ is graded by the degree $\bk = (k_1,...,k_n) \in \nn$ encoding the number of variables, and also by the total homogeneous degree $d\in \BZ$ of rational functions $R\in \CV$. We write:
$$
\deg R = \bk \qquad \qquad \homdeg R = d \qquad \qquad \bideg R = (\bk,d)
$$
We will henceforth say that $\CV$ is bigraded by $\nn \times \BZ$. \\

\subsection{}\label{sub:defshuf}

Define the \underline{positive shuffle algebra} $\CA^+$ as the subspace of $\CV$ consisting of rational functions of the form:
\begin{equation}
\label{eqn:shuf}
R(...,z_{i1},...,z_{ik_i},...) = \frac {r(...,z_{i1},...,z_{ik_i},...)}{\prod_{i=1}^{n} \prod_{1\leq b \leq k_{i+1}}^{1\leq a \leq k_{i}} (z_{ia} q - z_{i+1,b} q^{-1})}
\end{equation}
where $r$ is a color-symmetric Laurent polynomial that satisfies the \underline{wheel conditions} below. Specifically, we require that for all $i\in \{1,...,n\}$ we have:
\begin{equation}
\label{eqn:wheel}
r(...,z_{ia},...) \Big |_{z_{i1} \mapsto w, \ z_{i2} \mapsto wq^{\pm 2}, \ z_{i \mp 1,1} \mapsto w} = 0
\end{equation}
Let us make a remark on the denominator in \eqref{eqn:shuf}: if $R$ is a rational function in variables $z_{1a},...,z_{nb}$ for various integers $a,b$, then there will be factors in the denominator of \eqref{eqn:shuf} of the form $z_{na} q - z_{n+1,b} q^{-1} = z_{na} q - z_{1b} q^{-1} \oq^{-2}$, according to \eqref{eqn:identify}. The following Proposition is easy to prove, and we leave it as an exercise to the interested reader (the proof follows Proposition 2.3 of \cite{Shuf} almost word by word): \\

\begin{proposition}

$\CA^+$ is closed under the product \eqref{eqn:mult}, and is thus an algebra. \\

\end{proposition}

\noindent The shuffle algebra $\CA^+$ inherits the grading by $\nn \times \BZ$ from $\CV$, and we will denote the graded pieces by:
$$
\CA^+ = \bigoplus_{\bk \in \nn} \CA_\bk, \qquad \qquad \CA_\bk = \bigoplus_{d\in \BZ} \CA_{\bk,d}
$$
We define the \underline{negative shuffle algebra} as $\CA^- = \left(\CA^+\right)^{\text{op}}$. We will write $R^+$ and $R^-$ for elements of $\CA^+$ and $\CA^-$, although as rational functions they are identical. Define the following grading by $(-\nn) \times \BZ$ on the negative shuffle algebra:
$$
\deg R^- = - \bk \qquad \qquad \homdeg R^- = d \qquad \qquad \bideg R^- = (- \bk,d)
$$
The graded pieces of $\CA^-$ will be denoted by $\CA_{-\bk}$ and $\CA_{-\bk,d}$. The assignment:
\begin{equation}
\label{eqn:ordern}
\CA^\pm \longrightarrow \CA^\pm \qquad \qquad R(...,z_{ia},...) \longrightarrow R\left(...,z_{i+1,a} \oq^{\frac 2n},...\right)
\end{equation}
is an order $n$ algebra automorphism, although to be precise, one needs to adjoin an $n-$th root of $\oq$ to the ground field for it to be well-defined. In particular, this automorphism preserves rational functions of the form \eqref{eqn:shuf} and \eqref{eqn:wheel}, due to condition \eqref{eqn:identify}. We think of \eqref{eqn:ordern} as an action $\BZ/n\BZ \curvearrowright \CA^\pm$ by algebra automorphisms. \\

\begin{remark}

Our $\CA^+$ differs from other shuffle algebras incarnations of $\UU$ in the literature by a simple change of notation. For example, to go from our conventions to those of \cite{FT}, Section 1.7, one needs to consider:
$$
d = q^{-1} \oq^{-\frac 2n}
$$
and it is straightforward to see that:
$$
R(...,z_{ia},...) \mapsto  R(...,z_{ia},...) \prod_{i=1}^{n} \left[ q^{-\frac {k_i^2}2} \prod_{1\leq b \leq k_{i+1}}^{1\leq a \leq k_{i}} \frac {z_{ia} q^2 - z_{i+1,b}}{z_{ia} \oq^{-\frac 2n} - z_{i+1,b}} \right] \Big|_{z_{ia} \mapsto x_{ia} \oq^{-\frac {2i}n}}
$$
gives an isomorphism between $\CA^+$ and the shuffle algebra of \loccitt. \\

\end{remark}

\subsection{}\label{sub:coproduct}

Define the \underline{extended shuffle algebras} as:
\begin{align*}
&\CA^\geq = \Big \langle \CA^+ , \psi_i^+, \psi^+_{i,1}, \psi^+_{i,2}, ... , c^+ \Big \rangle_{1 \leq i \leq n} \Big/^{c \text{ central}}_{[\psi_{i,d}^+, \psi_{i',d'}^+] = 0 \text{ and relation \eqref{eqn:comm0} for }\pm = +} \\
&\CA^\leq = \Big \langle \CA^- , \psi_i^-, \psi^-_{i,1}, \psi^-_{i,2}, ..., c^- \Big \rangle_{1 \leq i \leq n} \Big/^{c \text{ central}}_{[\psi_{i,d}^-, \psi_{i',d'}^-] = 0 \text{ and relation \eqref{eqn:comm0} for }\pm = -}
\end{align*}
We assume $c^\pm$ and $\psi_i^\pm = \psi_{i,0}^\pm$ invertible, and include their inverses into the algebras above. As before, consider the currents:
$$
\psi^\pm_i(z) = \sum_{d = 0}^\infty \frac {\psi^\pm_{i,d}}{z^{\pm d}} 
$$
for all $1 \leq i \leq n$, and extend the index set to all $i \in \BZ$ via $\psi^\pm_{i}(z) = c^{\pm 1} \psi^\pm_{i-n}(z\oq^2)$. Finally, we impose the following commutation relations in $\CA^\leq$ and $\CA^\geq$: 
\begin{equation}
\label{eqn:comm0}
\psi_j^\pm(w) * R^\pm = \left[ R^\pm(...,z_{ia},...) \prod_{i=1}^{n} \prod_{a=1}^{k_{i}} \zeta \left( \frac {z_{ia}}w \right)^{\pm 1} \right] * \psi_j^\pm(w) 
\end{equation}
for any $R^\pm \in \CA^\pm$, where the $\zeta$ factor is expanded in non-negative powers of $w^{\mp 1}$. For the purpose of defining $\zeta$ in \eqref{eqn:comm0}, we think of $\col z_{ia} = i$ and $\col w = j$ and use \eqref{eqn:bigzeta}. One of the main reasons for defining the extended shuffle algebras is that they admit coproducts: \\

\begin{proposition}
\label{prop:coproduct}

There are bialgebra structures on $\CA^\leq$ and $\CA^\geq$, with coproduct:
\begin{equation}
\label{eqn:coproduct0}
\Delta(\psi^\pm_i(z)) = \psi^\pm_i(z) \otimes \psi^\pm_i(z)
\end{equation}
while for all $R^\pm \in \CA_{\pm \bk}$, we have:
\begin{align}
&\Delta(R^+) = \sum_{\bl \in \nn}^{\bl\leq \bk}  \frac {\left[ \prod_{1 \leq i \leq n}^{a>l_i} \ph^+_i(z_{ia}) \right] * R^+(z_{i, a \leq l_i} \otimes z_{i, a>l_i})}{\prod_{1 \leq i' \leq n}^{1 \leq i \leq n} \prod^{a \leq l_{i}}_{a' > l_{i'}} \zeta(z_{i'a'}/z_{ia})} \label{eqn:coproduct1} \\
&\Delta(R^-) = \sum_{\bl \in \nn}^{\bl\leq \bk}  \frac {R^-(z_{i, a \leq l_i} \otimes z_{i, a>l_i}) * \left[ \prod_{1 \leq i \leq n}^{a \leq l_i} \ph^-_i(z_{ia}) \right]}{\prod_{1 \leq i' \leq n}^{1 \leq i \leq n} \prod^{a \leq l_{i}}_{a' > l_{i'}} \zeta(z_{ia}/z_{i'a'})} \label{eqn:coproduct2}
\end{align}
where we set $\ph^\pm_i(z) = \psi_{i+1}^\pm(zq^2)/\psi_i^\pm(z)$, in accordance with \eqref{eqn:ph}. \\
\end{proposition}

\begin{remark}

To think of \eqref{eqn:coproduct1} as a well-defined tensor, we expand the right-hand side in non-negative powers of $z_{ia} / z_{i'a'}$ for $a\leq l_i$ and $a'>l_{i'}$, thus obtaining an infinite sum of monomials. In each of these monomials, we put the symbols $\ph^+_{i,d}$ to the very left of the expression, then all powers of $z_{ia}$ with $a\leq l_i$, then the $\otimes$ sign, and finally all powers of $z_{ia}$ with $a>l_i$. The resulting expression will be a power series, and therefore lies in a completion of $\CA^\geq \otimes \CA^\geq$. The same argument applies to \eqref{eqn:coproduct2}, still using non-negative powers of $z_{ia}/z_{i'a'}$ for $a\leq l_i$ and $a'>l_{i'}$, and keeping all the $\ph_{i,d}^-$ to the very right. Proposition \ref{prop:coproduct} is proved exactly like Proposition 4.1 of \cite{Shuf}, so we will leave it as an exercise to the interested reader. \\

\end{remark}

\subsection{}\label{sub:full}

Let us consider the following modification of the rational function \eqref{eqn:bigzeta}:
$$
\zeta_p \left( \frac zw \right) = \zeta \left( \frac zw \right) \cdot \left( \frac {z p \oq^{2 \frac {i-j}n } - w p^{-1}}{z q \oq^{2 \frac {i-j}n} - w q^{-1}} \right)^{\delta_{\bari}^{\barj}}
$$
for variables $z$, $w$ of any colors $i,j \in \BZ$. While $q,\oq$ are formal parameters, we will assume them to be complex numbers in what follows, for the purpose of defining contour integrals. The notation:
\begin{equation}
\label{eqn:assumption int}
\int^*_{|z| = r} f(z) Dz \qquad \text{where } Dz = \frac {dz}{2\pi i z} 
\end{equation} 
refers to the sum of all residues of $f(z)$ inside the circle of radius $r>0$ around the origin. The issue of which poles fall inside this circle is determined by the conditions $|q| > 1 > |\oq| > |p|$ and $|pq|=1$. The following will be proved in the Appendix: \\

\begin{proposition}
\label{prop:pairshuf}

There exists a bialgebra pairing $\CA^\leq \otimes \CA^\geq \rightarrow \BF$ given by:
\begin{equation}
\label{eqn:pairshuf0}
\left \langle \psi^-_j(w),\psi^+_i(z) \right \rangle = \text{right-hand side of \eqref{eqn:pairtor2}}
\end{equation}
while:
\begin{equation}
\label{eqn:pairshuf}
\left \langle R^-,R^+ \right \rangle = \frac {(1-q^{-2})^{|\bk|}}{\bk!} \int^*_{|z_{ia}|=\oq^{- \frac {2i}n}} \frac {R^-(...,z_{ia},...)R^+(...,z_{ia},...)}{\prod_{i,j=1}^{n} \prod^{(i,a) \neq (j,b)}_{a\leq k_i, b \leq k_j} \zeta_p(z_{ia}/z_{jb})} \prod^{1 \leq i \leq n}_{1 \leq a \leq k_i} Dz_{ia} \Big |^{p \mapsto q} \quad
\end{equation}
for any $R^\pm \in \CA_{\pm \bk}$. The notation in the right-hand side should be read as follows: ``compute the integral by residues assuming $|q|>1>|\oq| > |p|$ and $|pq|=1$, obtain an answer which is a rational function of $p$ and $q$, and then set $p = q$ in the answer". \\

\end{proposition}

\noindent The datum provided by Proposition \ref{prop:pairshuf} allows one to construct the Drinfeld double:
\begin{equation}
\label{eqn:triangular a}
\CA \ := \ \CA^\leq \otimes \CA^\geq \ \Big|_{\psi_i^+\psi_i^- = c^+ c^- = 1} 
\end{equation}
We will write $\Ao \subset \A$ for the Cartan subalgebra generated by the elements $\psi_{i,d}^\pm$. \\

\begin{proof} \emph{of Theorem \ref{thm:iso}:} We claim the existence of a bialgebra isomorphism:
\begin{equation}
\label{eqn:upsy}
\Upsilon:\UU \stackrel{\sim}\longrightarrow \A
\end{equation}
generated by the assignments $\Upsilon(c) = c$ and:
\begin{align}
&\Upsilon^+:\UUp \longrightarrow \Ap, \qquad \qquad \Upsilon^+\left(x^+_{i,d}\right)  = \frac {z_{i1}^d}{q^{-1} - q} \label{eqn:ass1} \\
&\Upsilon^0:\UUo \longrightarrow \Ao, \ \qquad \qquad \Upsilon^0\left(\psi^\pm_{i,d}\right) \ =  \psi^\pm_{i,d} \label{eqn:ass2} \\
&\Upsilon^-:\UUm \longrightarrow \Am, \qquad \qquad \Upsilon^-\left(x^-_{i,d}\right)  = \frac {z_{i1}^d}{1 - q^{-2}} \label{eqn:ass3} 
\end{align}
By analogy with \cite{E}, one proves that the assignments \eqref{eqn:ass1} and \eqref{eqn:ass3} extend to algebra homomorphisms, and the fact that they are compatible with \eqref{eqn:ass2} is immediate from our choice of relations. We conclude that we have homomorphisms:
\begin{equation}
\label{eqn:ass4}
\Upsilon^\geq:\UUg \longrightarrow \Ag \qquad \qquad \Upsilon^\leq : \UUl \longrightarrow \Al
\end{equation}
obtained by merging \eqref{eqn:ass1} with \eqref{eqn:ass2} and \eqref{eqn:ass2} with \eqref{eqn:ass3}, respectively. It is easy to check that the homomorphisms \eqref{eqn:ass4} match the coproduct and the pairing on the two sides, by comparing \eqref{eqn:deltator1}--\eqref{eqn:pairtor2} with \eqref{eqn:coproduct0}--\eqref{eqn:pairshuf} for $R^\pm = z_{i1}^d$. Since the Drinfeld double is uniquely determined from the coproduct and pairing, $\Upsilon$ extends to a bialgebra homomorphism \eqref{eqn:upsy} between the Drinfeld doubles. \\

\noindent Let us now prove that $\Upsilon$ is injective. There exist representations called Fock spaces, indexed by various $u \in \mathbb{C}$ and $\alpha \in \BZ/n\BZ$, and one can form their tensor products:
\begin{equation}
\label{eqn:two actions}
\UU \curvearrowright \bigotimes_i F^{\alpha_i}(u_i) \curvearrowleft \CA 
\end{equation}
where the action on the left is treated in \cite{BT} and the action on the right is treated in \cite{thesis}. Proposition IV.8. of \loccit deals with the incarnation of tensor products of Fock spaces as $K$--theory groups of Nakajima quiver varieties for the cyclic quiver, but it is elementary to check that the two actions in \eqref{eqn:two actions} are compatible under the homomorphism $\Upsilon$. Since the direct sum of the actions on the left of \eqref{eqn:two actions} is faithful as a $\UU$--module (Corollary 2.15 of \cite{BT}), this implies the map $\Upsilon$ is injective. \\

\noindent Therefore, all that remains in order to prove Theorem \ref{thm:iso} is the surjectivity of $\Upsilon$. This boils down to the following result, to be proved at the end of this Section. \\

\begin{proposition}
\label{prop:surj}
	
The homomorphisms $\Ups^\leq$ and $\Ups^\geq$ are surjective. 
	
\end{proposition} \end{proof}

\subsection{}\label{sub:slope}

In the remainder of this Section, we introduce technology which will allow us to prove Proposition \ref{prop:surj}. The main idea is to estimate the dimensions of the infinite-dimensional graded components $\CA_{\pm \bk,d}$ by introducing a suitable filtration. \\

\begin{definition}
\label{def:slope}	
	
For any shuffle elements $R^\pm\in \CA_{\pm \bk}$ and any $\mu \in \BR$, if the limit:
\begin{equation}
\label{eqn:limit1}
\lim_{\xi \rightarrow \infty} \frac {R^+(..., \xi z_{i1},..., \xi z_{i l_i}, z_{i, l_i+1},..., z_{ik_i},...)}{\xi^{\mu |\bl|}}
\end{equation}
or:
\begin{equation}
\label{eqn:limit2}
\lim_{\xi \rightarrow 0} \frac {R^-(..., \xi z_{i1},..., \xi z_{il_i}, z_{i, l_i+1},..., z_{ik_i},...)}{\xi^{-\mu |\bl|}} 
\end{equation}
exists and is finite for all $0\leq \bl \leq \bk$, then we say that $R^\pm$ has \underline{slope} $\leq \mu$. \\
	
\end{definition}

\noindent In \eqref{eqn:limit1} and \eqref{eqn:limit2}, we have $|\bl| = l_1+...+l_n$. We will write:
$$
\CA_{\leq \mu | \pm \bk} \subset \CA_{\pm \bk}
$$ 
for the subspace of shuffle elements of slope $\leq \mu$. By taking the case $\bl = \bk$ in \eqref{eqn:limit1} and \eqref{eqn:limit2}, respectively, we see that $R^\pm$ can only have slope $\leq \mu$ if:
\begin{equation}
\label{eqn:ineq}
\pm \homdeg R^\pm \leq \mu |\bk|
\end{equation}
Therefore, the bigraded pieces:
$$
\CA_{\leq \mu|\pm \bk, d} \ := \ \CA_{\leq \mu| \pm \bk} \cap \CA_{\pm \bk,d} 
$$
are non-zero only if $\pm d \leq \mu|\bk|$. The following Proposition is a simple exercise, based on the fact that $\lim_{x \rightarrow 0 \text{ or } \infty}\zeta(x) \in \{q,1,q^{-1}\}$. We leave its proof to the interested reader, and note that it is analogous to Proposition 2.3 of \cite{Shuf}. \\

\begin{proposition}
\label{prop:qarth}
	
For either choice of the sign $\pm$, the vector space:
$$
\CA^\pm_{\leq \mu} \ := \ \bigoplus_{\bk\in \nn} \CA_{\leq \mu| \pm \bk } \ \subset \ \CA^\pm
$$
is a subalgebra, i.e. the property of having slope $\leq \mu$ is preserved under $*$ of \eqref{eqn:mult}. \\
	
\end{proposition}

\subsection{}\label{sub:subalgebra}

We will say that a shuffle element has slope $<\mu$ if it has slope $\leq \mu-\e$ for a positive number $\e$. By the very definition of the coproduct in \eqref{eqn:coproduct1} and \eqref{eqn:coproduct2}, and the slope in \eqref{eqn:limit1} and \eqref{eqn:limit2}, we have:
\begin{equation}
\label{eqn:radu1}
\Delta(R^+) = \Delta_\mu(R^+) + (\textrm{anything}) \otimes (\text{slope} < \mu) 
\end{equation}
\begin{equation}
\label{eqn:radu2}
\Delta(R^-) = \Delta_\mu(R^-) + (\text{slope} < \mu) \otimes (\textrm{anything})
\end{equation}
\footnote{Formula \eqref{eqn:radu1} says that in $\Delta(R^+) = R_1^+ \otimes R_2^+$, every term of the form $R_2^+$ has slope $\leq \mu$. To see this, one needs to show that the second tensor factor of $\Delta(R_2^+)$ has \underline{naive slope}, i.e. the ratio of its total homogeneous degree by its total number of variables, $\leq \mu$. Since $R^+$ has slope $\leq \mu$, then the third tensor factor of $\Delta^{(2)}(R^+) = R_1^+ \otimes R_2^+ \otimes R_3^+$ always has naive slope $\leq \mu$. But the coassociativity of $\Delta$ implies that all second tensor factors of $\Delta(R_2^+)$'s are among the $R_3^+$'s.} for any $R^\pm \in \CA^\pm_{\leq \mu}$, where the \underline{leading terms} $\Delta_\mu$ are defined by:
\begin{align}
&\Delta_\mu(R^+) = \sum^{\bl \in \nn}_{\bl\leq \bk} \ph_{\bk-\bl} * \lim_{\xi \rightarrow \infty} \frac {R^+(z_{i,a\leq l_i} \otimes \xi\cdot z_{i,a>l_i})}{\xi^{\mu|\bk - \bl|} q^{\langle \bk - \bl, \bl \rangle}} \label{eqn:cristian1} \\
&\Delta_\mu(R^-) = \sum^{\bl \in \nn}_{\bl\leq \bk} \lim_{\xi \rightarrow 0} \frac {R^-(\xi \cdot z_{i,a\leq l_i} \otimes z_{i,a>l_i})}{\xi^{-\mu|\bl|} q^{- \langle \bl, \bk - \bl \rangle}} * \ph_{-\bl} \label{eqn:cristian2}
\end{align}
In the above, we write $\ph_{\pm \bl} = \prod_{i=1}^{n} \ph_i^{\pm l_i}$ for any $\bl \in \nn$. Recall that the tensor product inside the rational function $R^\pm$ means that all powers of $z_{i,a\leq l_i}$ go to the left of the $\otimes$ sign, while all powers of $z_{i,a>l_i}$ go to the right. The powers of $q$ in the denominators of formulas \eqref{eqn:cristian1} and \eqref{eqn:cristian2} arise from the fact that:
$$
\lim_{\xi \rightarrow \infty} \zeta \left(\frac {\xi x}y \right) = q^{\langle \col x, \col y \rangle} \qquad \text{and} \qquad \lim_{\xi \rightarrow 0} \zeta \left( \frac {\xi x}y \right) = q^{- \langle \col x, \col y \rangle}
$$
Let us now consider those elements of slope $\leq \mu$ for which the inequality \eqref{eqn:ineq} becomes an equality, and define the vector spaces:
$$
\CB^\pm_\mu \ = \ \bigoplus^{\mu|\bk| \in \BZ}_{\bk \in \nn} \CB_{\mu|\pm \bk} \ := \ \bigoplus^{\mu|\bk| = d}_{\bk \in \nn} \CA_{\leq \mu|\pm \bk, \pm d} \ \subset \ \CA^\pm
$$
By analogy with Proposition \ref{prop:qarth}, one shows that $\CB^\pm_\mu$ are subalgebras for any $\mu$. This is clear, since both having slope $\leq \mu$ and equality in \eqref{eqn:ineq} are properties preserved by the shuffle product. We may consider the extended subalgebras:
\begin{align}
&\CB^\geq_\mu = \Big \langle \CB^+_\mu, \psi_1,...,\psi_n, c \Big \rangle \ \qquad \subset \CA^\geq  \label{eqn:sub1} \\
&\CB^\leq_\mu = \Big \langle \CB^-_\mu, \psi^{- 1}_1,...,\psi^{- 1}_n, c^{- 1} \Big \rangle \subset \CA^\leq \label{eqn:sub2}
\end{align}
Using \eqref{eqn:cristian1} and \eqref{eqn:cristian2}, it is easy to see that $\Delta_\mu$ is a coproduct on the subalgebras \eqref{eqn:sub1} and \eqref{eqn:sub2}. Moreover, by analogy with Proposition V.3 of \cite{thesis}, the bialgebra pairing between $\CA^\leq$ and $\CA^\geq$ of Proposition \ref{prop:pairshuf} descends to a bialgebra pairing:
$$
\langle \cdot , \cdot \rangle \ : \ \CB^\leq_\mu \otimes \CB^\geq_\mu \longrightarrow \BF
$$
and the Drinfeld double $\CB_\mu := \CB^\leq_\mu \otimes \CB^\geq_\mu$ with respect to the above data embeds:
\begin{equation}
\label{eqn:subalgebra}
\CB_\mu \ \subset \ \CA
\end{equation}
into the double shuffle algebra. One of the main results of the current Section is: \\

\begin{lemma}
\label{lem:sub}
	
For any $\mu = \frac ba \in \BQ$ with $a$ and $b$ coprime, we have an isomorphism:
\begin{equation}
\label{eqn:smalliso}
\Xi: U_{q}(\dot{\fgl}_{\frac ng})^{\otimes g} \stackrel{\sim}\longrightarrow \CB_{\mu}
\end{equation}
where $g = \gcd(n,a)$. The isomorphism $\Xi$ preserves the bialgebra structures. \\
	
\end{lemma}

\subsection{}\label{sub:bound}

The isomorphism $\Xi$ will be constructed by exhibiting bialgebra morphisms:
\begin{align}
&\Xi^\geq : U^\geq_{q}(\dot{\fgl}_{\frac ng})^{\otimes g} \longrightarrow \CB^\geq_{\mu} \label{eqn:small1} \\
&\Xi^\leq : U^\leq_{q}(\dot{\fgl}_{\frac ng})^{\otimes g} \longrightarrow \CB^\leq_{\mu} \label{eqn:small2}
\end{align}
which preserve the bialgebra pairings, i.e. $\langle \Xi^\leq(x), \Xi^\geq(y) \rangle = \langle x,y \rangle$. This will imply the existence of a homomorphism \eqref{eqn:smalliso}, and then we must prove that this homomorphism is injective and surjective. To this end, we will need to estimate the dimensions of the graded pieces of the algebras $\CB^\geq_\mu$ and $\CB^\leq_\mu$. We will actually prove the following more general result in the Appendix: \\

\begin{lemma}
\label{lem:min}
	
For any $\mu\in \BR$, $\bk \in \nn$ and $d\in \BZ$, the dimension of $\CA_{\leq \mu| \pm \bk,\pm d}$ is at most equal to the number of unordered collections:
\begin{equation}
\label{eqn:collection}
\CC = \{(i_a,j_a,d_a)\}_{a\in \{1,...,t\}},
\end{equation}
such that the following conditions hold:
\begin{align}
&\sum_{a=1}^t [i_a;j_a) = \bk, \qquad \ \sum_{a=1}^t d_a = d \label{eqn:cond1} \\
&d_a \leq \mu (j_a-i_a), \quad \forall a\in \{1,...,t\} \label{eqn:cond2}
\end{align}
Recall that $[i;j) \in \nn$ was defined in \eqref{eqn:part}, and we identify $[i;j) = [i-n;j-n)$. \\
	
\end{lemma}

\noindent As a corollary of Lemma \ref{lem:min}, let us estimate the dimension of $\CB_{\mu|\pm \bk} = \CA_{\leq \mu|\pm \bk, \pm d}$ in the case $d = \mu|\bk|$. In this situation, \eqref{eqn:cond1} forces equality in \eqref{eqn:cond2}, so the numbers $d_1,...,d_t$ are completely determined from the intervals $[i_1;j_1),...,[i_t;j_t)$ via $d_a = \mu(j_a-i_a)$. One requires $d_a \in \BZ$ for all $a \in \{1,...,t\}$, so we conclude that: 
\begin{equation}
\label{eqn:bound}
\dim \CB_{\mu|\pm \bk} \leq \# \ \Big \{\mu-\text{integral partitions }C \vdash \bk \Big \}
\end{equation}
where a partition: 
\begin{equation}
\label{eqn:partition}
C \ = \ \Big \{[i_1;j_1),...,[i_t;j_t) \Big \}
\end{equation}
of $\bk \in \nn$ is called $\mu-$\underline{integral} if:
\begin{equation}
\label{eqn:integral}
\mu(j_a-i_a) \in \BZ, \quad \forall a\in \{1,...,t\}
\end{equation}

\subsection{}\label{sub:alphabeta}

For any arc $[i;j)$, define the following linear maps:
\begin{align}
&\alpha_{[i;j)} : \CA_{[i;j)} \longrightarrow \BF, \ \qquad \alpha_{[i;j)}(R) = \frac {R(1_{j-1},...,1_i)}{\prod_{i \leq a < b < j} \zeta(1_{b-a})} \label{eqn:alpha} \\
&\beta_{[i;j)} : \CA_{-[i;j)} \longrightarrow \BF, \qquad \beta_{[i;j)}(R) = \frac {R(q^{2i},...,q^{2j-2})}{\prod_{i \leq a < b < j} \zeta(q^{2a-2b})} \label{eqn:beta}
\end{align}
In the right-hand side of \eqref{eqn:alpha}, we write $1_a$ for the number 1 plugged into a variable of color $a$. In the right-hand side of \eqref{eqn:beta}, each constant $q^{2a}$ is plugged into a variable of color $a
$, and we note that one needs to cancel the poles of $\zeta(q^{2a-2b})$ against the poles of the evaluation $R(q^{2i},...,q^{2j-2})$ in order for the fraction to be well-defined. Refining the proof of Lemma \ref{lem:min} will allow us to prove: \\

\begin{proposition}
\label{prop:sam}

A shuffle element $R \in \CB_{\mu|\bk}$ is 0 if and only if:
\begin{equation}
\label{eqn:sam}
\alpha_{[i_1;j_1)} \otimes ... \otimes \alpha_{[i_t;j_t)} \left(\textrm{component of }\Delta_\mu^{(t-1)}(R)\textrm{ in } \CB_{\mu| [i_1;j_1)} \otimes ... \otimes \CB_{\mu| [i_t;j_t)}\right) = 0 \quad
\end{equation}
while a shuffle element $R \in \CB_{\mu|-\bk}$ is 0 if and only if:
\begin{equation}
\label{eqn:adams}
\beta_{[i_1;j_1)} \otimes ... \otimes \beta_{[i_t;j_t)} \left(\textrm{component of }\Delta_\mu^{(t-1)}(R)\textrm{ in } \CB_{\mu| - [i_1;j_1)} \otimes ... \otimes \CB_{\mu| - [i_t;j_t)}\right) = 0 \quad
\end{equation}
for any $\mu-$integral partition $C = \{[i_1;j_1),...,[i_t;j_t)\} \vdash \bk$. \\

\end{proposition}

\noindent The fact that we use the linear maps $\alpha$ to describe positive shuffle elements and the linear maps $\beta$ to describe negative shuffle elements is simply a choice we make to ensure consistency in the remainder of this Section. We could have effortlessly used either $\alpha$ or $\beta$ to describe both positive and negative shuffle elements. \\

\begin{lemma}
\label{lem:pseudo}

The maps $\alpha_{[i;j)}$ are \underline{pseudo-multiplicative}, in the sense that:
$$
\alpha_{[i;j)}(R_1 * R_2) = \begin{cases} \alpha_{[a;j)}(R_1) \alpha_{[i;a)}(R_2) & \text{if } \exists a \text{ s.t. } \deg R_1 = [a;j) \text{ and } \deg R_2 = [i;a) \\ 
0 & \text{otherwise} \end{cases}
$$
The same statement holds for the linear maps $\beta_{[i;j)}$. \\
\end{lemma}

\noindent It is straightforward to prove Lemma \ref{lem:pseudo} directly from the definition \eqref{eqn:alpha}, so we leave it as an exercise. It is proved along the lines of Exercise V.9 of \cite{thesis}, and it essentially follows from the fact that $\zeta(z)|_{z \mapsto 1} = 0$ whenever $\col z = - 1$. \\ 

\subsection{}\label{sub:min}

We will now construct explicit shuffle elements, which will allow us to interpret properties \eqref{eqn:sam} and \eqref{eqn:adams} via the bialgebra pairing \eqref{eqn:pairshuf}. These shuffle elements, which will be interpreted geometrically in \cite{Aff}, are:
\begin{align}
&S^\pm_m(z_i,...,z_{j-1}) = \sym \left[ \frac {m(z_i,...,z_{j-1})}{\left(1 - \frac {z_{i+1}}{z_{i}q^2}  \right) ... \left(1 - \frac {z_{j-1}}{z_{j-2}q^2} \right)} \prod_{i\leq a < b < j} \zeta \left( \frac {z_b}{z_a} \right)  \right] \in \CA^\pm  \label{eqn:s} \\ &T^\pm_m(z_i,...,z_{j-1}) = \sym \left[ \frac {m(z_i,...,z_{j-1})}{\left(1 - \frac {z_{i}}{z_{i+1}}  \right) ... \left(1 - \frac {z_{j-2}}{z_{j-1}} \right)} \prod_{i\leq a < b < j} \zeta \left( \frac {z_a}{z_b} \right)  \right] \in \CA^\pm \label{eqn:t} 
\end{align}
for any $i < j$ and any Laurent polynomial $m(z_i,...,z_{j-1})$. Let us explain how to think of the right-hand sides as shuffle elements. We regard each $z_a$ as a variable of color $a$, for all $a \in \{i,...,j-1\}$, and relabel them according to \eqref{eqn:identify}:
$$
z_a,z_{a+n},z_{a+2n}, ... \quad \leadsto \quad z_{a1}, z_{a2} \oq^{-2}, z_{a3} \oq^{-4},...
$$
for each $1 \leq a \leq n$. Then $S_m^\pm$ and $T_m^\pm$ are manifestly elements of $\CV$ from \eqref{eqn:big}. The following stronger result will be proved in the Appendix: \\

\begin{proposition}
\label{prop:belong}
	
For any $i < j$ and Laurent polynomial $m(z_i,...,z_{j-1})$, we have: 
$$
S_m^\pm, \ T_m^\pm \ \in \ \emph{Im }\Upsilon^\pm
$$
where $\Upsilon^\pm$ denote the homomorphisms of \eqref{eqn:ass1} and \eqref{eqn:ass3}. \\
	
\end{proposition}

\subsection{}\label{sub:george} 

We will now construct particular cases of the shuffle elements $S^+_m$ and $T^-_m$, which will correspond to the root generators $e_{[i;j)}$ and $f_{[i;j)}$ under the isomorphisms of Lemma \ref{lem:sub}. The following results will be proved in the Appendix: \\

\begin{proposition}
\label{prop:tommy}
	
For any $\mu \in \BQ$ and any $\mu-$integral arc $[i;j)$, the elements:
\begin{align}
&E_{[i;j)}^{\mu} = \emph{Sym} \left[ \frac {\prod_{a=i}^{j-1} z_a^{\lfloor \mu(a-i+1) \rfloor - \lfloor \mu(a-i) \rfloor}}{\left(1 - \frac {z_{i+1}}{z_{i}q^2}  \right) ... \left(1 - \frac {z_{j-1}}{z_{j-2}q^2} \right)} \prod_{i\leq a < b < j} \zeta \left( \frac {z_b}{z_a} \right)  \right] \in \CA^+  \label{eqn:min1} \\
&F_{[i;j)}^{\mu} = \emph{Sym} \left[ \frac {\prod_{a=i}^{j-1} z_a^{\lfloor \mu(a-i) \rfloor - \lfloor \mu(a-i+1) \rfloor}}{\left(1 - \frac {z_{i}}{z_{i+1}}  \right) ... \left(1 - \frac {z_{j-2}}{z_{j-1}} \right)} \prod_{i\leq a < b < j} \zeta \left( \frac {z_a}{z_b} \right)  \right] \in \CA^- \label{eqn:min2}
\end{align}
lie inside $\CB^\pm_{\mu}$. Moreover, we have: 
\begin{align}
&\Delta_\mu\left(E_{[i;j)}^\mu\right) = \sum_{a=i}^j E_{[a;j)}^\mu * \frac {\psi_a}{\psi_i} \otimes E_{[i;a)}^\mu \label{eqn:cop1} \\ 
&\Delta_\mu\left(F_{[i;j)}^\mu\right) = \sum_{a=i}^j F_{[i;a)}^\mu \otimes F_{[a;j)}^\mu * \frac {\psi_i}{\psi_a} \label{eqn:cop2}
\end{align}
where we set $E^\mu_{[i;i)} = F^\mu_{[i;i)} = 1$, and $E_{[i;j)}^\mu = F_{[i;j)}^\mu = 0$ if $[i;j)$ is not $\mu-$integral. \\
\end{proposition}

\begin{proposition}
\label{prop:pairing}
	
For any $R^\pm \in \CB^\pm_\mu$ and any $i<j$, we have:
\begin{align}
&\left \langle F_{[i;j)}^\mu, R^+ \right \rangle = (1-q^{-2})^{j-i} \cdot \alpha_{[i;j)}(R^+) \label{eqn:bonnie1} \\
&\left \langle R^-, E_{[i;j)}^\mu \right \rangle = (1 -q^{-2})^{j-i} \gamma \cdot  \beta_{[i;j)}(R^-)  \label{eqn:bonnie2}
\end{align}
where $\gamma = \prod_{a=i}^{j-1} q^{2a \left (\lfloor \mu(a-i+1) \rfloor - \lfloor \mu(a-i) \rfloor \right)}$. As a corollary, we obtain:
\begin{equation}
\label{eqn:tyler}
\left \langle E_{[i;j)}^\mu, F_{[i';j')}^\mu \right \rangle = \delta_{[i';j')}^{[i;j)} \left( 1 - q^{-2} \right)
\end{equation}
Any element of $\CB_\mu^-$ (resp. $\CB_\mu^+$) which pairs trivially with all products of $E_{[i;j)}^\mu$ (resp. $F_{[i;j)}^\mu$) is 0, hence the bialgebra pairing is non-degenerate between $\CB_\mu^-$ and $\CB_\mu^+$. \\
	
\end{proposition} 

\noindent The Proposition above gives us a quick proof of Lemma \ref{lem:pseudo} when $R_1, R_2 \in \CB^\pm_\mu$:
\begin{multline*}
\alpha_{[i;j)}(R_1 * R_2) = \frac {\left \langle F_{[i;j)}^\mu, R_1 * R_2 \right \rangle}{(1-q^{-2})^{j-i}} \stackrel{\eqref{eqn:bialg 2}}= \frac {\left \langle \Delta_\mu^{\text{op}} \left( F_{[i;j)}^\mu \right), R_1 \otimes R_2 \right \rangle}{(1-q^{-2})^{j-i}} = \\ = \sum_{a=i}^j \frac {\left \langle F_{[a;j)}^\mu * \frac {\psi_i}{\psi_a}, R_1 \right \rangle \left \langle F_{[i;a)}^\mu , R_2 \right \rangle}{(1-q^{-2})^{a-i}(1-q^{-2})^{j-a}} = \sum_{a=i}^j \delta_{\deg R_2}^{[i;a)} \alpha_{[a;j)}(R_1)\alpha_{[i;a)}(R_2)
\end{multline*}
and similarly for $\beta_{[i;j)}$. \\

\subsection{}\label{sub:g}

In the special case $\mu = 0$, we will construct another important class of shuffle elements $\in \CB_0^\pm$. They have degrees $k \bde = (k,...,k)$ for any $k\in \BN$, and are given by:
\begin{align}
&G_{k} = \frac {\oq^{k^2}}{(q^{-1}-q)^{nk}} \prod_{i=1}^n \prod_{1\leq a,b \leq k} \frac {z_{ia} q - z_{ib}q^{-1}}{z_{i-1,a} q - z_{ib}q^{-1}} \in \CA_{k \bde} \label{eqn:g1} \\
&G_{-k} = \frac {\oq^{k^2}}{(1-q^{-2})^{nk}} \prod_{i=1}^n \prod_{1\leq a,b \leq k} \frac {z_{ia} q - z_{ib}q^{-1}}{z_{i-1,a} q - z_{ib}q^{-1}} \in \CA_{- k \bde} \label{eqn:g2}
\end{align}
The fact that the above rational function satisfies the wheel conditions \eqref{eqn:wheel} is immediate, as is the fact that it has homogeneous degree 0. To show that $G_{\pm k}$ has slope $\leq 0$, we need to compute its degree in any subset of $\ba \leq k \bde$ variables:
$$
\deg_{...,z_{i1},...,z_{ia_i},...} G_{\pm k} = \sum_{i=1}^n \Big[ a_i^2 + 2a_i(k-a_i) - a_i a_{i+1} - a_i(k-a_{i-1}) - a_{i-1}(k-a_i) \Big] = 
$$
\begin{equation}
\label{eqn:gm}
= \sum_{i=1}^n \left( a_i a_{i-1} - a_i^2 \right) = - \frac 12 \sum_{i=1}^n (a_i - a_{i-1})^2 \leq 0
\end{equation}
We thus conclude that $G_{\pm k} \in \CB^\pm_0$. Moreover, we have equality in \eqref{eqn:gm} if and only if $a_0 = ... = a_{n-1} = a$ for some natural number $a$. Therefore:
\begin{align}
&\Delta_0 \left( G_{k} \right) \ \ = \ \ \sum^{a,b \geq 0}_{a+b = k} G_a c^b \otimes G_b \label{eqn:nicecop1} \\
&\Delta_0 \left( G_{-k} \right) = \sum^{a,b \geq 0}_{a+b = k} G_{-a} \otimes G_{-b} c^{-a} \label{eqn:nicecop2}
\end{align}
where $c = \ph_1...\ph_n$ is the central element. Moreover, it is clear that $G_{\pm k}$ is invariant under the action $\BZ/n\BZ \curvearrowright \CA^\pm$. The shuffle elements $G_{\pm k}$ (although under a different normalization) have been studied by \cite{FT}, who proved relations \eqref{eqn:sasha 1}--\eqref{eqn:sasha 2} below: \\

\begin{proposition}
\label{prop:sasha}

For all $k,k'\in \BN$ and $i\in \{1,...,n\}$, we have:
\begin{equation}
\label{eqn:sasha 1}
[G_{\pm k}, G_{\pm k'}] = 0 
\end{equation}
\begin{equation}
\label{eqn:sasha 2}
[G_{\pm k}, 1_i^\pm] = 0
\end{equation}
where $1_i^\pm \in \CA_{\pm \bs^i}$ is the rational function 1 in a single variable $z_{i1}$. Moreover, if:
\begin{equation}
\label{eqn:new heis}
\sum_{k=0}^\infty G_{\pm k} x^k  =\exp \left( \sum_{k=1}^\infty \frac {P_{\pm k}x^k}k \right)
\end{equation}
then the elements $P_{\pm k} \in \CB_0$ are primitive for the coproduct $\Delta_0$, and:
\begin{equation}
\label{eqn:new pair heis}
\langle P_{-k}, P_l \rangle = \delta_k^l k \cdot \frac {q^{nk}-q^{-nk}}{(\oq^k - \oq^{-k})(q^{-nk}\oq^{-k} - q^{nk} \oq^{k})}
\end{equation}

\end{proposition} 

\subsection{}\label{sub:primitive}

We will prove the Proposition above in the Appendix, as it relies on a certain criterion for when a shuffle element is 0. To state this criterion, let us study the vector subspace of primitive elements:
$$
\CB_{0|\pm \bk}^{\prim} \ \subset \ \CB_{0|\pm \bk}
$$
for the coproduct $\Delta_0$. By definition, a shuffle element $R$ is primitive if and only if the limits \eqref{eqn:cristian1} or \eqref{eqn:cristian2} vanish for $\mu = 0$ and all $\bl \notin \{0,\bk\}$. As an immediate corollary of Proposition \ref{prop:sam}, we have:
\begin{equation}
\label{eqn:vacuous}
\CB_{0|\pm \bk}^{\prim} \ni R \text{ equals }0 \qquad \Leftrightarrow \qquad \alpha_{[i;j)}(R) = 0
\end{equation}
for all $i<j$ such that $\bk = [i;j)$.\footnote{The analogous statement holds if we replace $\alpha$ by $\beta$} Depending on $\bk$, we are in one of three cases: \\

\begin{itemize} 
	
\item If $\bk = [i;j)$ with $j \not \equiv i$ mod $n$, then: 
$$
\dim \ \CB_{0|\pm \bk}^{\prim} \leq 1 
$$
Primitive elements are determined by their image under $\alpha_{[i;j)}$ \\

\item If $\bk = k\bde$ for some $k\in \BN$ then: 
$$
\dim \ \CB_{0|\pm \bk}^{\prim} \leq n
$$
Primitive elements are determined by their images under $\{\alpha_{[i;i+nk)}\}_{1 \leq i \leq n}$ \\

\item For all other $\bk$, we have $\CB_{0|\pm \bk}^{\prim} = 0$. The reason for this is that $\bk \neq [i;j)$ for all $i<j$, and therefore the condition $\alpha_{[i;j)}(R) = 0$ in \eqref{eqn:vacuous} is vacuous. \\

\end{itemize}

\noindent In the situation of the second bullet, let us recall the order $n$ automorphism $\BZ/n\BZ \curvearrowright \CA^\pm$. Since the homogeneous degree of $R$ is 0, we observe that $\alpha_{[i;i+nk)}$ does not depend on $i$ if $R$ is $\BZ/n\BZ-$invariant. We conclude the following: \\

\begin{corollary}
\label{cor:sam} 

Up to a constant multiple, $\CB_{0|\pm \bk}$ has: \\

\begin{itemize}
	
\item at most one primitive element if $\bk = [i;j)$ for $j\not \equiv i$ mod $n$ \\	

\item at most one primitive, $\BZ/n\BZ-$invariant element if $\bk = (k,...,k)$ \\

\item no primitive elements otherwise. \\
	
\end{itemize}

\end{corollary}

\begin{conjecture}
\label{conj:sam}

There exists a single primitive element in $\CB_{0|\pm (k,...,k)}$, $\forall k \in \BN$. Moreover, if $j\not \equiv i$ mod $n$, then $\CB_{0|\pm [i;j)}$ contains a primitive element iff $j=i+1$. \\

\end{conjecture}

\begin{proposition} 
\label{prop:vasile}
	
There is a bialgebra isomorphism:
\begin{equation}
\label{eqn:yota}
\su \otimes \uui \ \stackrel{\sim}\longrightarrow \ \CB_0
\end{equation}
given by $c \mapsto c$,
\begin{equation}
\label{eqn:hot}
x_i^+ \mapsto \frac {1^+_i}{q^{-1} - q}, \qquad x_i^- \mapsto \frac {1^-_i}{1 - q^{-2}}, \qquad \ph_i \mapsto \frac {\psi_{i+1}}{\psi_i} 
\end{equation}
and the group-like elements $g_{\pm k} \in \uui$ go to the elements $G_{\pm k}$ of \eqref{eqn:g1}--\eqref{eqn:g2} (after applying a plethysm \eqref{eqn:plethysm}, in order to make \eqref{eqn:pie} correspond to \eqref{eqn:new pair heis}). \\

\end{proposition}

\begin{proof} First of all, note that the assignments \eqref{eqn:hot} give rise to a homomorphism:
\begin{equation}
\label{eqn:subalg 1}
\su \hookrightarrow \CB_0
\end{equation}
of bialgebras. Indeed, that one obtains an algebra homomorphism is a consequence of the fact that the shuffle elements appearing in \eqref{eqn:hot} satisfy relations \eqref{eqn:sug1}, \eqref{eqn:sug2}, \eqref{eqn:sug3}, \eqref{eqn:basiccomm}, which are all easy to check. The fact that the coproduct is respected can be seen by comparing \eqref{eqn:cop 1} with \eqref{eqn:cristian1} for $\mu = 0$, and the fact that the pairing is respected follows by comparing \eqref{eqn:pair 1} with \eqref{eqn:pairshuf} for $R^\pm = 1^\pm_i$. \\

\noindent Furthermore, it is clear that the elements $G_{\pm k} \in \CB_0$ give rise to a homomorphism:
\begin{equation}
\label{eqn:subalg 2}
\uui \hookrightarrow \CB_0
\end{equation}
because \eqref{eqn:nicecop1}, \eqref{eqn:nicecop2}, \eqref{eqn:sasha 1}, \eqref{eqn:new pair heis} imply that the elements $P_{\pm k} \in \CB_0$ defined by \eqref{eqn:new heis} are primitive and satisfy the commutation relations:
\begin{equation}
\label{eqn:new comm heis}
[P_k, P_l] = \delta_{k+l}^0 k \cdot \frac {(q^{nk}-q^{-nk})(c^k - c^{-k})}{(\oq^k - \oq^{-k})(q^{nk} \oq^{k} - q^{-nk}\oq^{-k})}
\end{equation}
The fact that the subalgebras \eqref{eqn:subalg 1} and \eqref{eqn:subalg 2} commute follows from \eqref{eqn:sasha 2} and:
$$
\langle G_{-k}, 1_i^+ \rangle = \langle 1_i^-, G_k \rangle = 0
$$
which holds for degree reasons, and the assumption $n>1$. Therefore, we obtain a bialgebra homomorphism \eqref{eqn:yota}, which preserves the triangular decompositions of its domain and target. This homomorphism is injective, because it is injective on all the factors in its triangular decomposition (the latter fact follows from Lemma \ref{lem:inj} and the fact that the tensor product of non-degenerate pairings is non-degenerate). To show that \eqref{eqn:yota} is an isomorphism, it is therefore enough to prove that the dimension of the positive/negative halves of its codomain are at most equal to the dimension of the positive/negative halves of its domain, in every degree. This follows from \eqref{eqn:mugabe} and \eqref{eqn:bound} for $\mu = 0$. \\	
\end{proof}

\subsection{}\label{sub:iso}

In this Section, we will prove the following result, akin to Proposition \ref{prop:vasile}. \\

\begin{proposition}
\label{prop:dumitru}

There is a bialgebra epimorphism:
\begin{equation}
\label{eqn:donny}
\Xi : \uu \ \twoheadrightarrow \ \CB_0 
\end{equation}
preserving the bialgebra pairing and the actions of $\BZ/n\BZ$ on the two sides, given by:
\begin{equation}
\label{eqn:buscemi}
e_{[i;j)} \mapsto E^0_{[i;j)} \qquad \qquad f_{[i;j)} \mapsto F^0_{[i;j)} \qquad \qquad \psi_i \mapsto \psi_i
\end{equation}

\end{proposition}

\begin{proof} The main observation is that the \eqref{eqn:buscemi} respects the coproduct and pairing of the generators, simply by comparing \eqref{eqn:quant2}, \eqref{eqn:quant3}, \eqref{eqn:quant5} with \eqref{eqn:cop1}, \eqref{eqn:cop2}, \eqref{eqn:tyler}. By the construction of the Drinfeld double, it is then sufficient to show that the assignment \eqref{eqn:buscemi} gives rise to an algebra epimorphism on half of the algebra:
\begin{equation}
\label{eqn:steve}
\uup \twoheadrightarrow \CB_0^+
\end{equation}
(the case when the sign is $-$ is analogous). To prove that \eqref{eqn:steve} preserves all multiplicative relations between the generators $e_{[i;j)}$, we must show that:
\begin{multline}
\text{if } r  := \sum_{i_1<j_1,...,i_t<j_t} \gamma \cdot e_{[i_1;j_1)}...e_{[i_t;j_t)} \quad \text{equals } 0, \label{eqn:answer} \\ 
\text{then } R := \sum_{i_1<j_1,...,i_t<j_t} \gamma \cdot E^0_{[i_1;j_1)} * ... * E^0_{[i_t;j_t)} \quad \text{is also } 0
\end{multline} 
where $\gamma$ is a place-holder for various coefficients. If $r = 0$, then $r$ pairs trivially with all products of $f_{[i;j)}$'s. Because of properties \eqref{eqn:bialg 1} and \eqref{eqn:bialg 2} of the bialgebra pairing, this implies that $R$ pairs trivially with all products of various $F_{[i;j)}$. According to Proposition \ref{prop:pairing}, this implies that $R = 0$. \\

\noindent The argument provided above implies not only that \eqref{eqn:steve} is a well-defined algebra homomorphism, but also the fact that:
$$
\text{if } R \in \CB_0^+ \text{ pairs trivially with Im } \Xi \quad \Rightarrow \quad R = 0
$$
Because the bialgebra pairing between $\CB_0^+$ and $\CB_0^-$ is non-degenerate (Prop. \ref{prop:vasile}) and the graded pieces of $\CB_0^+$ are finite-dimensional, this implies that $\Xi$ is surjective.

\end{proof}

\noindent In fact, the map \eqref{eqn:donny} is an isomorphism, as shown below. \\

\begin{proof} \emph{of Proposition \ref{prop:iso}:} Composing the epimorphism \eqref{eqn:donny} with the inverse of the isomorphism \eqref{eqn:yota}, we obtain an epimorphism:
$$
\uu \twoheadrightarrow \su \otimes \uui \ \cong \ \CB_0
$$
which preserves degrees. Then \eqref{eqn:dimboundaff} and \eqref{eqn:mugabe} imply that the above must be an isomorphism, as must be \eqref{eqn:donny} and hence \eqref{eqn:decomposition}. The inverse images:
$$
\Xi^{-1}(G_{\pm k}) \in \uu
$$
are $\BZ/n\BZ-$invariant group-like elements, which generate $\uui \subset \uu$. The uniqueness statement follows from Corollary \ref{cor:sam}. 

\end{proof}

\subsection{}\label{sub:dimension}

The fact that $\Xi$ of \eqref{eqn:donny} is an isomorphism, which was proved in the previous Subsection, establishes Lemma \ref{lem:sub} for $\mu = 0$. Let us now discuss the general case: 
\begin{equation}
\label{eqn:mu}
\mu \ = \ \frac ba
\end{equation}
with $\gcd(a,b) = 1$. A partition $C$ as in \eqref{eqn:partition} is $\mu-$integral if and only if $a$ divides the length of its constituent arcs $[i_s;j_s)$, for all $s\in \{1,...,t\}$. All these constituent arcs must therefore be obtained by stringing together the \underline{basic arcs}: 
\begin{equation}
\label{eqn:widehat}
\widehat{i} \ := \ [i;i+a) 
\end{equation}
as $i$ varies over all residues modulo $n$. We draw an oriented edge for all $i$:
$$
\widehat{i} \ \longrightarrow \ \widehat{i+a}
$$
modulo $n$, and this has the effect of dividing up the set $\{\widehat{1},...,\widehat{n}\}$ into $g := \gcd(n,a)$ disjoint cycles $C_1,...,C_g$. Any $\mu-$integral arc of length $la$ in the cyclic quiver $\{1,...,n\}$ simply corresponds to an arc of length $l$ in one of the cycles $C_1,...,C_g$:
$$
[i;i+la) = [i;i+a) \sqcup ...\sqcup [i+(l-1)a;i+la) \ \leftrightsquigarrow \ \left\{ \widehat{i}, \widehat{i+a}, ..., \widehat{i+(l-1)a} \right\}
$$ 
Then \eqref{eqn:bound} implies that:
$$
\dim  \CB_{\mu|\pm \bk} \leq \# \ \Big\{ \text{partitions of }\bk \text{ into arcs from the cycles }C_1,...,C_g \Big \}
$$
where each node $\widehat{i} \in C_1 \sqcup ... \sqcup C_g$ contributes the degree vector $[i;i+a)$ to the partition of $\bk$. Then Proposition \ref{prop:iso} and \eqref{eqn:mugabe} imply that:
\begin{equation}
\label{eqn:dimensionbound}
\dim  \CB_{\mu| \pm \bk} \leq \dim \ U^\pm_q(\dot{\fgl}_{\frac ng})^{\otimes g} \ \text{in degree } \pm \bk 
\end{equation}
The grading in the right-hand side assigns degree $[p+ia;p+(i+1)a) \in \nn$ to the $i-$th simple root in the $p-$th tensor factor, which corresponds to $\widehat{p+ia}$ of \eqref{eqn:widehat}. \\

\begin{proof} \emph{of Lemma \ref{lem:sub}:} We claim that the assignment:
\begin{equation}
\label{eqn:roz}
e_{[i;i+l)}^{(p)} \mapsto E_{[p+ia;p+(i+l)a)}^{\mu}, \qquad \qquad \psi_i^{(p)} \mapsto \psi_{p+ia} 
\end{equation}
yields an bialgebra homomorphism:
\begin{equation}
\label{eqn:pink}
\Xi: U^\geq_{q}(\dot{\fgl}_{\frac n{g}})^{\otimes g} \longrightarrow \CB^\geq_{\frac ba}
\end{equation}
Indeed, comparing \eqref{eqn:quant2}, \eqref{eqn:quant3}, \eqref{eqn:quant5} with \eqref{eqn:cop1}, \eqref{eqn:cop2}, \eqref{eqn:tyler} shows that the assignment \eqref{eqn:roz} respects the coproduct and pairings of the generators. To show that $\Xi$ is a well-defined homomorphism, we must prove the analogue of statement \eqref{eqn:answer}. The argument is identical to the case $\mu = 0$, from the proof of Proposition \ref{prop:dumitru}, and so we leave it to the interested reader (the key fact is the last sentence of Proposition \ref{prop:pairing}). We have thus constructed a bialgebra homomorphism $\Xi$. \\

\noindent Lemma \ref{lem:inj} implies that $\Xi$ is injective, since the bialgebra pairing between the halves of $\uu$ is non-degenerate, as a consequence of the isomorphism \eqref{eqn:decomposition}. Then inequality \eqref{eqn:dimensionbound} implies that $\Xi$ is surjective, and also that we have equality:
\begin{equation}
\label{eqn:dimensionbound 2}
\dim  \CB_{\mu| \pm \bk} = \dim \ U^\pm_q(\dot{\fgl}_{\frac ng})^{\otimes g} \ \text{in degree } \pm \bk 
\end{equation}

\end{proof}

\noindent An equivalent reformulation of \eqref{eqn:dimensionbound 2} is that \eqref{eqn:bound} becomes an equality:
\begin{equation}
\label{eqn:dimensionbound 3}
\dim \CB_{\mu|\pm \bk} = \# \ \Big \{\mu-\text{integral partitions }C \vdash \bk \Big \}
\end{equation}
for all $\mu \in \BQ$ and $\bk \in \nn$. \\

\subsection{}\label{sub:final}

We will now complete the proof of Theorem \ref{thm:iso}, by proving Proposition \ref{prop:surj}. \\ 

\begin{proof} \emph{of Proposition \ref{prop:surj}:} We need to show that:
$$
\CA^\pm  = \text{Im }\Upsilon^\pm
$$
As was shown in Proposition \ref{prop:belong}, the particular shuffle elements $S^+_m, T^-_m$ belong to $\text{Im }\Upsilon^\pm$, for all Laurent polynomials $m$. Hence so do the shuffle elements $E_{[i;j)}^\mu, F_{[i;j)}^\mu$ for all $\mu \in \BQ$ and all $i<j$, and Lemma \ref{lem:sub} therefore implies that:
$$
\CB^\pm_\mu \subset \text{Im }\Upsilon^\pm
$$
All that remains to prove is that the subalgebras $\{ \CB^\pm_\mu \}_{\mu \in \BQ}$ generate $\CA^\pm$. Since $\CA^\pm$ is infinite-dimensional in most degrees, we will actually prove a finer result:
\begin{equation}
\label{eqn:foxy}
\CA^\pm_{\leq \mu} = \prod^{\rightarrow}_{\BQ \ni \mu' \leq \mu} \CB^\pm_{\mu'}
\end{equation}
for all $\mu \in \BQ$, where we recall that $\CA^\pm_{\leq \mu}$ is the subalgebra of shuffle elements of slope $\leq \mu$. To prove \eqref{eqn:foxy}, pick a pair of dual bases $\{v_{\pm i}^\mu\}$ of the subalgebras $\CB^{\pm}_\mu$, where $i$ ranges over some indexing set. Then it is enough to show that the elements:
\begin{equation}
\label{eqn:joplin}
\Big \{v^{\mu_1}_{\pm i_1} * ... * v^{\mu_t}_{\pm i_t} \Big\}_{\mu_1<\mu_2<...<\mu_t\leq \mu}^{i_1,...,i_t} 
\end{equation}
are dual bases of $\CA_{\leq \mu}^\pm$. The number of elements \eqref{eqn:joplin} in degree $\pm \bk \in \pm \nn$ is:
\begin{equation}
\label{eqn:number}
\sum_{\bk_1 + ... + \bk_t = \bk}^{\mu_1 < ... < \mu_t \leq \mu} \dim \CB_{\mu_1|\pm \bk_1} \cdot ... \cdot \dim \CB_{\mu_t | \pm \bk_t} \stackrel{\eqref{eqn:dimensionbound 2}}= 
\end{equation}
$$
= \# \Big \{\text{partitions }C \vdash \bk \text{ into arcs, each of which is } \tilde{\mu}\text{--integral for some } \tilde{\mu} \leq \mu \Big \}
$$
The right-hand side is equal to the upper bound on $\dim \CA_{\leq \mu|\pm \bk}$ that we discussed in Lemma \ref{lem:min} (strictly speaking, one needs to fix the homogeneous degree $d \in \BZ$ in order for the number \eqref{eqn:number} to be finite and for the argument above to be completely precise). Therefore, all that remains is to show that:
\begin{equation}
\label{eqn:lady}
\left \langle v^{\mu_1}_{-i_1} * ... * v^{\mu_t}_{-i_t}, v^{\mu'_1}_{+j_1} * ... * v^{\mu'_s}_{+j_s} \right \rangle = \text{Kronecker delta}
\end{equation}
which would prove that the products \eqref{eqn:joplin} are linearly independent, and therefore yield a basis. We may assume without loss of generality that $\mu_1' \geq \mu_1$, otherwise we switch the roles of $+$ and $-$. Then applying property \eqref{eqn:bialg 1}, we see that:
\begin{equation}
\label{eqn:hendrix}
\text{LHS of \eqref{eqn:lady}} = \left \langle v^{\mu_1}_{-i_1} \otimes v^{\mu_2}_{-i_2}...v^{\mu_t}_{-i_t}, \Delta\left(v^{\mu'_1}_{+j_1}... v^{\mu'_s}_{+j_s}\right)\right \rangle
\end{equation}
For a homogeneous shuffle element $R^\pm \in \CA_{\pm \bk}$, let us denote its \underline{naive slope} as:
$$
\frac {\homdeg R^\pm}{\pm |\bk|}
$$
Consider the coproduct of any element $v_{i}^\mu \in \CB_\mu^+$: formula \eqref{eqn:radu1} implies that all second tensor factors in the coproduct have slope $ \leq \mu$, hence naive slope $\leq \mu$, hence all first tensor factors in the coproduct have naive slope $\geq \mu$. Therefore, the coproduct in \eqref{eqn:hendrix} only has first tensor factors of naive slope $\geq \mu_1'$. So if $\mu_1' > \mu_1$, the pairing \eqref{eqn:hendrix} is trivial. By assumption, the only remaining case is $\mu_1' = \mu_1$, so we may assume that the degree $-\bk \in -\nn$ of $v^{\mu_1}_{-i_1}$ is no less than the degree $\bk'$ of $v^{\mu_1}_{+j_1}$ (i.e. $\bk \not < \bk'$). Then the only non-trivial term in the pairing \eqref{eqn:hendrix} is:
\begin{equation}
\label{eqn:jared}
\left \langle v^{\mu_1}_{-i_1} \otimes v^{\mu_2}_{-i_2}...v^{\mu_t}_{-i_t}, v^{\mu_1}_{+j_1} \otimes v^{\mu'_2}_{+j_2} ... v^{\mu'_s}_{+j_s} \right \rangle
\end{equation}
(the reader may note that there should be a product of Cartan elements in the second argument of the pairing, but this product does not change the value of \eqref{eqn:jared}). Since $\{v^{\mu_1}_{\pm i}\}$ are dual bases, the pairing \eqref{eqn:jared} is non-trivial only if $i_1=j_1$. We may repeat the argument to prove that non-triviality of the pairing forces $\mu_2'=\mu_2$ and $i_2=j_2$ etc. This proves that \eqref{eqn:joplin} are dual bases of $\CA^\pm_{\leq \mu}$ for any $\mu$. Thus $\CA^\pm_{\leq \mu}$ is generated by $\{\CB^\pm_{\mu'}\}_{\mu' \leq \mu}$, and hence is contained in the image of $\Upsilon^\pm$. \\	
\end{proof}

\subsection{}\label{sub:R}

The proof of Proposition \ref{prop:surj} in the previous Subsection not only shows that:
$$
\CA^\pm \ = \ \prod_{\mu \in \BQ}^\rightarrow \CB^\pm_\mu 
$$
(note that the decomposition above only holds as vector spaces, not as algebras, i.e. the various factors $\CB_\mu^\pm$ do not commute; their commutators are studied in \cite{PBW}) but it shows that the above equality preserves the pairing. In other words:
\begin{equation}
\label{eqn:pairings}
\Big \langle \prod_{\mu \in \BQ}^\rightarrow v^\mu_-, \prod_{\mu \in \BQ}^\rightarrow v^\mu_+ \Big \rangle_{\CA} = \prod_{\mu \in \BQ}^\rightarrow \langle v^\mu_-, v^\mu_+ \rangle_{\CB_\mu}
\end{equation}
which implies that dual bases of $\CA^\pm$ are given by products of dual bases of $\{\CB_\mu\}_{\mu \in \BQ}$. We may also include $\mu = \infty$ in the aforementioned analysis, and infer that the bialgebra pairing is non-degenerate on the subalgebras:
\begin{equation}
\label{eqn:pairshuffinal}
\left( \prod_{\mu \in \BQ \sqcup \infty}^{\rightarrow} \CB^+_\mu \right) \bigotimes \left( \prod_{\mu \in \BQ \sqcup \infty}^{\rightarrow} \CB^-_\mu \right)  \rightarrow \BF
\end{equation}
The two algebras in the brackets match $\CA^\geq$ and $\CA^\leq$ from Proposition \ref{prop:pairshuf}, except that they lack the invertible Cartan elements $\psi_i$ and their various products. \\

\begin{definition}

The canonical tensor of the pairing \eqref{eqn:pairshuffinal}, denoted by: 
$$
R_{\CA} \in \CA\woo \CA
$$
will be called the \underline{universal$^*$ $R$--matrix}. \\

\end{definition}

\begin{remark}

We have $R_{\CA} = \CR_1 \CR_0$, in the notation in Subsection 2.2 of \cite{FJM}. As shown in \loccitt, the actual universal $R$--matrix of $\CA$ would be equal to $R_{\CA}$ times $q$ raised to a certain linear combination of Cartan elements, but one would need to make certain modifications in order for this to be well-defined: set $q = e^{\hbar}$ and work over $\BQ((\hbar))$, and introduce certain additional Cartan elements into the algebra $\CA$. \\

\end{remark}

\noindent One can consider the analogous notion of universal$^*$ $R$--matrix for the individual subalgebras $\CB_\mu$, and then \eqref{eqn:pairings} implies that:
\begin{equation}
\label{eqn:r a}
R_{\CA} = \prod^{\rightarrow}_{\mu \in \BQ \cup \{\infty\}} R_{\CB_{\mu}} \ \in \ \CA \widehat{\otimes} \CA
\end{equation}
Using Theorem \ref{thm:iso} and Lemma \ref{lem:sub}, we obtain the factorization result \eqref{eqn:R}:
\begin{equation}
\label{eqn:older}
R_{\UU} = \prod_{\frac ba \in \BQ \cup \{\infty\}}^{\rightarrow} R_{\uux^{\otimes g}} \ \in \ \UU \ \widehat{\otimes} \ \UU
\end{equation}
where we write $g = \gcd(n,a)$ for each factor in the right-hand side. The factorization \eqref{eqn:older} may be interpreted as the toroidal version of the well-known product formulas for $R$--matrices of quantum groups (see \cite{Drin, KR, LS, R}). In particular, the quantum affine case considered in \cite{KT, LSS} reads:
\begin{equation}
\label{eqn:kt}
R_{\su} = \prod_{\text{positive root }\alpha} R_{U_q(\fsl_2)-\text{subalgebra corresponding to }\alpha}
\end{equation}
and in fact \eqref{eqn:kt} allows us to further factor the products in the RHS of \eqref{eqn:older}. \\

\section{Appendix} 
\label{sec:app}

\begin{proof} \emph{of Proposition \ref{prop:rtt}:} Let us first check the basic properties of the coproduct $\Delta$ of \eqref{eqn:copy1}--\eqref{eqn:copy2}. The fact that it is coassociative is immediate from the definition. To show that $\Delta$ extends multiplicatively from the generators to the entire $\uug$ and $\uul$, we need to prove that it respects relation \eqref{eqn:rtt}. Let us take care of the $\geq$ case and leave the $\leq$ case to the interested reader:
$$
\Delta\left(R \left( \frac xy \right) T^+_1(x)T^+_2(y) \right) = R\left(\frac xy \right) \Delta(T^+_1(x))\Delta(T^+_2(y)) = 
$$
$$
= R\left(\frac xy \right) \left(T_1^+(x) \otimes T_1^+(x/c_1) \right) \left(T^+_2(y) \otimes T^+_2(y/c_1) \right) = 
$$
$$
= R\left(\frac xy \right) T_1^+(x) T_2^+(y) \otimes T^+_1(x/c_1) T^+_2(y/c_1) = 
$$
$$
= T^+_2(y) T_1^+(x) R\left(\frac xy \right) \otimes T^+_1(x/c_1) T_2^+(y/c_1) =
$$
$$
= T^+_2(y) T^+_1(x) \otimes R\left(\frac xy \right) T^+_1(x/c_1) T^+_2(y/c_1) =
$$
$$
= T^+_2(y) T^+_1(x) \otimes T^+_2(y/c_1) T^+_1(x/c_1) R\left(\frac xy \right) =
$$
$$
= \left(T_2^+(y) \otimes T_2^+(y/c_1)\right) \left(T_1^+(x) \otimes T_1^+(x/c_1) \right)R\left(\frac xy \right) = 
$$
$$
= \Delta(T^+_2(y))\Delta(T^+_1(x)) R\left(\frac xy \right) = \Delta\left(T^+_2(y)  T^+_1(x) R \left( \frac xy \right)\right)
$$
In order to prove that \eqref{eqn:bialg2} extends to a bialgebra pairing on the whole algebras, we need to check that it preserves the defining relations \eqref{eqn:rtt}. In other words, we need to check that:
$$
\left \langle R_{12}\left(\frac xy \right) T_1^-(x) T_2^-(y) ,T_3^+(z) \right \rangle = \left \langle  T_2^-(y) T_1^-(x) R_{12}\left(\frac xy \right) ,T_3^+(z) \right \rangle 
$$
where now $R_{12} = R \boxtimes 1$, and $T_1=T \boxtimes 1 \boxtimes 1$, $T_2=1 \boxtimes T \boxtimes 1$, $T_3=1 \boxtimes 1 \boxtimes T$. The situation where the roles of $+$ and $-$ are switched is analogous, and left as an exercise to the interested reader. Using the bialgebra property \eqref{eqn:bialg 1}, we see that: 
$$
\left \langle R_{12}\left(\frac xy \right) T_1^-(x) T_2^-(y) ,T_3^+(z) \right \rangle = R_{12}\left(\frac xy \right) \left \langle T_1^-(x) \otimes T_2^-(y) , \Delta \left( T_3^+(z) \right) \right \rangle = 
$$
$$
= R_{12}\left(\frac xy \right) \left \langle T_1^-(x), T_3^+(z) \right \rangle \left \langle T_2^-(y), T_3^+(z) \right \rangle = R_{12}\left(\frac xy \right)  R_{13} \left(\frac xz \right)  R_{23}\left(\frac yz \right) = 
$$
$$
= R_{23}\left(\frac yz \right) R_{13}\left(\frac xz \right) R_{12}\left(\frac xy \right) = \left \langle T_2^-(y), T_3^+(z) \right \rangle \left \langle T_1^-(x), T_3^+(z) \right \rangle  R_{12}\left(\frac xy \right) = 
$$
$$
= \left \langle T_2^-(y) \otimes T_1^-(x) , \Delta \left( T_3^+(z) \right) \right \rangle R_{12}\left(\frac xy \right)  = \left \langle  T_2^-(y) T_1^-(x) R_{12}\left(\frac xy \right) ,T_3^+(z) \right \rangle
$$
where the equality between the second and third rows is simply the quantum Yang-Baxter equation \eqref{eqn:yb}. The central element $c$ did not enter the above computation, because:
\begin{equation}
\label{eqn:pair triv}
\langle ca,b \rangle = \langle a, bc \rangle = \langle a,b \rangle \qquad \forall a \in \uul, b \in \uug
\end{equation}
Finally, the commutation relations between the positive and negative halves of the Drinfeld double are prescribed by \eqref{eqn:drinfeld}. For $a = T_1^-(x)$ and $b = T_2^+(y)$, we obtain:
$$
\left \langle T_1^-(x/c), T_2^+(y) \right \rangle T_1^-(x) T_2^+(y) = T_2^+(y) T_1^-(x) \left \langle T_1^-(x), T_2^+(y/c) \right \rangle
$$
(we ignore the central elements $c$ which enter the pairings, due to \eqref{eqn:pair triv}). As a consequence of formula \eqref{eqn:bialg2}, the equality above implies \eqref{eqn:unwind}.

\end{proof}

\begin{proof} \emph{of Proposition \ref{prop:dim}:} We will prove the statement for $\uup$, since the case when the sign $+$ is replaced with $-$ is analogous. For any root generator $e_{[i;i+l)}$, we define its \emph{length} as the natural number $l$. We will refer to products of root generators as \underline{monomials}, and define the \underline{signature} of such a monomial:
\begin{equation}
\label{eqn:mony}
m = e_{[i_1;i_1+l_1)} ... e_{[i_d;i_d+l_d)} \qquad \text{as the tuple} \qquad (l_1,...,l_d)
\end{equation}
We restrict attention to all monomials of fixed total length $l_1+...+l_d$. Such signatures can be ordered lexicographically, i.e. write $(l_1,...,l_d) <  (l_1,...,l_{d'})$:
\begin{align*}
&\text{if } d<d' \text{ or} \\
&\text{if }d = d' \text{ and } a_1 = a_1', \ ..., \ a_{i-1}=a_{i-1}', \ a_i<a_i' \text{ for some }i 
\end{align*}
where $\{ a_1\geq ... \geq a_d\} = \{l_1,...,l_d\}$ and $\{a'_1\geq ... \geq a'_{d'}\} = \{l'_1,...,l'_{d'}\}$. We define the partial \underline{lexicographic} ordering on monomials to be given by the lexicographic ordering of their signatures. As for a single monomial \eqref{eqn:mony}, we call it \underline{ordered} if:
$$
l_s \leq l_{s+1} \quad \forall s \qquad \text{with equality only if} \qquad i_s \leq i_{s+1}
$$
Clearly, the number of ordered monomials of a given degree $\bk \in \nn$ equals the number of partitions of $\bk$ into arcs $[i;i+l)$. Therefore, the Proposition reduces to the claim that any product of root generators can be expressed as a linear combination of ordered monomials. To prove this, it suffices to prove that:
\begin{equation}
\label{eqn:refined}
m \in \text{span} \Big \{\text{ordered monomials}, \text{monomials lexicographically} < m \Big \}  
\end{equation}
for any monomial $m$ as in \eqref{eqn:mony}, and then proceed by induction in lexicographic ordering. To prove \eqref{eqn:refined}, let us rewrite \eqref{eqn:expre 3} as:
$$
e_{[i;i+l)}e_{[i';i'+l')} \in \gamma e_{[i';i'+l')}e_{[i;i+l)} + \sum \gamma e_{\text{length}<l'} e_{\text{length}>l} + \sum \gamma  e_{\text{length}>l}e_{\text{length}<l'}
$$
The symbols $\gamma$ in front of the monomials indicate constants we will not need to bother with. Assuming $l\geq l'$, we can iterate the above formula to reorder the two $e$'s in the last sum:
\begin{equation}
\label{eqn:ca}
e_{[i;i+l)}e_{[i';i'+l')} \in \gamma e_{[i';i'+l')}e_{[i;i+l)} + \sum \gamma e_{\text{length}<l'} e_{\text{length}>l} + \sum \gamma e_{\text{arc}}
\end{equation}
where the last summand consists of a single root generator $e_{\text{arc}}$, instead of a product of two such root generators. If the monomial $m$ of \eqref{eqn:mony} is not ordered, then choose the smallest $s$ such that $l_s>l_{s+1}$ or $l_s = l_{s+1}$ and $i_s>i_{s+1}$. We use formula \eqref{eqn:ca} to reorder the root generators $e_{[i_s;i_s+l_s)}$ and $e_{[i_{s+1};i_{s+1}+l_{s+1})}$, and observe that in the right hand side we obtain monomials which are smaller lexicographically. Repeating this argument allows us to completely order the monomial $m$, thus proving \eqref{eqn:refined}. \\
\end{proof}

\begin{proof} \emph{of Proposition \ref{prop:pairtor}:} We need to check that the properties of being a bialgebra pairing \eqref{eqn:bialg 1}--\eqref{eqn:bialg 2} are compatible with the defining relations \eqref{eqn:reltor1}--\eqref{eqn:serre} of the quantum toroidal algebra. For \eqref{eqn:reltor1}, relations \eqref{eqn:bialg 1} and \eqref{eqn:deltator2} imply that:
\begin{equation}
\label{eqn:floyd}
\left \langle \psi_j^-(w) x_i^-(z), x_{i'}^+(y) \right \rangle  = \left \langle \psi_j^-(w) \otimes x_i^-(z), \Delta(x_{i'}^+(y)) \right \rangle = 
\end{equation}
$$
= \left \langle \psi_j^-(w) \otimes x_i^-(z), \ph^+_{i'}(y) \otimes x_{i'}^+(y) \right \rangle = \left \langle \psi_j^-(w), \ph^+_{i'}(y) \right \rangle \cdot \left \langle x_i^-(z), x_{i'}^+(y) \right \rangle
$$
where $\ph^+_{i'}(y) = \psi_{i'+1}^+(yq^2) / \psi_{i'}^+(y)$. Using \eqref{eqn:bialg 2}, \eqref{eqn:deltator1} and \eqref{eqn:pairtor2}, we see that:
$$
\left \langle \psi_j^-(w), \ph^+_{i'}(y) \right \rangle = \left \langle \psi_j^-(w), \frac {\psi_{i'+1}^+(yq^2)}{\psi_{i'}^+(y)} \right \rangle = \frac {\left \langle \psi_j^-(w), \psi_{i'+1}^+(y q^2)\right \rangle}{\left \langle \psi_j^-(w), \psi_{i'}^+(y)\right \rangle} =
$$
\begin{equation}
\label{eqn:ohm1}
= \frac {\zeta \left( \frac {yq^2}w \right) \zeta \left( \frac {yq^4}w \right) \zeta \left( \frac {yq^6}w \right)...}{\zeta \left( \frac yw \right) \zeta \left( \frac {yq^2}w \right) \zeta \left( \frac {yq^4}w \right)...} = \zeta \left( \frac yw \right)^{-1} 
\end{equation}
Therefore, \eqref{eqn:floyd} becomes:
\begin{equation}
\label{eqn:lloyd}
\left \langle \psi_j^-(w) x_i^-(z), x_{i'}^+(y) \right \rangle  = \zeta \left( \frac yw \right)^{-1}  \left \langle x_i^-(z), x_{i'}^+(y) \right \rangle = \zeta \left( \frac zw \right)^{-1}  \frac {\delta_{i'}^i\delta\left(\frac yz\right)}{q^{-1} - q} 
\end{equation}
where in the last equality we used \eqref{eqn:pairtor1} and the following property of $\delta-$functions:
\begin{equation}
\label{eqn:prop delta}
f(y) \delta\left(\frac yz\right) \ = \ f(z) \delta\left(\frac yz\right)  \qquad \forall \ \text{power series } f 
\end{equation}
(indeed, according to Remark \ref{rem:note}, the functions $\zeta$ that appear in \eqref{eqn:lloyd} should be expanded as power series in $w$). Similarly, relations \eqref{eqn:bialg 1} and \eqref{eqn:deltator2} imply that:
$$
\left \langle x_i^-(z) \psi_j^-(w), x_{i'}^+(y)\right \rangle \cdot \zeta \left( \frac zw \right)^{-1}  = \left\langle x_i^-(z) \otimes \psi_j^-(w) , \Delta(x_{i'}^+(y)) \right \rangle = 
$$
\begin{equation}
\label{eqn:lloyds}
= \zeta \left( \frac zw \right)^{-1}  \cdot \left\langle x_i^-(z) \otimes \psi_j^-(w) , x_{i'}^+(y) \otimes 1 \right\rangle = \zeta \left( \frac zw \right)^{-1}  \cdot \frac {\delta_{i'}^i \delta\left(\frac yz\right)}{q^{-1} - q} 
\end{equation}
Comparing \eqref{eqn:lloyd} with \eqref{eqn:lloyds}, we see that the pairing respects relation \eqref{eqn:reltor1} when the sign is $-$ (to be completely thorough, one should prove that \eqref{eqn:lloyd} and \eqref{eqn:lloyds} are equal when one replaces $x_{i'}^+(y)$ by its product with an arbitrary number of factors $\psi_{j'}(x)$, but we leave this straightforward generalization to the interested reader). The computation when the signs $+$ and $-$ are switched is analogous, once one uses the following relation instead of \eqref{eqn:ohm1}:
$$
\left \langle \ph_i^-(y), \psi^+_{j}(w) \right \rangle = \left \langle \frac {\psi_{i+1}^-(yq^2)}{\psi_{i}^-(y)}, \psi^+_j(w) \right \rangle = \frac {\left \langle \psi_{i+1}^-(yq^2), \psi_j^+(w)\right \rangle}{\left \langle \psi_{i}^-(y), \psi_j^+(w)\right \rangle} =
$$
\begin{equation}
\label{eqn:ohm2}
= \frac {\zeta \left( \frac w{yq^2} \right) \zeta \left( \frac {wq^2}{yq^2} \right) ...}{\zeta \left( \frac wy \right) \zeta \left(\frac {wq^2}y \right)...} = \underline{\zeta \left( \frac w{yq^2} \right)} = \underline{\underline{\zeta \left( \frac yw \right)}} 
\end{equation}
where in the last equality we used \eqref{eqn:identity}. Note that for the purpose of defining the color-dependent rational functions $\zeta$ in \eqref{eqn:ohm2}, $y$ has color $i+1$ in the term with a single underline, and color $i$ in the term with a double underline. Then either \eqref{eqn:ohm1} or \eqref{eqn:ohm2} allow us to prove:
\begin{equation}
\label{eqn:ohm}
\left \langle \ph_i^-(y), \ph^+_{j}(w) \right \rangle = \frac {\left \langle \psi_{i+1}^-(yq^2), \ph_j^+(w)\right \rangle}{\left \langle \psi_{i}^-(y), \ph_j^+(w)\right \rangle} = \frac {\zeta (w/y)}{\zeta(w/yq^2)} = \frac {\zeta(w/y)}{\zeta(y/w)}
\end{equation}
where in the last equality, we again used \eqref{eqn:identity}. Before we set out to prove that \eqref{eqn:reltor2} and \eqref{eqn:serre} are preserved by the bialgebra pairing, let us compute the pairing of arbitrary products of currents by using \eqref{eqn:bialg 1}:
\begin{equation}
\label{eqn:vache}
\left \langle x_{j_1}^-(w_1)... x_{j_k}^-(w_k), x_{i_1}^+(z_1)... x_{i_k}^+(z_k) \right \rangle = \left \langle \bigotimes_{s=1}^k x_{j_s}^-(w_s) , \prod_{s=1}^k \Delta^{(k-1)} \left( x_{i_s}^+(z_s) \right) \right \rangle \qquad
\end{equation}
Iterating the coproduct of \eqref{eqn:deltator2} gives us:
$$
\Delta^{(k-1)}(x^+_i(z)) = \sum_{s=1}^k \ph_i^+(z) \otimes ... \otimes \ph_i^+(z) \otimes \underbrace{x_i^+(z)}_{s-\text{th factor}} \otimes 1 \otimes ... \otimes 1 
$$
and so we obtain:
\begin{equation}
\label{eqn:quirit}
\text{LHS of \eqref{eqn:vache}} = \sum_{\sigma \in S(k)} \prod_{s=1}^k \left \langle x_{j_s}^-(w_s), \prod^\rightarrow_{s < t \leq k} \ph^+_{i_{\sigma(t)}}(z_{\sigma(t)}) \cdot x^+_{i_{\sigma(s)}}(z_{\sigma(s)}) \right \rangle
\end{equation}
where the symbol $\prod^\rightarrow$ means that the $\ph$ and $x$ factors following it must be multiplied in increasing order of the index of $z$. We may use \eqref{eqn:reltor1} to obtain:
$$
\ph_j^+(w) x_i^+(z) = \frac {\psi_{j+1}^+(wq^2)}{\psi_{j}^+(w)} x_i^+(z) = x_i^+(z) \frac {\psi_{j+1}^+(wq^2)}{\psi_{j}^+(w)} \cdot \frac {\zeta \left(\frac z{wq^2} \right)}{\zeta \left( \frac zw \right)} = x_i^+(z) \ph_j^+(w) \cdot \frac {\zeta \left( \frac wz \right)}{\zeta \left( \frac zw \right)}
$$ 
in order to place all the $\ph$'s in \eqref{eqn:quirit} to the very right:
$$
\text{LHS of \eqref{eqn:vache}} = \sum_{\sigma \in S(k)} \prod_{s=1}^k \left \langle x_{j_s}^-(w_s), x^+_{i_{\sigma(s)}}(z_{\sigma(s)}) \prod^\rightarrow_{s < t \leq k} \ph^+_{i_{\sigma(t)}}(z_{\sigma(t)}) \right \rangle \prod_{s < t \leq k}^{\sigma(t)<\sigma(s)} \frac {\zeta \left( \frac {z_{\sigma(t)}}{z_{\sigma(s)}} \right)}{\zeta \left( \frac {z_{\sigma(s)}}{z_{\sigma(t)}} \right)}
$$
Then we can use \eqref{eqn:bialg 2} and \eqref{eqn:pairtor1} to evaluate the right-hand side above:
$$
\left \langle x_{j_1}^-(w_1)... x_{j_k}^-(w_k), x_{i_1}^+(z_1)... x_{i_k}^+(z_k) \right \rangle = \sum_{\sigma \in S(k)} \prod_{s=1}^k \frac {\delta^{j_s}_{i_{\sigma(s)}} \delta \left(\frac {w_s}{z_{\sigma(s)}} \right)}{q^{-1} - q} \prod_{s < t \leq k}^{\sigma(t)<\sigma(s)}  \frac {\zeta \left( \frac {z_{\sigma(t)}}{z_{\sigma(s)}} \right)}{\zeta \left( \frac {z_{\sigma(s)}}{z_{\sigma(t)}} \right)}  = 
$$
\begin{equation}
\label{eqn:rachel}
= \sum_{\sigma \in S(k)} \frac {\delta \left(\frac {w_1}{z_{\sigma(1)}} \right) ...  \delta \left(\frac {w_k}{z_{\sigma(k)}} \right)}{(q^{-1} - q)^k} 
\prod^{s < t}_{\sigma(t)<\sigma(s)} \frac {\zeta \left( \frac {w_t}{w_s} \right)}{\zeta \left( \frac {w_s}{w_t} \right)} 
\end{equation}
In the interest of space, in the second line above we make the convention that $\delta(w/z) = 0$ if $\col w \neq \col z$. The ratio of rational functions $\zeta$ must be interpreted as an equality of power series expanded in the range $|w_1| \ll ... \ll |w_k|$, in order for the second equality of \eqref{eqn:rachel} to be a correct application of \eqref{eqn:prop delta}. \\


\noindent Let us now prove that relation \eqref{eqn:reltor2} (interpreted as in the last sentence of Remark \ref{rem:note}) is preserved by the bialgebra pairing. Applying \eqref{eqn:rachel} for $k=2$ implies that:
\begin{equation}
\label{eqn:gog}
\left \langle x_i^-(z)x_j^-(w), x^+_{i'}(y_1) x^+_{j'}(y_2) \right \rangle = \frac {\delta\left(\frac z{y_1}\right)\delta\left(\frac w{y_2}\right) }{(q^{-1} - q)^2}+ \frac {\delta\left(\frac z{y_2}\right)\delta\left(\frac w{y_1}\right)}{(q^{-1} - q)^2} \frac {\zeta \left( \frac wz \right)}{\zeta \left( \frac zw \right)} \qquad
\end{equation}
At the same time, we have:
\begin{equation}
\label{eqn:gogo}
\left \langle x_i^-(w)x_j^-(z), x^+_{i'}(y_1) x^+_{j'}(y_2) \right \rangle \cdot \frac {\zeta \left( \frac wz \right)}{\zeta \left( \frac zw \right)} = \frac {\delta\left(\frac w{y_1}\right)\delta\left(\frac z{y_2}\right)}{(q^{-1} - q)^2} \frac {\zeta \left( \frac wz \right)}{\zeta \left( \frac zw \right)} + \frac {\delta\left(\frac w{y_2}\right)\delta\left(\frac z{y_1}\right) }{(q^{-1} - q)^2} \qquad
\end{equation}
Comparing the right-hand sides of the above formulas shows that they are equal, and so \eqref{eqn:reltor2} is preserved by the bialgebra property \eqref{eqn:bialg 1} (to be completely thorough, one should prove that the left-hand sides of \eqref{eqn:gog} and \eqref{eqn:gogo} are equal when one replaces $x^+_{i'}(y_1) x^+_{j'}(y_2)$ by its product with an arbitrary number of factors $\psi_{k'}(w)$, but we leave this as an exercise to the interested reader). The analogous computation when the signs $+$ and $-$ are switched, as well as the computation which checks \eqref{eqn:bialg 2}, are also left to the interested reader. \\

\noindent As for \eqref{eqn:serre}, let us use \eqref{eqn:rachel} for $k=3$:
$$
\left \langle x^-_{i + 1}(w) x^-_i(z) x^-_i(z'), x_{j_1}^+(y_1) x_{j_2}^+(y_2) x_{j_3}^+(y_3) \right \rangle = \frac {\delta \left(\frac {w}{y_3}\right) \delta \left(\frac {z}{y_2}\right) \delta \left(\frac {z'}{y_1}\right)}{(q^{-1} - q)^3} \frac {\zeta \left( \frac zw \right)\zeta \left(\frac {z'}w \right)\zeta \left( \frac {z'}z \right)}{\zeta \left( \frac wz \right) \zeta \left( \frac w{z'} \right) \zeta \left( \frac z{z'} \right)} +
$$
$$
+ \frac {\delta \left(\frac {w}{y_1}\right) \delta \left(\frac {z}{y_2}\right) \delta \left(\frac {z'}{y_3}\right)}{(q^{-1} - q)^3} + \frac {\delta \left(\frac {w}{y_2}\right) \delta \left(\frac {z}{y_1}\right) \delta \left(\frac {z'}{y_3}\right)}{(q^{-1} - q)^3} \frac {\zeta \left( \frac zw \right)}{\zeta \left( \frac wz \right)} + \frac {\delta \left(\frac {w}{y_1}\right) \delta \left(\frac {z}{y_3}\right) \delta \left(\frac {z'}{y_2}\right)}{(q^{-1} - q)^3} \frac {\zeta \left( \frac {z'}z \right)}{\zeta \left( \frac z{z'} \right)} +
$$
\begin{equation}
\label{eqn:1}
+ \frac {\delta \left(\frac {w}{y_2}\right) \delta \left(\frac {z}{y_3}\right) \delta \left(\frac {z'}{y_1}\right)}{(q^{-1} - q)^3} \frac {\zeta \left( \frac {z'}w \right)\zeta\left(\frac {z'}z \right)}{\zeta \left( \frac w{z'} \right)\zeta \left( \frac z{z'} \right)} + \frac {\delta \left(\frac {w}{y_3}\right) \delta \left(\frac {z}{y_1}\right) \delta \left(\frac {z'}{y_2}\right)}{(q^{-1} - q)^3} \frac {\zeta \left( \frac zw \right)\zeta \left( \frac {z'}w \right)}{\zeta \left( \frac wz \right) \zeta \left( \frac w{z'} \right)} \qquad
\end{equation}
$$
\left \langle x^-_i(z) x^-_{i + 1}(w) x^-_i(z'), x_{j_1}^+(y_1) x_{j_2}^+(y_2) x_{j_3}^+(y_3) \right \rangle = \frac {\delta \left(\frac {z}{y_3}\right) \delta \left(\frac {w}{y_2}\right) \delta \left(\frac {z'}{y_1}\right)}{(q^{-1} - q)^3} \frac {\zeta \left( \frac wz \right) \zeta \left( \frac {z'}z \right) \zeta \left( \frac {z'}w \right)}{\zeta \left( \frac zw \right) \zeta \left( \frac z{z'} \right) \zeta \left( \frac w{z'} \right)} +
$$
$$
+ \frac {\delta \left(\frac {z}{y_1}\right) \delta \left(\frac {w}{y_2}\right) \delta \left(\frac {z'}{y_3}\right)}{(q^{-1} - q)^3} + \frac {\delta \left(\frac {z}{y_2}\right) \delta \left(\frac {w}{y_1}\right) \delta \left(\frac {z'}{y_3}\right)}{(q^{-1} - q)^3} \frac {\zeta \left( \frac wz \right)}{\zeta \left( \frac zw \right)} + \frac {\delta \left(\frac {z}{y_1}\right) \delta \left(\frac {w}{y_3}\right) \delta \left(\frac {z'}{y_2}\right)}{(q^{-1} - q)^3} \frac {\zeta \left( \frac {z'}w \right)}{\zeta \left( \frac w{z'} \right)} +
$$
\begin{equation}
\label{eqn:2}
+ \frac {\delta \left(\frac {z}{y_2}\right) \delta \left(\frac {w}{y_3}\right) \delta \left(\frac {z'}{y_1}\right)}{(q^{-1} - q)^3} \frac {\zeta \left( \frac {z'}z \right) \zeta \left( \frac {z'}w \right)}{\zeta \left( \frac z{z'} \right) \zeta \left( \frac w{z'} \right)} + \frac {\delta \left(\frac {z}{y_3}\right) \delta \left(\frac {w}{y_1}\right) \delta \left(\frac {z'}{y_2}\right)}{(q^{-1} - q)^3} \frac {\zeta \left( \frac wz \right)\zeta \left( \frac {z'}z \right)}{\zeta \left( \frac zw \right)\zeta\left( \frac z{z'} \right)} \qquad
\end{equation}
$$
\left \langle x^-_i(z) x^-_i(z') x^-_{i + 1}(w), x_{j_1}^+(y_1) x_{j_2}^+(y_2) x_{j_3}^+(y_3) \right \rangle = \frac {\delta \left(\frac {z}{y_3}\right) \delta \left(\frac {z'}{y_2}\right)\delta \left(\frac {w}{y_1}\right)}{(q^{-1} - q)^3} \frac {\zeta \left( \frac {z'}z \right)\zeta \left( \frac wz \right)\zeta \left( \frac w{z'} \right)}{\zeta \left( \frac z{z'} \right)\zeta \left( \frac zw \right)\zeta \left( \frac {z'}w \right)} +
$$
$$
+ \frac {\delta \left(\frac {z}{y_1}\right) \delta \left(\frac {z'}{y_2}\right) \delta \left(\frac {w}{y_3}\right)}{(q^{-1} - q)^3} + \frac {\delta \left(\frac {z}{y_2}\right) \delta \left(\frac {z'}{y_1}\right) \delta \left(\frac {w}{y_3}\right)}{(q^{-1} - q)^3} \frac {\zeta \left( \frac {z'}z \right)}{\zeta \left( \frac z{z'} \right)} + \frac {\delta \left(\frac {z}{y_1}\right) \delta \left(\frac {z'}{y_3}\right) \delta \left(\frac {w}{y_2}\right)}{(q^{-1} - q)^3} \frac {\zeta \left( \frac w{z'} \right)}{\zeta \left( \frac {z'}w \right)} +
$$
\begin{equation}
\label{eqn:3}
+ \frac {\delta \left(\frac {z}{y_2}\right) \delta \left(\frac {z'}{y_3}\right) \delta \left(\frac {w}{y_1}\right)}{(q^{-1} - q)^3} \frac {\zeta \left( \frac wz \right)\zeta \left( \frac w{z'} \right)}{\zeta \left( \frac zw \right)\zeta \left( \frac {z'}w \right)} + \frac {\delta \left(\frac {z}{y_3}\right) \delta \left(\frac {z'}{y_1}\right) \delta \left(\frac {w}{y_2}\right)}{(q^{-1} - q)^3} \frac {\zeta \left( \frac {z'}z \right)\zeta \left( \frac wz \right)}{\zeta \left( \frac z{z'} \right)\zeta\left( \frac zw \right)} \qquad
\end{equation}
The terms in the right-hand sides of the expressions above are products of the form:
$$
\frac {\zeta \left(\frac ab \right)}{\zeta \left(\frac ba \right)}
$$
and any such product must be expanded as a power series in the range $|a| \gg |b|$ in order for relations \eqref{eqn:1}--\eqref{eqn:3} to hold. To prove that \eqref{eqn:serre} is preserved by the bialgebra pairing, one needs to show that:
$$
\Big(\text{RHS of }\eqref{eqn:1}\Big) - (q+q^{-1})\Big(\text{RHS of }\eqref{eqn:2}\Big) + \Big(\text{RHS of }\eqref{eqn:3}\Big) +
$$
\begin{equation}
\label{eqn:blade}
+\left\{\text{same expression with }z\text{ and }z' \text{ switched}\right \} = 0
\end{equation}
for each choice of $\{j_1,j_2,j_3\} = \{i,i,i+1\}$. One proves the above equality for each coefficient of the $\delta$ functions. For example, when $j_1=i+1$ and $j_2=j_3=i$, the coefficient of $\delta(\frac w{y_1})\delta(\frac z{y_2}) \delta(\frac {z'}{y_3})$ in the LHS of \eqref{eqn:blade} equals $(q^{-1} - q)^{-3}$ times:
$$
1 -(q+q^{-1}) \frac {\zeta \left( \frac wz \right)}{\zeta \left( \frac zw \right)} + \frac {\zeta \left( \frac wz \right)\zeta \left( \frac w{z'} \right)}{\zeta \left( \frac zw \right)\zeta \left( \frac {z'}w \right)} + \frac {\zeta \left( \frac z{z'} \right)}{\zeta \left( \frac {z'}z \right)} - (q+q^{-1})\frac {\zeta \left( \frac w{z'} \right) \zeta \left( \frac z{z'} \right)}{\zeta \left( \frac {z'}w \right)\zeta \left( \frac {z'}z \right)} +\frac {\zeta \left( \frac {z}{z'} \right)\zeta \left( \frac wz \right) \zeta \left( \frac w{z'} \right)}{\zeta \left( \frac {z'}{z} \right)\zeta \left( \frac {z'}w \right)\zeta \left( \frac zw \right)}
$$
\begin{align*}
&= 1 - (q+q^{-1})\frac {zq-wq^{-1}}{z-w} + \frac {(zq-wq^{-1})(z'q-wq^{-1})}{(z-w)(z'-w)} + \frac {z'-zq^2}{z'q^2-z} - \\ &- (q+q^{-1}) \frac {(z'q-wq^{-1})(z'-zq^2)}{(z'-w)(z'q^2-z)} + \frac {(z'-zq^2)(z'q-wq^{-1})(zq-wq^{-1})}{(z'q^2-z)(z'-w)(z-w)} = 0
\end{align*}
(since the equation above is an equality of rational functions, it is also an equality of power series, when expanded in the range $|w| \gg |z| \gg |z'|$). The analogous computations when $i+1$ is replaced by $i-1$, or when the signs $+$ and $-$ are switched, are proved in similar fashion. \\	
\end{proof}

\begin{proof} \emph{of Proposition \ref{prop:pairshuf}:} The fact that the pairing respects \eqref{eqn:comm0} is proved just like the analogous statement that \eqref{eqn:floyd} equals \eqref{eqn:lloyd} in the proof of Proposition \ref{prop:pairtor}. We will therefore focus on proving the fact that the pairing respects the shuffle product within $\CA^\pm$. Specifically, we will prove that the pairing \eqref{eqn:pairshuf} satisfies:
\begin{equation}
\label{eqn:lou}
\left \langle R_1^- * R_2^-, R^+ \right \rangle = \left \langle R_1^- \otimes R_2^-, \Delta(R^+) \right \rangle
\end{equation}
\begin{equation}
\label{eqn:reed}
\left \langle R^-, R_1^+ * R_2^+ \right \rangle = \left \langle \Delta(R^-), R_2^+ \otimes R_1^+ \right \rangle
\end{equation}
for all $R^\pm_1\in \CA_{\pm \bk^1}$, $R^\pm_2 \in \CA_{\pm \bk^2}$ and $R^\pm \in \CA_{\pm \bk}$, where $\bk = \bk^1+\bk^2$. Plugging in the formula for $\Delta(R^+)$ from \eqref{eqn:coproduct1} tells us that: 
$$
\text{RHS of \eqref{eqn:lou} } = \left \langle R^-_1(z_{ia}) \otimes R^-_2(w_{jb}), \frac {\prod_{j,b} \ph^+_{j}(w_{jb}) * R^+(z_{ia} \otimes w_{jb})}{\prod_{j,b}^{i,a} \zeta(w_{jb}/z_{ia})} \right \rangle \stackrel{\eqref{eqn:bialg 2}}=
$$
$$
= \left \langle \prod_{i,a} \ph^-_{i}(z_{ia}) \otimes R^-_1(z_{ia}) \otimes R^-_2(w_{jb}), \frac {\prod_{j,b} \ph^+_{j}(w_{jb}) \otimes R^+(z_{ia} \otimes w_{jb})}{\prod_{j,b}^{i,a} \zeta(w_{jb}/z_{ia})} \right \rangle = 
$$
$$
= \frac {(1-q^{-2})^{|\bk^1|+|\bk^2|}}{\bk^1! \cdot \bk^2!} \int^{*, |\alpha| \ll 1}_{|z_{ia}| = \alpha \oq^{- \frac {2i}n}, |w_{jb}| = \oq^{- \frac {2j}n}} \prod_{j,b}^{i,a} \Big \langle \ph^-_{i}(z_{ia}),  \ph^+_{j}(w_{jb}) \Big \rangle \cdot 
$$
$$
\frac {R^-_1(z_{ia}) R^-_2(w_{jb}) R^+(z_{ia},w_{jb}) \prod Dz_{ia} \prod Dw_{jb}}{\prod \zeta_p(z_{jb}/z_{ia})\prod \zeta_p(w_{jb}/w_{ia}) \prod \zeta(w_{jb}/z_{ia})} \Big |_{p \mapsto q}
$$
for any $\alpha \ll 1$. Use \eqref{eqn:ohm} to compute the pairing of $\ph^\pm$ in the above:
$$
\text{RHS of \eqref{eqn:lou} } = \frac {(1-q^{-2})^{|\bk|}}{\bk^1! \cdot \bk^2!} \int^{*, |\alpha| \ll 1}_{|z_{ia}| = \alpha \oq^{-\frac {2i}n}, |w_{jb}| = \oq^{-\frac {2j}n}} R^-_1(z_{ia}) R^-_2(w_{jb}) \prod_{j,b}^{i,a} \zeta \left( \frac {w_{jb}}{z_{ia}} \right)
$$
\begin{equation}
\label{eqn:worms}
\frac {R^+(z_{ia},w_{jb}) \prod Dz_{ia} \prod Dw_{jb}}{\prod \zeta_p(z_{jb}/z_{ia}) \prod \zeta_p(w_{jb}/w_{ia}) \prod \zeta_p(z_{ia}/w_{jb}) \prod \zeta_p(w_{jb}/z_{ia})}  \Big |_{p \mapsto q}
\end{equation}
Note that we have changed $\zeta$ to $\zeta_p$ in the last two products in the denominator, which is allowed since the assumption $|z_{ia}| \ll |w_{jb}|$ means that one can compute the integral without picking up any residues of the form $z_{ia} q - w_{ib} q^{-1}$ or $z_{ia} p - w_{ib} p^{-1}$, regardless of the relative sizes of $q$, 1 and $p$. One observes that the only poles of \eqref{eqn:worms} which involve both the $z_{ia}$ and the $w_{jb}$ are of the form:
\begin{equation}
\label{eqn:residues}
z_{i+1,a} q^{-1} - w_{ib} q \quad \qquad z_{ia} - w_{i+1,b} \quad \qquad z_{ia} p - w_{ib} p^{-1}
\end{equation}
There is also a simple pole at $z_{ia}p^{-1} - w_{ib}p$, but it is canceled by the zero of the numerator at $z_{ia}q^{-1} - w_{ib}q$, upon setting $p=q$ at the end. We observe that none of the poles \eqref{eqn:residues} hinder us from moving the contours from $|\alpha| \ll 1$ to $|\alpha| = 1$, because of the assumptions $|q| > 1 > |\oq|> |p|$ and $|pq|=1$. Since all $z$ and $w$ variables are now integrated over the same contours, we can symmetrize (i.e. average over all $\bk! = \prod_{i=1}^n k_i!$ permutations which preserve the colors of the variables) the integrand without changing the value of the integral:
$$
\text{RHS of \eqref{eqn:lou} } = \frac {(1-q^{-2})^{|\bk|}}{\bk!} \int^*_{|z_{ia}| = |w_{ib}| = \oq^{-\frac {2i}n}} \sym \left[ \frac {R^-_1(z_{ia})R^-_2(w_{jb})}{\bk^1!\cdot \bk^2!} \prod_{j,b}^{i,a} \zeta \left( \frac {w_{jb}}{z_{ia}} \right) \right]
$$
$$
\frac {R^+(z_{ia},w_{jb}) \prod Dz_{ia} \prod Dw_{jb}}{\prod \zeta_p(z_{jb}/z_{ia}) \prod \zeta_p(w_{jb}/w_{ia}) \prod \zeta_p(z_{ia}/w_{jb}) \prod \zeta_p(w_{jb}/z_{ia})}  \Big |_{p \mapsto q}
$$
Recall that the symmetrization on the first line is precisely the shuffle product $R_1^- * R_2^-$, and so the expression above matches the LHS of \eqref{eqn:lou}, as required. \\

\noindent Let us now prove \eqref{eqn:reed}. Formula \eqref{eqn:coproduct2} for $\Delta(R^-)$ implies that:
$$
\text{RHS of \eqref{eqn:reed} } = \left \langle \frac {R^-(z_{ia} \otimes w_{jb}) * \prod_{i,a} \ph^-_{i}(z_{ia})}{\prod_{j,b}^{i,a} \zeta(z_{ia}/w_{jb})}, R^+_2(z_{ia}) \otimes R^+_1(w_{jb}) \right \rangle = 
$$
$$
= \left \langle \frac {R^-(z_{ia} \otimes w_{jb}) \otimes \prod_{i,a} \ph^-_{i}(z_{ia})}{\prod_{j,b}^{i,a} \zeta(z_{ia}/w_{jb})}, R^+_2(z_{ia}) \otimes R^+_1(w_{jb}) \otimes 1 \right \rangle = \frac {(1-q^{-2})^{|\bk^1|+|\bk^2|}}{\bk^1! \cdot \bk^2!} \cdot
$$
$$
\int^{*, |\alpha| \ll 1}_{|z_{ia}| = \alpha \oq^{-\frac {2i}n}, |w_{jb}| = \oq^{- \frac {2j}n}} \frac {R^+_1(w_{jb}) R^+_2(z_{ia}) R^-(z_{ia},w_{jb}) \prod Dz_{ia} \prod Dw_{jb}}{\prod \zeta_p(z_{jb}/z_{ia})\prod \zeta_p(w_{jb}/w_{ia}) \prod \zeta(z_{ia}/w_{jb})}  \Big |_{p \mapsto q} =
$$
$$
= \frac {(1-q^{-2})^{|\bk|}}{\bk^1! \cdot \bk^2!} \int^{*, |\alpha| \ll 1}_{|z_{ia}| = \alpha \oq^{-\frac {2i}n}, |w_{jb}| = \oq^{-\frac {2j}n}} R^+_1(w_{jb}) R^+_2(z_{ia}) \prod_{j,b}^{i,a} \zeta \left( \frac {w_{jb}}{z_{ia}} \right) 
$$
$$
\frac {R^-(z_{ia},w_{jb}) \prod Dz_{ia} \prod Dw_{jb}}{\prod \zeta_p(z_{jb}/z_{ia}) \prod \zeta_p(w_{jb}/w_{ia}) \prod \zeta_p(z_{ia}/w_{jb}) \prod \zeta_p(w_{jb}/z_{ia})}  \Big |_{p \mapsto q}
$$
As before, we have changed $\zeta$ to $\zeta_p$ in the denominator because of the assumption $|z_{ia}| \ll |w_{jb}|$. Moreover, the only poles which involve both $z_{ia}$ and $w_{jb}$ are of the form \eqref{eqn:residues}, and they do not hinder us from moving the contours from $|\alpha| \ll 1$ to $|\alpha| = 1$. Since all $z$ and $w$ variables are now integrated over the same contours, we can symmetrize the integrand without changing the value of the integral:
$$
\text{RHS of \eqref{eqn:reed} } = \frac {(1-q^{-2})^{|\bk|}}{\bk!} \int^{|q| >  1 > |\oq|, |p|}_{|z_{ia}| = |w_{ib}| = \oq^{-\frac {2i}n}} \sym \left[ \frac {R^+_1(w_{jb})R^+_2(z_{ia})}{\bk^1!\cdot \bk^2!} \prod_{j,b}^{i,a} \zeta \left( \frac {w_{jb}}{z_{ia}} \right) \right]
$$
$$
\frac {R^-(z_{ia},w_{jb}) \prod Dz_{ia} \prod Dw_{jb}}{\prod \zeta_p(z_{jb}/z_{ia}) \prod \zeta_p(w_{jb}/w_{ia}) \prod \zeta_p(z_{ia}/w_{jb}) \prod \zeta_p(w_{jb}/z_{ia})}  \Big |_{p \mapsto q}
$$
Recalling that the symmetrization on the first line is precisely the shuffle product $R_1^+ * R_2^+$, we note that the above matches the LHS of \eqref{eqn:reed}, as we set out to prove. \\
\end{proof}

\begin{proof} \emph{of Lemma \ref{lem:min}:} This Lemma and its proof are a generalization of Lemma 2.14 from \cite{F}, along the lines of Proposition 2.4 of \cite{Shuf}. We will only prove the case when the sign is $+$, as the case of $-$ is analogous. Any shuffle element $R$ of bidegree $(\bk,d)$ is determined by the Laurent polynomial $r$ of \eqref{eqn:shuf}, which has total degree $d+\sum_{i=1}^n k_{i}k_{i-1} $. The finiteness of the limits \eqref{eqn:limit1} implies that for all $0 \leq \bl \leq \bk$:
\begin{equation}
\label{eqn:degree}
\deg_{z_{i1},...,z_{il_i}} r(z_{i1},...,z_{ik_i}) \leq \mu \sum_{i=1}^n l_i + \sum_{i=1}^n \left(k_i l_{i-1} + l_i k_{i-1} - l_i l_{i-1} \right)
\end{equation}
Let $A^\mu_{\bk,d}$ denote the set of symmetric Laurent polynomials $r$ which satisfy the wheel conditions \eqref{eqn:wheel} and the degree constraints \eqref{eqn:degree}. Then we need to prove the desired bound for the dimension of $A^\mu_{\bk,d}$. Consider any unordered set: 
$$
C = \Big\{ [i_1,j_1),...,[i_t,j_t) \Big\}
$$ 
of arcs $[i,j)$. If $\bk = \sum_{s=1}^t [i_s;j_s)$, then we will call such a set a \underline{partition} of $\bk$ and write $C \vdash \bk$. We order the constituent intervals from 1 to $t$ in descending order of length $j_s - i_s$, where those of the same length can be placed in any order. Given $\bk$ variables, which means $k_i$ variables of color $i$ for each $i\in \{1,...,n\}$, we can play the following game. Split up the set $\{z_{ia}\}$ into groups consisting of variables of colors $i_s,...,j_s-1$, and set the variables in each group equal to an indeterminate $y_s$. Thus, we assign to each variable an indeterminate, and we write this as $z_{ia} \leadsto y_s$. This evaluation gives rise to a linear operator:
\begin{equation}
\label{eqn:evalgarden}
\phi_{C}: A^\mu_{\bk,d} \longrightarrow \BF[y_1^{\pm 1},...,y_t^{\pm 1}], \qquad r(...,z_{ia},...) \mapsto r|_{z_{ia} \leadsto y_s} 
\end{equation}
Because $r$ is color-symmetric, the way we split up the $\bk$ variables among the intervals $[i_s;j_s)$ does not matter. The image of $\phi_C$ is a \underline{partially symmetric} polynomial in $y_1,...,y_t$, i.e. is symmetric in $y_s$ and $y_{s'}$ whenever $(i_s,j_s) - (i_{s'},j_{s'}) \in \BZ (n,n)$. We can filter the vector space $A^\mu_{\bk,d}$ by using the evaluation maps \eqref{eqn:evalgarden}:
$$
A^{\mu,C}_{\bk,d} := \bigcap_{C' > C} \text{Ker} \phi_{C'} 
$$ 
where $C'>C$ is the partial ordering on partitions, lexicographic in the lengths of the constituent intervals\footnote{Given two sets of natural numbers $X=\{x_1\geq x_2 \geq ...\}$ and $X' = \{x'_1\geq x'_2\geq ...\}$, we say that $X>X'$ in lexicographic order if there exists some $i$ such that $x_1=x'_1$, ..., $x_{i-1}=x'_{i-1}$ but $x_i>x'_i$. In our case, $C$ and $C'$ are sets of intervals $[i;j)$, so we consider $X = \{j-i|[i;j) \in C\}$ and $X' = \{j'-i'|[i';j') \in C'\}$, and define $C> C'$ if $X>X'$} (if the partition $C$ is maximal, then we set $A^{\mu,C}_{\bk,d} = A^\mu_{\bk,d}$). It is elementary to prove the fact that:
$$
\dim A^\mu_{\bk,d} \leq \sum_C \dim \phi_C(A^{\mu,C}_{\bk,d})	
$$
Then the conclusion of Lemma \ref{lem:min} follows from the fact that:
\begin{equation}
\label{eqn:claim2}
\dim \phi_C(A^{\mu,C}_{\bk,d}) \leq \# \text{ of tuples }(d_1,...,d_t) 
\end{equation}
such that conditions \eqref{eqn:cond1} and \eqref{eqn:cond2} hold. The above collections $(d_1,...,d_t)$ are partially ordered: we ignore the order of $d_s$ and $d_{s'}$ if $(i_s,j_s) - (i_{s'},j_{s'}) \in \BZ(n,n)$. Given $r\in A^{\mu,C}_{\bk,d}$, the Laurent polynomial $p = \phi_C(r)$ has total degree $d+\sum_{i=1}^n k_{i}k_{i-1}$, while the degree in each variable is controlled by \eqref{eqn:degree}:
\begin{equation}
\label{eqn:will do}
\deg_{y_s}(p) \leq \mu |\bl^s| + \sum_{i=1}^{n} \left(k_i l^s_{i-1} + l^s_i k_{i-1} - l^s_i l^s_{i-1} \right), \ \forall s\in \{1,...,t\}
\end{equation}
where $\bl^s : =[i_s;j_s)$. We claim that for all $s < s'$, the polynomial $p$ vanishes for: \\
	
\begin{itemize}
		
\item $y_{s'} = y_{s} q^{-2}$ \ for all $x \in [i_s,j_s-1)$ and $x'\in [i_{s'},j_{s'})$ with $x'\equiv x$ mod $n$ \\
		
\item $y_{s'} = y_{s} q^{2}$ \quad for all $x\in [i_s+1,j_s)$ and $x'\in [i_{s'},j_{s'})$ with $x'\equiv x$ mod $n$ \\
		
\item $y_{s'} = y_{s}$ \qquad for all $x' \in [i_{s'},j_{s'}) $ with $x' \equiv i_s-1$ mod $n$ \\
		
\item $y_{s'} = y_{s}$ \qquad for all $x' \in [i_{s'},j_{s'})$ with $x' \equiv j_s$ mod $n$ \\
		
\end{itemize} 
	
\noindent The above zeroes of $p$ are counted with the correct multiplicities. The first two bullets take place because of the wheel conditions \eqref{eqn:wheel}. As for the last two bullets, if we could set the variables equal to each other as prescribed therein, we could splice together the longer interval $[i_s;j_s)$ to part of the shorter interval $[i_{s'};j_{s'})$, thus obtaining a collection $C'$ which is larger than $C$ in lexicographic order. The assumption $r \in \cap_{C' > C} \ker \phi_{C'} $ implies that $p$ vanishes upon the specialization associated to the collection $C'$, and therefore also vanishes when we further specialize the variables to $y_{s'} = y_s$, as prescribed in the last two bullets. In conclusion, $p$ is divisible by the polynomial $p_0$ which is the product of the linear factors $y_{s'} - y_s q^{-2 \text{ or }0 \text{ or } 2}$ in the four bullets above. One sees that:
\begin{equation}
\label{eqn:degree2}
\deg_{y_s}(p_0) =  \sum_{s \neq s'} \sum_{i=1}^n \left( l^s_{i}l^{s'}_{i-1} +  l^s_{i-1}l^{s'}_{i} \right) = \sum_{i=1}^n \left( l^s_{i}k_{i-1} + l^s_{i-1} k_{i} - 2 l^s_{i-1} l^s_{i} \right)
\end{equation}
while the total degree of $p_0$ equals:
$$
\deg (p_0) = \sum_{i=1}^n \left(k_ik_{i-1} - \sum_{s=1}^t l^s_i l^s_{i-1} \right)
$$
We conclude that the Laurent polynomial $p/p_0$ has total degree:
$$
\sum_{s=1}^t \deg_{y_s} \left( \frac p{p_0} \right) = d + \sum_{i=1}^n \sum_{s=1}^t l^s_i l^s_{i-1}
$$
while \eqref{eqn:will do} and \eqref{eqn:degree2} imply that the degree in each variable is bounded by:
$$
\deg_{y_s} \left( \frac p{p_0} \right) \leq \mu|\bl^s| + \sum_{i=1}^n l^s_i l^s_{i-1}
$$
Therefore, the Laurent polynomial $p/p_0$ is a linear combination of the monomials:
\begin{equation}
\label{eqn:monomial}
\prod_{s=1}^t  y_s^{d_s + \sum_{i=1}^n l^s_i l^s_{i-1}}
\end{equation}
for various $d_s \leq \mu|\bl^s| = \mu (j_s-i_s)$. Since $p/p_0$ is symmetric in $y_s$ and $y_{s'}$ if $(i_s,j_s) - (i_{s'},j_{s'}) \in \BZ(n,n)$, this proves that $p/p_0$ belongs to a vector space of dimension equal to the RHS of \eqref{eqn:claim2}. This proves \eqref{eqn:claim2}, and with it, the Lemma. \\

\end{proof}

\begin{proof} \emph{of Proposition \ref{prop:sam}:} Let us deal with the case when the sign is $+$, and prove the statement about $\alpha$. The situations $+ \leftrightarrow -$ and $\alpha \leftrightarrow \beta$ are analogous, and we leave them to the interested reader. We will appeal to the proof of Lemma \ref{lem:min}, and recall the linear maps $\phi_C$ of \eqref{eqn:evalgarden}. Let $C_0$ be the finest partition of $\bk$, i.e. the partition into simple roots $\bs^i$. Because $\phi_{C_0}$ is the identity, we have $R = 0$ if and only if $\phi_{C_0}(r) = 0$, where $R$ and $r$ are connected by \eqref{eqn:shuf}. Then tautologically:
$$
R = 0 \qquad \Leftrightarrow \qquad \phi_C(r) = 0 \qquad \forall C \vdash \bk
$$
In the proof of Lemma \ref{lem:min}, we showed that the dimension of $\phi_C(A_{\bk,d}^{\mu,C})$ is at most the number of collections $(d_1,...,d_t)$ which satisfy \eqref{eqn:cond2} and $d = d_1+...+d_t$, where $[i_1;j_1),...,[i_t;j_t)$ are the constituent intervals of the partition $C$. But since: 
$$
R \in \CB_{\mu|\bk} = \CA_{\leq \mu|\bk,d} \qquad \text{where} \qquad d = \mu|\bk|
$$
the number of such collections is 1 if the partition $C \vdash \bk$ is $\mu-$integral, and 0 otherwise. We conclude that:
$$
R = 0 \qquad \Leftrightarrow \qquad \phi_C(r) = 0 \qquad \forall \mu-\text{integral } C \vdash \bk
$$ 
If the partition $C$ is $\mu-$integral, it was shown in the proof of Lemma \ref{lem:min} that $\phi_C(r)$ lies in the one-dimensional vector subspace spanned by (let $\bl^s = [i_s;j_s)$):
$$
p_C(y_1,...,y_t) = \Big[p_0 \text{ from the proof of Lemma \ref{lem:min}}\Big]\prod_{s=1}^t  y_s^{\mu|\bl^s| + \sum_{i=1}^n l^s_i l^s_{i-1}}
$$
Therefore, $\phi_C(r)$ is 0 if and only if it $\phi_C(r) <p_C$ in the limit $y_1 \gg ... \gg y_t$, by which we mean that either the $y_1$--degree of $\phi_C(r)$ is less than that of $p_C$, or if these degrees are equal then the joint $(y_1,y_2)$--degree of $\phi_C(r)$ is less than that of $p_C$, etc. From the proof of Lemma \ref{lem:min}, it is easy to see that for all $u \in \{1,...,t\}$, the joint $(y_1,...,y_u)$--degree of $p_C$ is equal to the right-hand side of \eqref{eqn:degree} with $l_i = l^1_i+...+l^u_i$. Therefore, we conclude that:
\begin{equation}
\label{eqn:cond above}
R = 0 \Leftrightarrow \text{leading term}_{y_1 \gg ... \gg y_t}  \text{ of } R(...,z_{ia},...)\Big |_{z_{ia} \leadsto y_s} = 0
\end{equation}
for all $\mu-$integral partitions $C \vdash \bk$, where: \\

\begin{itemize}
	
\item the phrase ``leading term" refers to the term of order $y_1^{\mu|\bl^1|}...y_t^{\mu|\bl^t|}$. \\

\item the meaning of $z_{ia} \leadsto y_s$ is that we specialize the variables as explained in Lemma \ref{lem:min}: associated to the partition $C$, we divide the set $\{z_{ia}\}$ into groups of variables of color $i_s,...,j_s-1$, and specialize all the variables in such a group to a common value $y_s$. \\

\end{itemize} 

\noindent Recalling the definition of the coproduct $\Delta_\mu$ from \eqref{eqn:cristian1} and of the linear maps $\alpha$ from \eqref{eqn:alpha}, then \eqref{eqn:cond above} is readily seen to be equivalent to \eqref{eqn:sam}, as required.

\end{proof}

\begin{proof} \emph{of Proposition \ref{prop:belong}:} For simplicity, let us only prove the case of $S_m^+$, since the cases of $S_m^-, T_m^+, T_m^-$ are analogous. Our first task is to show that $S^+_m$ is indeed a shuffle element. We can write $S^+_m$ in the form \eqref{eqn:shuf}, with the numerator:  
\begin{equation}
\label{eqn:yalla}
r(z_i,...,z_{j-1}) = \oq^{\frac{2\#}n} \cdot \sym \left[ \frac {m(z_i,...,z_{j-1})}{\left(1-\frac {z_{i+1}}{z_iq^2}\right)...\left(1-\frac {z_{j-1}}{z_{j-2}q^2}\right)} \right.
\end{equation}
$$
\left. \frac {\prod^{a \equiv b}_{a < b} (z_bq \oq^{\frac {2b}n}-z_aq^{-1} \oq^{\frac {2a}n})\prod^{a\equiv b+1}_{a < b} (z_b  \oq^{\frac {2(b+1)}n} -z_a \oq^{\frac {2a}n}) \prod_{a < b}^{a\equiv b-1} (z_aq \oq^{\frac {2(a+1)}n} - z_bq^{-1} \oq^{\frac {2a}n})}{\prod^{a \equiv b}_{a < b} (z_b \oq^{\frac {2b}n} - z_a \oq^{\frac {2a}n})} \right]
$$
where $\#$ is a certain integer that will not be important for us, and all the products go over all $i \leq a < b < j$. The RHS of \eqref{eqn:yalla} is indeed a Laurent polynomial, for the following two reasons: firstly, the factors $1-\frac {z_{a+1}}{z_{a}q^2}$ in the denominator are canceled by like factors in the numerator. Moreover, the simple pole at: 
$$
z_b \oq^{\frac {2b}n} - z_a \oq^{\frac {2a}n}
$$
for $a\equiv b$ disappears after we take the symmetrization, due to the convention \eqref{eqn:identify}. Therefore, all we need to do is to show that \eqref{eqn:yalla} vanishes at the specializations \eqref{eqn:wheel}. In other words, we need to show that $r$ vanishes when we set: 
\begin{equation}
\label{eqn:wheel 2}
z_a \oq^{\frac {2a}n} = w, \quad z_b \oq^{\frac {2b}n} = w q^{\pm 2}, \quad z_c \oq^{\frac {2c}n} = w
\end{equation}
for three indices that satisfy the congruences $a\equiv b \equiv c \pm 1$. Let us treat the case when the sign is $\pm = +$ and leave the opposite case as an exercise to the interested reader. In each of the summands in the symmetrization \eqref{eqn:yalla}, the indices $a,b,c$ are permuted. Depending on their order under this permutation, we fall into one or more of the following situations: \\

\begin{itemize}

\item If $b<a$, then the first product in the numerator vanishes. \\

\item If $a<c$, then the second product in the numerator vanishes. \\

\item If $c<b-1$, then the third product in the numerator vanishes. 

\end{itemize}

\text{} \\
Regardless of the relative order of $a \neq b \neq c \neq a$ under the aformentioned permutation, we will be in at least one of the above cases, so $r$ vanishes at the specialization \eqref{eqn:wheel 2}, thus establishing the wheel conditions. Note that if we added any more factors to the denominator of \eqref{eqn:yalla}, one of the bullets above might fail to be true and our argument would break down. 

\text{} \\
We will now prove that $S^+_m \in \text{Im }\Upsilon^+$, and to do so it is enough to assume the Laurent polynomial $m$ to be homogeneous. We will prove the statement by induction on $j-i$, where the base case $j-i=1$ is trivial. Note that for any $a\in \{i+1,...,j-1\}$ and any homogeneous Laurent polynomials $m_1(z_i,...,z_{a-1})$ and $m_2(z_{a},...z_{j-1})$,
$$
S^+_{m_2} * S^+_{m_1} = \sym \left[ \frac {m_1(z_i,...,z_{a-1})m_2(z_{a},...,z_{j-1}) \left(1 -  \frac {z_{a}}{z_{a-1}q^2} \right)}{\left(1 - \frac {z_{i+1}}{z_{i}q^2}  \right) ... \left(1 - \frac {z_{j-1}}{z_{j-2}q^2} \right)} \prod_{i\leq a < b < j} \zeta \left( \frac {z_b}{z_a} \right) \right]
$$
By the induction hypothesis, the above shuffle element lies in $\textrm{Im }\Upsilon^+$ for all $m_1, m_2$. Therefore, so does $S^+_m$ for all Laurent polynomials $m$ in the homogeneous ideal:
$$
(z_{j-1} - z_{j-2} q^2,...,z_{i+1} - z_iq^2) \in \BF[z_i^{\pm 1},...,z_{j-1}^{\pm 1}] 
$$
This ideal consists of all homogeneous polynomials such that $m(q^{2i},...,q^{2j-2})=0$. So in order to prove that $S^+_m \in \textrm{Im } \Upsilon^+$ for \underline{all} homogeneous Laurent polynomials $m$, it is enough to do so for a \underline{single} homogeneous polynomial $m$ of any given degree such that $m(q^{2i},....,q^{2j-2})\neq 0$. To this end, consider the shuffle element:
$$
z_i^{c_i} * ... * z_{j-1}^{c_{j-1}} = \sym\left[z_i^{c_i}...z_{j-1}^{c_{j-1}} \prod_{i\leq a < b < j} \zeta \left( \frac {z_a}{z_b} \right) \right] 
$$
It lies in $\textrm{Im }\Upsilon^+$ for all $c_i,...,c_{j-1} \in \BZ$, by the very definition of $\Upsilon^+$. A suitable linear combination of these elements allows us to add any number of linear factors to the right-hand side, while still obtaining a shuffle element of $\text{Im } \Upsilon^+$. In particular:
$$
R = \sym\left[z_i^{c_i}...z_{j-1}^{c_{j-1}} \prod_{a < b}^{a\equiv b} \frac {z_aq \oq^{\frac {2a}n} - z_bq^{-1} \oq^{\frac {2b}n}}{z_a \oq^{\frac {2a}n} - z_b \oq^{\frac {2b}n}} \prod_{a < b}^{a\equiv b-1} \frac {z_a \oq^{\frac {2(a+1)}n} - z_b \oq^{\frac {2b}n}}{z_aq \oq^{\frac {2(a+1)}n} - z_bq^{-1} \oq^{\frac {2b}n}} \right.
$$
$$
\left. \prod_{a < b}^{a\equiv b} \left( z_bq \oq^{\frac {2b}n} - z_aq^{-1} \oq^{\frac {2a}n} \right)  \prod_{a < b}^{a\equiv b+1} \left( z_b \oq^{\frac {2(b+1)}n} - z_a \oq^{\frac {2a}n} \right) \prod_{a < b-1}^{a\equiv b-1} \left( z_bq^{-1} \oq^{\frac {2b}n} - z_aq \oq^{\frac {2(a+1)}n} \right) \right]  
$$
lies in $\textrm{Im }\Upsilon^+$. $R$ is a shuffle element, which can be written in the form \eqref{eqn:shuf}, with:
$$
r = \oq^{\frac {2\#}n} \cdot \sym\left[z_i^{c_i}...z_{j-1}^{c_{j-1}}  \prod_{a < b}^{a\equiv b} \frac {z_bq^{-1} \oq^{\frac {2b}n}- z_aq\oq^{\frac {2a}n}}{z_b \oq^{\frac {2b}n} - z_a \oq^{\frac {2a}n}}\prod_{a < b}^{a\equiv b-1} \left(z_a \oq^{\frac {2(a+1)}n} - z_b \oq^{\frac {2b}n} \right)   \right.
$$
$$
\prod_{a < b}^{a\equiv b+1}\left(z_b q \oq^{\frac {2(b+1)}n} - z_aq^{-1} \oq^{\frac {2a}n} \right) \prod_{a < b}^{a\equiv b} \left( z_bq \oq^{\frac {2b}n} - z_aq^{-1} \oq^{\frac {2a}n} \right)  
$$
$$
\left. \prod_{a < b}^{a\equiv b+1} \left( z_b \oq^{\frac {2(b+1)}n} - z_a \oq^{\frac {2a}n} \right)  \prod_{a < b-1}^{a\equiv b-1} \left( z_bq^{-1} \oq^{\frac {2b}n} - z_aq \oq^{\frac {2(a+1)}n} \right) \right]
$$
for some $\# \in \BZ$. The above symmetrization is precisely of the form \eqref{eqn:yalla}, with:
$$
m (z_i,...,z_{j-1})= \pm q^{\#'} \oq^{\frac {2\#''}n} \cdot z_i^{c_i-1}...z_{j-2}^{c_{j-2}-1} z_{j-1}^{c_{j-1}} \cdot
$$
$$
\prod_{a < b}^{a\equiv b} \left(z_bq^{-1} \oq^{\frac {2b}n}- z_aq\oq^{\frac {2a}n} \right) \prod_{a < b}^{a\equiv b-1} \left(z_a \oq^{\frac {2(a+1)}n} - z_b \oq^{\frac {2b}n} \right)   \prod_{a < b}^{a\equiv b+1}\left(z_b q \oq^{\frac {2(b+1)}n} - z_aq^{-1} \oq^{\frac {2a}n} \right)
$$
for some $\pm \in \{+,-\}$ and some $\#',\#'' \in \BZ$. Since $m(q^{2i},...,q^{2j-2}) \neq 0$, and one can choose the exponents $c_i,...,c_{j-1}$ arbitrarily, this concludes our proof. \\
\end{proof} 

\begin{proof} \emph{of Proposition \ref{prop:tommy}:} We will prove the required statements for $E_{[i;j)}^\mu$, and leave the analogous case of $F_{[i;j)}^\mu$ as an exercise to the interested reader. Set:
$$
R := E_{[i;j)}^\mu \qquad \text{and assume that} \qquad d := \mu(j-i) \in \BZ
$$
otherwise the Proposition is vacuous. It is clear that $\homdeg R = d$, so in order to prove that $R \in \CB_\mu^+$, we must show that $R$ has slope $\leq \mu$. Looking back to \eqref{eqn:limit1}, we must show that for all subsets $A \subset \{i,...,j-1\}$, we have:
\begin{equation}
\label{eqn:want}
\text{total degree of }\{z_a\}_{a\in A} \text{ in } R(z_i,...,z_{j-1}) \leq \mu|A|
\end{equation}
Looking back at formula \eqref{eqn:min1} for $R$, and recalling the fact that $\zeta$ has finite limit $q$, 1 or $q^{-1}$ when one of the variables are sent off to $\infty$, we conclude that: 	
$$
\text{total degree of }\{z_a\}_{a\in A} \text{ in } R = \sum_{a\in A} \Big( \lfloor \mu(a-i+1) \rfloor - \lfloor \mu(a-i) \rfloor \Big) - \sum_{a\in A}^{a-1 \notin A} 1
$$	
Let us assume that $A$ can be partitioned into distinct blocks of consecutive integers:
$$
A = \{i_1,...,j_1-1\} \sqcup \{i_2,...,j_2-1\} \sqcup ... \sqcup \{i_t,...,j_t-1\}
$$
where $i\leq i_1<j_1<i_2<j_2<...<i_t<j_t \leq j$. Then we have:
$$
\text{total degree of }\{z_a\}_{a\in A} \text{ in } R = \sum_{k=1}^t \Big( \lfloor \mu(j_k-i)\rfloor - \lfloor \mu(i_k-i) \rfloor \Big) - t+\delta_{i_1}^i
$$	
Then \eqref{eqn:want} follows from the general inequality:
\begin{equation}
\label{eqn:melcs}
\lfloor x \rfloor - \lfloor y \rfloor - 1 + \delta_{y}^0 \leq x-y
\end{equation}	
Moreover, since one has equality in \eqref{eqn:melcs} only when $x\in \BZ$ and $y = 0$, we conclude that one has equality in \eqref{eqn:want} only when $t=1$ and $i_1 = i$, i.e. when:
$$
A = \{i,...,a-1\} \quad \text{for some } a \text{ such that } \mu(a-i) \in \BZ
$$
These are precisely the terms which survive in the coproduct $\Delta_\mu$ from \eqref{eqn:cristian1}: the limit is non-zero if and only if the ``small" variables (those to the left of $\otimes$) are precisely $z_a...,z_{j-1}$ and the ``large" variables (those to the right of $\otimes$) are precisely $z_i,...,z_{a-1}$. Taking the limit, one obtains precisely formula \eqref{eqn:cop1} for $\Delta_\mu(R)$. \\	
\end{proof} 

\begin{proof} \emph{of Proposition \ref{prop:pairing}:} We will use formula \eqref{eqn:pairshuf} to compute the pairing: 
$$
\langle F_{[i;j)}^\mu, R^+ \rangle \ = \ (1 - q^{-2})^{j-i} \int^*_{|z_a|=\oq^{-\frac {2a}n}}  \frac {\prod_{a=i}^{j-1} z_a^{\lfloor \mu(a-i) \rfloor - \lfloor \mu(a-i+1) \rfloor}}{\left(1 - \frac {z_i}{z_{i+1}}  \right) ... \left(1 - \frac {z_{j-2}}{z_{j-1}} \right)} \cdot
$$
\begin{equation}
\label{eqn:cc}
\cdot \frac {\prod_{i\leq a < b < j} \zeta \left( \frac {z_a}{z_b} \right) R^+(z_i,...,z_{j-1})}{\prod_{i \leq a \neq b < j} \zeta_p \left( \frac {z_{a}}{z_{b}} \right)} \prod_{a=i}^{j-1} Dz_{a} \Big |_{p \mapsto q}
\end{equation}
We have foregone the symmetrization because all the variables $z_a$ with a given $a \text{ mod }n$ are integrated over the same contour. The poles of the integrand are: \\

\begin{itemize}

\item $z_a = z_{a+1}$ \\

\item $z_a = z_b \cdot q^{-2}\oq^{\frac {2(b-a-1)}n}$ for $\col a \equiv \col b - 1$ and $a < b$ \\

\item $z_a = z_b \cdot \oq^{\frac {2(b-a+1)}n}$ for $\col a \equiv \col b + 1$ and $a < b$ \\

\item $z_a = z_b \cdot p^2 \oq^{\frac {2(b-a)}n}$ for $\col a \equiv \col b$ and $a<b$ \\

\end{itemize}

\noindent The apparent pole at $z_a = z_b p^2 \oq^{\frac {2(b-a)}n}$ for $a>b$ vanishes, due to $z_a - z_b q^2 \oq^{\frac {2(b-a)}n}$ in the numerator (since after evaluating the integral, we set $p = q$). Because of the assumption on the sizes of $q,\oq,p$ that we made in Subsection \ref{sub:full}, we observe that the only poles that show up as we move the contours to $|z_i| \gg ... \gg |z_{j-1}|$ are the ones in the first bullet, namely $z_a = z_{a+1}$. In other words, the residues one picks up are of the form:
$$
(z_i,...,z_{j-1}) = \Big(\underbrace{y_1,...,y_1}_{k_1 \text{ terms}}, \underbrace{y_2,..., y_2}_{k_2 \text{ terms}},...,\underbrace{y_t,...,y_t}_{k_t \text{ terms}} \Big) 
$$
for certain natural numbers $k_1+...+k_t = j-i$ and $|y_1| \gg ... \gg |y_t|$. Because $R^+ \in \CB^+_\mu$, then the integrand of \eqref{eqn:cc} has order:
$$
\leq - \lfloor \mu k_1 \rfloor - 1 + \delta_t^1 + \mu k_1
$$
in the variable $y_1$. The quantity above is $<0$ and so produces a zero residue when integrated around $\infty$, unless $t = 1$ (which implies $\mu k_1 = d \in \BZ$). We conclude that the only residue which contributes to \eqref{eqn:cc} non-trivially is $z_a = 1$, $\forall a\in \{i,...,j-1\}$, hence:
$$
\langle F_{[i;j)}^\mu, R^+ \rangle = (1 - q^{-2})^{j-i} \cdot \frac {R^+(z_i,...,z_{j-1})}{\prod_{i \leq a < b < j} \zeta(z_{b}/z_{a})} \Big|_{z_a \mapsto 1}
$$
Recalling the definition of $\alpha$ in \eqref{eqn:alpha}, one obtains precisely \eqref{eqn:bonnie1}. Meanwhile:
$$
\langle R^- , E_{[i;j)}^\mu \rangle \ = \ (1 - q^{-2})^{j-i} \int^*_{|z_a|=\oq^{-\frac {2a}n}}  \frac {\prod_{a=i}^{j-1} z_a^{\lfloor \mu(a-i+1) \rfloor - \lfloor \mu(a-i) \rfloor}}{\left(1 - \frac {z_{i+1}}{z_{i}q^2}  \right) ... \left(1 - \frac {z_{j-1}}{z_{j-2}q^2} \right)} \cdot
$$
\begin{equation}
\label{eqn:catch}
\cdot \frac {\prod_{i\leq a < b < j} \zeta \left( \frac {z_b}{z_a} \right) R^-(z_i,...,z_{j-1})}{\prod_{i \leq a \neq b < j} \zeta_p \left( \frac {z_{a}}{z_{b}} \right)} \prod_{a=i}^{j-1} Dz_{a} \Big |_{p \mapsto q}
\end{equation}
We have foregone the symmetrization because all the variables $z_a$ with a given $a \text{ mod }n$ are integrated over the same contour. As before, because of the assumption on the sizes of $q,\oq,p$, as we move the contours to $|z_i| \ll ... \ll |z_{j-1}|$, the only poles we pick up arise from the factors:
$$
1 - \frac {z_{a+1}}{z_aq^2} \qquad \Rightarrow \qquad \text{poles at }z_{a+1} = z_{a}q^2
$$
In other words, the residues one picks up are of the form:
$$
(z_i,...,z_{j-1}) = \Big(y_1q^{2k_0},...,y_1q^{2k_1-2}, y_2q^{2k_1},..., y_2q^{2k_2-2} ,...,y_tq^{2k_{t-1}},...,y_tq^{2k_t-2} \Big) 
$$
for certain $i = k_0 < k_1 < ... < k_t = j$ and $|y_1| \ll ... \ll |y_t|$. Because $R^- \in \CB^-_\mu$, then property \eqref{eqn:limit2} implies that the integrand of \eqref{eqn:catch} has order:
$$
\leq d - \lfloor \mu(k_{t-1}-i) \rfloor - 1 + \delta_t^1 - \mu(j - k_{t-1})
$$
in the variable $y_t$. Because $d = \mu(j-i) \in \BZ$, the above quantity is $<0$ and so produces a zero residue when integrated around $\infty$, unless $t = 1$. We conclude that the only residue which contributes to \eqref{eqn:catch} is $z_a = q^{2a}$, $\forall a\in \{i,...,j-1\}$, hence:
$$
\langle R^-, E_{[i;j)}^\mu \rangle = (1 - q^{-2})^{j-i} \prod_{a=i}^{j-1} q^{2a \left(\lfloor \mu(a-i+1) \rfloor - \lfloor \mu(a-i) \rfloor \right)} \cdot \frac {R^-(z_i,...,z_{j-1})}{\prod_{i \leq a < b < j} \zeta(z_a / z_b)} \Big|_{z_a \mapsto q^{2a}}
$$
Recalling the definition of $\beta$ in \eqref{eqn:beta}, one obtains precisely \eqref{eqn:bonnie2}. Relation \eqref{eqn:tyler} follows from \eqref{eqn:bonnie1} and the fact that:
$$
\alpha_{[i';j')} \left(\sym \left[ \frac {\prod_{a=i}^{j-1} z_a^{\lfloor \mu(a-i+1) \rfloor - \lfloor \mu(a-i) \rfloor}}{\left(1 - \frac {z_{i+1}}{z_{i}q^2}  \right) ... \left(1 - \frac {z_{j-1}}{z_{j-2}q^2} \right)} \prod_{i\leq a < b < j} \zeta \left( \frac {z_b}{z_a} \right)  \right] \right) = \frac {\delta_{[i';j')}^{[i;j)}}{(1-q^{-2})^{j-i-1}}
$$
The equality above follows easily from \eqref{eqn:alpha} and the fact that $\zeta(z)|_{z \mapsto 1} = 0$ if $\col z = - 1$. Indeed, this means that the only summand of the symmetrization which survives the evaluation $z_i = ... = z_{j-1} = 1$ is the identity permutation. Finally, let us prove by induction on $\bk \in \nn$ that if:
$$
R^+ \in \CB_{\mu|\bk}
$$
pairs trivially with all products of $F_{[i;j)}^\mu$, then $R^+ = 0$. Note that \eqref{eqn:bonnie1} implies:
$$
\text{LHS of \eqref{eqn:sam} } = (1-q^{-2})^{-|\bk|} \left \langle F_{[i_1;j_1)}^\mu \otimes ... \otimes F_{[i_t;j_t)}^\mu, \Delta_\mu^{(t-1)}(R^+) \right \rangle \stackrel{\eqref{eqn:bialg 1}}=
$$
$$
= (1-q^{-2})^{-|\bk|} \left \langle F_{[i_1;j_1)}^\mu * ... * F_{[i_t;j_t)}^\mu, R^+ \right \rangle = 0 \quad \stackrel{\text{Proposition \ref{prop:sam}}}\Longrightarrow \quad R^+ = 0
$$

\end{proof}

\begin{proof} \emph{of Proposition \ref{prop:sasha}:} Without loss of generality, we will prove \eqref{eqn:sasha 1} and \eqref{eqn:sasha 2} only when the sign is $\pm = +$. We proceed by induction on $k+k'$. Write:
$$
R := [G_k, G_{k'}] \in \CB_0^+
$$
and let us compute:
\begin{align*}
\Delta_0(R) &= \Delta_0(G_k)\Delta_0(G_{k'}) - \Delta_0(G_{k'})\Delta_0(G_k) = \\
&= \sum^{a+b=k}_{a'+b'=k'}  \left( G_aG_{a'}c^{b+b'} \otimes G_b G_{b'} - G_{a'} G_a c^{b+b'} \otimes G_{b'} G_b\right)
\end{align*}
By the induction hypothesis, $[G_a,G_{a'}] = 0$ whenever $a+a' < k+k'$, so the only terms which survive in the coproduct are:
$$
\Delta_0(R) = c^{k+k'} \otimes G_{k}G_{k'} - c^{k+k'} \otimes G_{k'} G_k + G_kG_{k'} \otimes 1 - G_{k'}G_k \otimes 1 = c^{k+k'} \otimes R + R \otimes 1
$$
We conclude that $R$ is a primitive element of $\CB_0^+$. Its images under the linear maps $\alpha_{[i;j)}$ can be computed by using Lemma \ref{lem:pseudo}:
$$
\alpha_{[i;i+n(k+k'))}(R) = \alpha_{[i;i+nk)}(G_k) \alpha_{[i;i+nk')}(G_{k'}) - \alpha_{[i;i+nk')}(G_{k'})  \alpha_{[i;i+nk)}(G_k) = 0
$$ 
for all $i\in \{1,...,n\}$. Since all the $G_k$ are $\BZ/n\BZ$--invariant, we conclude that $R$ also is $\BZ/n\BZ$--invariant. Then Corollary \ref{cor:sam} implies that $R = 0$, as expected in \eqref{eqn:sasha 1}. \\

\noindent To prove \eqref{eqn:sasha 2}, let us consider $R' = [G_k, 1^+_i] \in \CB_0^+$, and so we may compute:
\begin{align*}
\Delta_0(R') &= \Delta_0(G_k)\Delta_0(1^+_i) - \Delta_0(1^+_i)\Delta_0(G_k) = \\ 
&= \sum_{a+b=k} \left( G_a c^{b}  \ph_i \otimes G_b 1^+_i + G_a c^{b} 1^+_i \otimes G_b - \ph_i  G_a c^{b} \otimes 1^+_i G_b -  1^+_i G_a c^{b}  \otimes G_b  \right)
\end{align*}
By the induction hypothesis, $[G_a,1^+_i] = 0$ for all $a<k$. Therefore, we infer that:
$$
\Delta_0(R') = \ph_i c^k \otimes [G_k,1^+_i] + [G_k,1^+_i] \otimes 1 = \ph_ic^k \otimes R' + R' \otimes 1
$$
We conclude that $R'$ is a primitive element of $\CB_0^+$. It has degree $[i;i+nk+1)$, so its image under the linear maps $\alpha_{[i;i+nk+1)}$ can be computed by using Lemma \ref{lem:pseudo}:
$$
\alpha_{[i;i+nk+1)}(R') = \alpha_{[i;i+nk)}(G_k) \cdot 1 - 1 \cdot \alpha_{[i+1;i+1+nk)}(G_k) = 0
$$
In the first equality, we used the fact that $\alpha_{[j;j+1)}(1^+_i) = \delta_j^i$. In the second equality, we used the $\BZ/n\BZ-$invariance of $G_k$ to conclude that $\alpha_{[j;j+nk)}(G_k)$ does not depend on $j$. By Corollary \ref{cor:sam}, we obtain $R' = 0$, as expected in \eqref{eqn:sasha 2}. \\

\noindent It is possible to compute the pairings $\langle G_{-k}, G_k \rangle$ using formula \eqref{eqn:pairshuf}, and we obtain:
\begin{multline}
\langle G_{-k}, G_k \rangle = \frac {1}{k!^n} \int^{|q| > 1 > |\oq| > |q|^{-1}}_{|z_{ia}|=\oq^{- \frac {2i}n}}  \label{eqn:pair group} \\ \prod_{i=1}^n \left[  \frac {\prod_{1\leq a \neq b \leq k} \left( 1 - \frac {z_{ib}}{z_{ia}} \right)}{\prod_{1\leq a,b \leq k} \left( 1 - \frac {z_{i-1,b}}{z_{ia}} \right)} \frac {\prod_{1\leq a, b \leq k} \left(1 - \frac {z_{ib}}{z_{ia}q^2} \right)}{\prod_{1\leq a,b \leq k} \left(1 - \frac {z_{i+1,b}}{z_{ia}q^2} \right)} \prod_{1\leq a \leq k} Dz_{ia}  \right] 
\end{multline}
(we may replace the rational function $\zeta_p$ by $\zeta$ in \eqref{eqn:pairshuf}, because plugging in $G_{\pm k}$ for $R^\pm$ annihilates any pole of the form $z_{ia} q - z_{ib} q^{-1}$). Let us compute the integral above by the following procedure: take the variable $z_{01}$, rename it $s$, and let us move the $s$ contour toward 0. As we do so, we encounter two kinds of poles: \\

\begin{itemize} 
	
\item $s = z_{-1a}$ for some $a$ \\

\item $sq^2 = z_{1b}$ for some $b$ \\

\end{itemize}

\noindent It is easy to see that there is no pole at $s=0$, because $n>1$. For the first (resp. second) type of pole, we may rename the variable $z_{-1a}$ (resp. $z_{1b}$) by the symbol $s$ (resp. $sq^2$) and move the $s$ contour toward 0. Repeating the procedure entails a totally ordered set of variables $V \subset \{z_{ia}\}^{1\leq i \leq n}_{1\leq a \leq k}$ being specialized to $s$ times various powers of $q$ (we will refer to this as the $V$--specialization of \eqref{eqn:pair group}). If we let $\bl = (l_1,...,l_n) \in \nn$, where $l_i$ is the number of variables of color $i$ in the set $V$, then we claim the degree of the $V$--specialization at $s=0$ is equal to: 
$$
- \langle \bl, k\bde - \bl\rangle - \langle k\bde - \bl, \bl \rangle \geq 0
$$
The $V$--specialization yields a non-zero contribution to \eqref{eqn:pair group} only when the degree above is 0, so that we have a non-trivial residue at $s = 0$, and this happens precisely when $\bl = l\bde$ for some $l \in \{1,...,k\}$. To summarize, the discussion above implies:
\begin{equation}
\label{eqn:induction}
\langle G_{-k}, G_k \rangle = \frac {1}{k} \sum_{l=1}^k f_l \cdot \underbrace{\langle G_{-k+l}, G_{k-l} \rangle}_{\text{coefficient of }s^0\text{ in the }V\text{--specialization of \eqref{eqn:pair group}}}
\end{equation}
(the denominator $k$ right before the summation sign is there because we fixed a choice of starting variable $z_{01}$), where $f_k$ denotes the sum of all $V$--specializations of \eqref{eqn:pair group}, where $V$ coincides with the full set of variables. \\

\begin{claim}
\label{claim:rea}
	
For any $k$, we have $f_k = \frac {\oq^{2k}(1-q^{2nk})}{(1-\oq^{2k})(1-q^{2nk} \oq^{2k})}$. \\

\end{claim} 

\noindent Let us show how Claim \ref{claim:rea} implies the proof of the Proposition. Let:
$$
\Gamma(x) = \sum_{k=0}^\infty \langle G_{-k}, G_k \rangle x^k
$$
As easy computation in the Heisenberg algebra (which we leave as an exercise to the interested reader) shows that:
\begin{equation}
\label{eqn:u}
\Gamma(x) = \exp \left( \sum_{k=1}^\infty \frac {\langle P_{-k}, P_k \rangle x^k}{k^2} \right) 
\end{equation}
where $P_{\pm k}$ are defined in \eqref{eqn:new heis}. It follows from property \eqref{eqn:induction} implies that:
$$
x\Gamma(x)' = \left(\sum_{k=1}^\infty f_k x^k \right) \Gamma(x) \quad \Rightarrow \quad x \ln(\Gamma(x))' = \sum_{k=1}^\infty f_k x^k
$$
which implies \eqref{eqn:new pair heis}. \\

\begin{proof} \emph{of Claim \ref{claim:rea}:} As we compute the integral \eqref{eqn:pair group} by residues, the succession of poles entails specializing the variables $\{z_{11},...,z_{nk}\}$ to $s$ times various powers of $q$. Moreover, any such a specialization remembers the order in which the variables were encountered, for example the first, second, ..., last variables are:
\begin{equation}
\label{eqn:sequence}
s q^{2\alpha_1} \text{ of color }i_1, s q^{2\alpha_2} \text{ of color }i_2, ..., sq^{2\alpha_{nk}} \text{ of color }i_{n k}
\end{equation}
where $\alpha_1 = 0$, $i_1 = 0$. Because of the assumptions $|q| > 1 > |\oq| > |q|^{-1}$, we have:
$$
\alpha_{t+1} > \alpha_t \qquad \text{or} \qquad \alpha_{t+1} = \alpha_t \text{ and } i_{t+1} = i_t - 1
$$
for all $t$, in any specialization \eqref{eqn:sequence}. Moreover, an analysis of the zeroes and poles of the linear factors in \eqref{eqn:pair group} shows the following: as we specialize the variables in the integral to \eqref{eqn:sequence}, we only encounter simple poles, and we do so only if the sequence reads:
$$
sq^0 \text{ of color } 0, sq^0 \text{ of color } -1, ..., sq^0 \text{ of color } - \lambda_0 + 1 
$$
$$
sq^2 \text{ of color } 1, sq^2 \text{ of color } \ \ \ 0 ,..., sq^2 \text{ of color } - \lambda_1 + 2 
$$
$$
sq^4 \text{ of color } 2, sq^2 \text{ of color } \ \ \ 1 ,..., sq^2 \text{ of color } - \lambda_2 + 3 
$$
\begin{equation}
\label{eqn:...}
...
\end{equation}
where $\lambda_0 \geq \lambda_1 \geq ...$. It is clear that such collections are in one-to-one correspondence with partitions $\lambda = (\lambda_0, \lambda_1,...)$ of the number $nk$, whose Young diagrams contain exactly $k$ boxes with content $i$ modulo $n$, for any $i \in \{1,...,n\}$. \footnote{Recall that the content of a box at coordinates $(x,y)$ in a Young diagram is the number $x-y$} Furthermore, such partitions are in one-to-one correspondence with the $\BC^* \times \BC^*$ fixed points of the Nakajima quiver variety $\CN_{n,k}$ corresponding to the $n$ vertex cyclic quiver, corresponding to the degree vectors ${\mathbf{v}} = (k,...,k)$ and ${\mathbf{w}} = (1,0,...,0)$, and the dominant stability condition. It is well-known that:
$$
\CN_{n,k} \text{ is a connected component of } \text{Hilb}_{nk}(\BA^2)^{\BZ/n\BZ}
$$
where the $\BZ/n\BZ$ action on the Hilbert scheme of points is induced by the anti-diagonal action of $\BZ/n\BZ$ on the affine plane $\BA^2$. The tautological rank $nk$ vector bundle $\CT$ on the Hilbert scheme has a trivial summand:
$$
\CT \cong \tilde{\CT} \oplus \CO 
$$
and the isomorphism above above descends to the $\BZ/n\BZ$--invariant part, which is a rank $k$ vector bundle on the connected component $\CN_{n,k}$:
$$
\CT^{\BZ/n\BZ} \cong \tilde{\CT}^{\BZ/n\BZ} \oplus \CO 
$$
In virtue of the well-known formulas for the equivariant weights of the $\BC^* \times \BC^*$ fixed points of $\CN_{n,k}$ and the tautological vector bundle $\CT$ (we refer to \cite{thesis} for the formulas), the $K$--theoretic equivariant localization formula implies that:
\begin{equation}
\label{eqn:whoa}
f_k = \chi_{\BC^* \times \BC^*} \left(\CN_{n,k}, \wedge^\bullet(\tilde{\CT}^{\BZ/n\BZ} ) \right)
\end{equation}
where the equivariant parameters of $\BC^* \times \BC^*$ are set to $q_1 = \oq^{-\frac 2n}$ and $q_2 = q^2 \oq^{\frac 2n}$. The fact that the right-hand side of \eqref{eqn:whoa} is given by:
$$
\chi_{\BC^* \times \BC^*} \left(\CN_{n,k}, \wedge^\bullet(\tilde{\CT}^{\BZ/n\BZ} ) \right) = \frac {(1 - q_1^{-nk}q_2^{-nk})}{(1 - q_1^{-nk})(1 - q_2^{-nk})}
$$
is well-known, thus concluding the proof of the Claim. 
\end{proof}

\end{proof}

\begin{proposition} We have the following equalities in $\CB_0 \cong \uu$:
\begin{equation}
\label{eqn:pairing 1}
\langle F_{[i;i+nk)}, G_k \rangle = (-q^{-1})^{nk} \prod_{l=1}^k \frac {q \oq^{1-l}- q^{-1} \oq^{l-1}}{\oq^{l}-\oq^{-l}}
\end{equation}
\begin{equation}
\label{eqn:pairing 2}
\langle G_{-k}, E_{[i;i+nk)} \rangle = \prod_{l=1}^k \frac {q (\oq q^n)^{1-l}- q^{-1} (\oq q^n)^{l-1}}{(\oq q^n)^{l}-(\oq q^n)^{-l}}
\end{equation}
for all $i \in \{1,...,n\}$. \\

\end{proposition} 

\noindent As a consequence of \eqref{eqn:pairing 1} and \eqref{eqn:pairing 2}, and the multiplicativity of pairing with $E$'s and $F$'s (see Lemma \ref{lem:pseudo}), we conclude that:
\begin{equation}
\label{eqn:pairing 3}
\langle F_{[i;i+nk)}, P_k \rangle = (-q^{-1})^{nk} (-1)^{k-1} \frac {q^k-q^{-k}}{\oq^k-\oq^{-k}}
\end{equation}
\begin{equation}
\label{eqn:pairing 4}
\langle P_{-k}, E_{[i;i+nk)} \rangle = (-1)^{k-1} \frac {q^k-q^{-k}}{\oq^k q^{nk} - \oq^{-k} q^{-nk}}
\end{equation}
Indeed, the right-hand sides of formulas \eqref{eqn:pairing 3}--\eqref{eqn:pairing 4} are (when combined into a generating series by multiplying with $z^k$ and summing over all $k$) is simply the logarithm of the right-hand sides of formulas \eqref{eqn:pairing 1}, \eqref{eqn:pairing 2}, respectively. \\

\begin{proof} As a consequence of \eqref{eqn:alpha}, we have:
$$
\alpha_{[i;i+nk)}(G_k) = G_k(1,...,1,\oq^{-2},...,\oq^{-2},...,\oq^{-2k+2},...,\oq^{-2k+2}) \prod_{l=1}^k \frac {q \oq^{2l} - q^{-1}}{\oq^{2l}-1}
$$
where $G_k$ is thought of having $k$ variables of each color $i,i+1,...,i+n-1$. Plugging in the explicit formula \eqref{eqn:g1}, we have:
$$
\alpha_{[i;i+nk)}(G_k) = \frac {\oq^{k^2}}{(q^{-1}-q)^{nk}} \prod_{l=1}^k \frac {q \oq^{2-2l} - q^{-1}}{\oq^{2l}-1} \quad \stackrel{\eqref{eqn:bonnie1}}\Longrightarrow \quad \eqref{eqn:pairing 1}
$$
Similarly, as a consequence of \eqref{eqn:beta} (with $\mu = 0$) we have:
\begin{multline*}
\beta_{[i;i+nk)}(G_{-k}) = G_{-k}(q^{2j},...,q^{2j}\oq^{-2k+2})_{j \in \{i,...,i+n-1\}} \\ \prod_{1 \leq a,b \leq k}  \left[ \frac {(\oq q^n)^{2(a-b)}-1}{q (\oq q^n)^{2(a-b)}-q^{-1}} \right]^n \cdot \prod_{l=0}^{k-1} \frac {q (\oq q^n)^{2l} - q^{-1}}{(\oq q^n)^{2l} - 1}
\end{multline*}
where $G_k$ is thought of having $k$ variables of each color $i,i+1,...,i+n-1$. Plugging in the explicit formula \eqref{eqn:g2}, we have:
$$
\beta_{[i;i+nk)}(G_{-k}) = \frac {(\oq q^n)^{k^2}}{(1-q^{-2})^{nk}} \prod_{l=1}^k \frac {q(\oq q^n)^{2-2l} - q^{-1}}{(\oq q^n)^{2l} - 1} \quad \stackrel{\eqref{eqn:bonnie2}}\Longrightarrow \quad \eqref{eqn:pairing 2}
$$

\end{proof}


\begin{thebibliography}{XXX}

\bibitem{BT} Bershtein M., Tsymbaliuk A., {\em Homomorphisms between different quantum toroidal and affine Yangian algebras}, J. Pure Appl. Algebra., Volume 223, Issue 2 (2019), 867--899

\bibitem{FD} Ding J., Frenkel I. {\em Isomorphism of two realizations of quantum affine algebra $U_q(\hgl_n)$}, Comm. Math. Phys. 156 (1993), no. 2, 277-300

\bibitem{Drin} Drinfeld V., {\em Quantum groups}, Proceedings of the ICM (1987), American Mathematical Society, Rhode Island

\bibitem{E} Enriquez B. {\em On correlation functions of Drinfeld currents and shuffle algebras}, Transform. Groups 5 (2000), no. 2, 111 - 120

\bibitem{FRT} Faddeev L., Reshetikhin N., Takhtajan L. {\em Quantization of Lie groups and Lie algebras, Yang-Baxter equation in Integrable Systems}, Adv. Ser. Math. Phys. 10 (1989), 299-309

\bibitem{F} Feigin B., Hashizume K., Hoshino A., Shiraishi J., Yanagida S. {\em A Commutative Algebra on Degenerate $\BC^1$ and MacDonald Polynomials},  J. Math. Phys. 50 (2009), no. 9

\bibitem{FJM} Feigin B., Jimbo M., Mukhin E., {\em Integrals of motion from quantum toroidal algebras}, Journal of Physics A: Mathematical and Theoretical, Volume 50, Number 46

\bibitem{FJMM} Feigin B., Jimbo M., Miwa T., Mukhin E. {\em Representations of quantum toroidal $\fgl_n$}, preprint arXiv:1204.5378

\bibitem{FO} Feigin B., Odesskii A. {\em Quantized moduli spaces of the bundles on the elliptic curve and their applications}, Integrable structures of exactly solvable two-dimensional models of quantum field theory (Kiev, 2000), 123-137, NATO Sci. Ser. II Math. Phys. Chem., 35, Kluwer Acad. Publ., Dordrecht, 2001

\bibitem{FT} Feigin B., Tsymbaliuk A. {\em Bethe subalgebras of quantum affine $U_q(\widehat{\fgl}_n)$ via shuffle algebras}, Sel. Math. New Ser., vol. 22 (2016), 979--1011

\bibitem{H} Hernandez D. {\em Quantum toroidal algebras and their representations}, Selecta Math. 14 (2009), no. 3-4, 701-725

\bibitem{GM} Gow L., Molev A. {\em Representations of twisted $q-$Yangians}, Selecta Math. 16 (2010), 439-499

\bibitem{Jan} Jantzen J. {\em Lectures on quantum groups}, Graduate Studies in Mathematics, 6. American Mathematical Society, Providence, RI, 1996. viii+266 pp. ISBN: 0-8218-0478-2

\bibitem{K} Kassel C. {\em Quantum groups}, Springer-Verlag, 1995

\bibitem{KR} Kirillov A., Reshetikhin N., {\em $q$--Weyl group and a multiplicative formula for universal $R$--matrices}, Comm. Math. Phys. 184 (1990), 421--431

\bibitem{KT} Khoroshkin S., Tolstoy V. {\em The universal $R-$matrix for quantum untwisted affine Lie algebras}, Functional Analysis and Its Applications, January-March, 1992, Vol 26, Issue 1, pp 69-71

\bibitem{L} Lusztig G., {\em Introduction to Quantum Groups}, Boston, Birkh\"auser, 1993

\bibitem{LS} Levendorsky S., Soibelman Ya., {\em Some applications of quantum Weyl groups}, J. Geom. Phys. 7 (1990), 241--254

\bibitem{LSS} Levendorsky S., Soibelman Ya., Stukopin V., {\em The Quantum Weyl group and the universal quantum R-Matrix for affine Lie algebra $A_1^{(1)}$}, Lett Math Phys 27, 253--264 (1993)

\bibitem{Aff} Negu\cb t A. {\em Affine Laumon spaces and a conjecture of Kuznetsov}, ar$\chi$iv:1811.01011

\bibitem{CM} Negu\cb t A. {\em Affine Laumon Spaces and Integrable Systems}, ar$\chi$iv:1112.1756

\bibitem{thesis} Negu\cb t A. {\em Quantum Algebras and Cyclic Quiver Varieties}, PhD thesis, Columbia University, 2015, ar$\chi$iv:math/1504.06525

\bibitem{Shuf} Negu\cb t A. {\em The shuffle algebra revisited}, Int Math Res Not, 2013, doi: 10.1093/imrn/rnt156

\bibitem{PBW} Negu\cb t A. {\em The PBW basis of $U_{q,\bar{q}}(\ddot{\fgl}_n)$}, ar$\chi$iv:arXiv:1905.06277	

\bibitem{R} Rosso M., {\em An analogue of P.B.W. theorem and the universal R-matrix for $U_h\fsl(N+1)$}, Commun. Math. Phys. vol. 124, 307--318 (1989)

\bibitem{RS} Reshetikhin N., Semenov-Tian-Shansky M. {\em Central Extensions of Quantum Current Groups}, Lett. Math. Phys. 19 (1990), 133 - 142

\bibitem{T} Tsymbaliuk A. {\em Quantum affine Gelfand-Tsetlin bases and quantum toroidal algebra via $K-$theory of affine Laumon spaces}, Selecta Math. (N.S.) 16 (2010), no. 2, 173-200

\end{thebibliography}
\end{document}